\documentclass[review]{elsarticle}

\usepackage{lineno,hyperref}
\modulolinenumbers[5]

\journal{Journal of \LaTeX\ Templates}
\usepackage{color}
\usepackage{subfig}
 \usepackage[utf8]{inputenc} 
 \usepackage{amsmath}
\usepackage{amssymb}
\usepackage{siunitx}
\usepackage{array}
\usepackage{lineno}
\usepackage{booktabs}
\usepackage{graphicx}

\usepackage[]{algorithm2e}
\usepackage{mathabx}

\usepackage{mathtools}

\newtheorem{mytheo}{Theorem}
\newtheorem{mydef}{Definition}

\usepackage{mathrsfs}

\newcommand{\mvecfun}[1]{ \underline{#1}}
\newcommand{\mpoint}[1]{ \mathbf{#1}}

\newcommand{\matent}[3]{ [\mathbf #1]_{#2#3}}

\newcommand{\domain}[0]{\Omega}
\newcommand{\fracset}[0]{\mathcal F}

\newcommand{\cardfracset}[0]{n^\fracset}
\newcommand{\cardall}[0]{n^H}
\newcommand{\cardreg}[0]{n^r}

\newcommand{\dirset}[0]  {\Gamma_D}
\newcommand{\neuset}[0]{\Gamma_N}
\newcommand{\cdirset}[0] {\overline\Gamma_D}
\newcommand{\cneuset}[0]{\overline \Gamma_N}
\newcommand{\idirset}[0]  {\mathring \Gamma_D}
\newcommand{\ineuset}[0]{\mathring\Gamma_N}

\newcommand{\mesh}[0]{\mathcal T}
\newcommand{\hh}[1]{ {#1}^h }

\newcommand{\mflux}[2]{q_{#2}(#1)}

\newcommand{\iset}[0]{J}
\newcommand{\isetdir}[0]{D}
\newcommand{\isetrest}[0]{I}
\newcommand{\isetall}[0]{J^H}
\newcommand{\isetreg}[0]{J^r}
\newcommand{\isethang}[0]{J^{this should not be used}}

\newcommand{\mb}[1]{ N_{#1} }
\newcommand{\mbconf}[1]{N^r_{#1} }
\newcommand{\mbnc}[1]{ N^H_{#1} }

\newcommand{\basissetconf}[0]{ B^r }
\newcommand{\basissetnc}[0]{ B^H }
\usepackage{diagbox}
\newcommand{\mS}[0]{S}
\newcommand{\mU}[0]{U}
\newcommand{\mV}[0]{V}

\newcommand{\mSnc}[0]{W^H}

\newcommand{\randomset}[0]{\Xi}

\newcommand{\lineb}{\text{B}\text{B}^{\prime}}
\newcommand{\linea}{\text{A}\text{A}^{\prime}}
\newcommand{\linec}{\text{C}\text{C}^{\prime}}
\newcommand{\lined}{\text{D}\text{D}^{\prime}}
\newcommand{\mG}[0]{G}

\newcommand{\mrestrictionoperator}[1]{\mathcal{R}_{#1}}
\newcommand{\mprolongationoperator}[1]{\mathcal{P}_{#1}}
\newcommand{\mrestrictionmatrix}[1]{\mathcal{R}_{#1}}
\newcommand{\mprolongationmatrix}[1]{\mathcal{P}_{#1}}

\newcommand{\mvertex}[1]{ \mathbf x_{#1} }

\begin{document}

\begin{frontmatter}

\title{A novel equi-dimensional finite element method for\\
flow and transport in fractured porous media satisfying\\
discrete maximum principle and conservation properties
}
%
%

\author[a,b]{Maria Giuseppina Chiara Nestola\corref{mycorrespondingauthor}}
\cortext[mycorrespondingauthor]{Corresponding author}
\ead{nestom@usi.ch}
\author[c]{Marco Favino}
\ead{marco.favino@unil.ch}

\address[a]{Institute of Computational Science, Center for Computational Medicine in Cardiology (CCMC), Universit\`a della Svizzera italiana, Via Giuseppe Buffi 13, 9600 
Lugano, Switzerland.}
\address[b]{Institute of Geochemistry and Petrology, ETH Zurich Clausiusstrasse 25, NW F72 8092 Zurich, Switzerland.}
\address[c]{Institute of Earth Science, University of Lausanne, G\'eopolis - CH-1015  Lausanne, Switzerland.}

\begin{abstract}
Numerical simulations of flow and transport in porous media usually rely on hybrid-dimensional models, i.e., the fracture is considered as objects of a lower dimension compared to the embedding matrix.
Such models are usually combined with non-conforming discretizations as they avoid the inherent difficulties associated with the generation of meshes that explicitly resolve fractures-matrix interfaces.
However, non-conforming discretizations demand a more complicated coupling of different sub-models and may require special care to ensure conservative fluxes.
We propose a novel approach for the simulation of flow and transport problems in fractured porous media based on an equi-dimensional representation of the fractures.
The major challenge for these types of representation is the creation of meshes which resolve the several complex interfaces between the fractures and the embedding matrix.
To overcome this difficulty, we employ a strategy based on adaptive mesh refinement (AMR).
The idea at the base of the proposed AMR is to start from an initially uniform coarse mesh and refine the elements which have non-empty overlaps with at least one of the fractures.
Iterating this process allows to create non-uniform non-conforming meshes, which do not resolve the interfaces but can approximate them with arbitrary accuracy.
We demonstrate that low-order finite element (FE) discretizations on adapted meshes are globally and locally conservative and
we suitably adapt an algebraic flux correction technique to ensure the discrete maximum principle.
In particular, we show that the notorious conditions on M-matrices have to be adapted to the basis functions defined on non-conforming meshes.
Although the proposed applications come from geophysical applications, the obtained results could be applied to any diffusion and transport problems, on both conforming and non-conforming meshes.
\end{abstract}



\end{frontmatter}

\section{Introduction}

Numerical simulations of flow and transport problems in fractured porous media are fundamental for numerous geophysical applications,
such as geothermal energy production, hydrocarbon exploration, nuclear waste disposal, and CO$_2$ storage~\cite{tester2006,mcclure2014b,bond2003,bonnet2001,rasmuson1986,amann2018}.
Fractures can be regarded as heterogeneities inside an embedding matrix and can create preferential paths and/or barriers for flow and transport.
Fractures can be arranged in complicated networks and their aperture is several orders of magnitude smaller than the characteristic size of the embedding matrix.
For these reasons, models that represent the fractures as inclusions with the same geometrical dimension as the embedding matrix (referred to as equi-dimensional models) are usually considered unfeasible~\cite{berre2019flow}.
In particular, the creation of meshes which explicitly resolve the complex interfaces between fractures and the embedding matrix is one of the main difficulties
and represents one of the major bottlenecks in the use of equi-dimensional models.

For overcoming such difficulties, alternative approaches based on implicit and/or lower-dimensional representations of the fractures have been derived.
Implicit representations consider the effect of the fractures as material parameters in an effective continuum.
This class includes single-continuum approaches~\cite{eikemo2009discontinuous,liu2016mathematical} and multi-continuum approaches where governing equations are defined for each \text{block} and are coupled by interaction terms~\cite{chung2017coupling}.
Lower-dimensional representations give rise to hybrid-dimensional approaches which model fractures as objects of a dimension lower than the one of the embedding matrix and the fracture aperture is treated as a parameter.
The most common approach currently employed is based on the so-called discrete fracture-matrix (DFM) models,
which implicitly treats the smaller fractures and the larger fractures are usually considered in a hybrid-dimensional approach.

The numerical simulation of both equi- and hybrid-dimensional approaches may make use of either conforming or non-conforming meshes.
The generation of conforming meshes is particularly complicated and difficult to make automatic for realistic fracture networks~\cite{karimi2003efficient,sandve2012efficient,ahmed2015control,angot2009asymptotic,fumagalli2019conforming,antonietti2019mixed}.

On the other hand, non-conforming meshes have to be employed in combination with specific discretization strategies in order to couple the equations defined on the matrix domain and on the fracture domain {such as the unfitted finite element methods}~\cite{koeppel2019,SCHADLE201942,frih2012,boon2018,burman2019simple,li2008efficient,odsaeter2019simple,formaggia2014reduced,schwenck2015dimensionally,flemisch2016review,burman2020cut,burman2019cut,antonietti2020polytopic}.

For what concerns flow problems in fractured porous media,
they have been approximated by
finite volume (FV)~\cite{karimi2003efficient,sandve2012efficient,ahmed2015control,angot2009asymptotic},
finite differences (FD)~\cite{antonietti2016mimetic,formaggia2018analysis},
virtual finite element (VFE)~\cite{fumagalli2019dual},
and finite element (FE) methods~\cite{boon2018robust,formaggia2014reduced,larson2004conservative}.
FV method is well-known to ensure local conservation properties by construction,
i.e., ensure the balance of fluxes between two arbitrary, complementary subdomains.
Conservation properties for FE method hold naturally for mixed formulations based on Raviart-Thomas discretizations~\cite{boon2018robust,martin2005modeling,d2012mixed,chave2018hybrid,douglas1982mixed,frih2012,boon2018},
or discontinuous Galerkin (DG) methods~\cite{rivie2000part},
but they can also be proven for enriched continuous Galerkin methods~\cite{sun2009locally,larson2004conservative},
the extended finite elements (XFEM)~\cite{formaggia2014reduced,schwenck2015dimensionally,flemisch2016review},
and the embedded finite element methods (EFEM)~\cite{odsaeter2019simple}.
Moreover, the virtual elements method (VEM) and mimetic finite difference (MFD) have been adapted~\cite{antonietti2016mimetic,fumagalli2019dual} to preserve conservation properties
even in the presence of a highly anisotropic and distorted grid.
For FE methods based on continuous Galerkin (CG) formulations,
the computation of the approximated fluxes relies, in general, on post-processing techniques of the numerical solution~\cite{wheeler1973simulation,carey2002some,cockburn2007locally,odsaeter2017postprocessing}
and the resulting fluxes are as accurate as the ones computed employing discretizations of the same order based on mixed formulations~\cite{cockburn2007locally}.
In~\cite{hughes}, the authors demonstrated that low-order FE discretizations are conservative by employing an equivalent formulation of the problem where an auxiliary flux is associated with the Dirichlet boundary conditions has been introduced.
This alternative formulation allowed to prove that continuous FE formulations are conservative for any patch of elements in case of conforming meshes.

The transport of a concentration through a fractured porous medium
is described by a pure advection problem, whose velocity field is computed through the solution of a flow problem.
The solution of a transport problem can be proven to respect the maximum principle and, hence, to be positive for positive boundary and initial conditions.
On the other hand, a violation of the discrete maximum principle (DMP) for the solutions of discretizations of transport equations
manifests itself in spurious oscillations and negative values.
Satisfying the DMP is, hence, a key requirement when dealing with transport problems and, more in general, with differential problems involving positive physical quantities,
such as concentrations.
Moreover, discretizations whose solution respect the DMP can be proven to converge uniformly~\cite{ciarlet1973maximum}.

Discretizations that provide positive solutions are usually characterized by stiffness matrices which are monotone or M-matrices.
{In particular low-order FE discretizations are characterized by matrices that have positive diagonal entries and non-negative extra-diagonal entries.}
A typical example of matrices that respect these properties come from the discretization of Laplace operators on grids whose elements do not present obtuse angles.
The FE discretization based on CG formulations of transport problems typically gives rise to stiffness matrices which violate the DMP.
This is due to the presence of positive extra-diagonal entries in the mass matrix and in the discrete advection operator and of negative diagonal entries in the discrete advection operator.

Several approaches have been introduced to try to ensure the DMP and
preserve the positivity of the solution such as
1) artificial diffusion~\cite{fischer1999parameter},
2) stream-line diffusion ~\cite{hughes1979multidimentional},
3) residual-free bubbles~\cite{brezzi1998further},
4) local projection stabilizations~\cite{john2008finite}. 
The idea at the base of all those methods is to introduce an additional diffusion term to the formulation of the transport problem.
On the other hand, algebraic flux correction~\cite{kuzmin2009explicit} introduces an algebraic diffusion operator to compensate for the entries which { does not allow to ensure the DMP}.
Instead, FV discretizations naturally provide monotone matrices when employed with upwind schemes.

Most of the discretization methods for flow and transport problems are based on hybrid-dimensional representations,
while few are based on equi-dimensional representations.
Equi-dimensional representations of fractures are less popular mainly due to the complexity of meshing complicated fractures networks~\cite{antonietti2016mimetic}.
The creation of meshes that resolve the interfaces between fractures and the embedding matrix is a time-consuming and difficult-to-make-automatic process.
Moreover, the resulting meshes are characterized by elongated and/or distorted elements, due to the small aperture of the fractures.
On the other hand, hybrid dimensional approaches require the development of
1) specific strategies to couple the equations defined on the fracture domain and matrix domain to impose the continuity of the unknowns at the interface;
2) specific algorithms to detect the intersection between the elements of the fracture and the matrix domain;
3) solutions strategies designed to treat the coupling terms~\cite{formaggia2014reduced,flemisch2016review,odsaeter2019simple,koppel2019stabilized}.
Moreover, hybrid-dimensional models
neglect cross fracture phenomena,
do not always guarantee local conservation at the matrix-fracture interfaces~\cite{xu2020hybrid},
may not ensure the DMP.

The focus of this work is to propose a discretization method for the simulation of the coupled flow and transport problems in fractured media
based on an equi-dimensional representation of the fractures. 
For creating suitable meshes, we employ the adaptive mesh refinement (AMR) strategy proposed in~\cite{Holliger2000}.
The idea of this strategy is to start from a uniform coarse mesh that is not related to the fracture distribution and to refine
the elements which have a non-empty overlap with at least one fracture.
Iterating this step allows generating adapted non-conforming meshes which are refined at the interfaces between the fractures and the embedding matrix.
The resulting meshes are composed of squared and cubic elements which reduce the ill-conditioning due to the presence of skewed elements,
at the price of having non-constant material properties over some elements.
This approach allows to explicitly take fractures into account by assigning different material properties during the assembly of the stiffness matrices.

For transport problems with shock waves,
linear and non-linear algebraic stabilizations
have been compared on adapted meshes to preserve monotonicity on adapted meshes~\cite{bonilla2020monotonicity}.
Differently, from our geometric approach, the authors employed a specific error indicator to select the elements to refine.

We show that the proposed discretization is conservative on non-conforming adapted meshes and for any subdomain.
In combination with an algebraic flux correction scheme,
such discretization ensures the DMP.
Moreover, the CG formulation does not present the matrix-fracture coupling terms
and, hence, produces linear systems with a positive definite matrix.

The paper is organized as follows.
In Section 2, we discuss the DMP on conforming meshes and introduce the governing equations for the coupled flow and transport problems.
In Section 3, we introduce the FE formulation and discuss the DMP on adapted (non-conforming) meshes while in Section 4 we study the conservation properties of the proposed approach.
In Section 5, we present some numerical results, including a real problem with a complex stochastic fracture network.
Finally, we make some concluding remarks in Section 6.

\section{Discrete maximum principle on conforming meshes}

We introduce the FE discretizations of the coupled models which describe the single-phase flow and the advective transport of a solute in a fractured porous medium.
and their relevant properties in terms of continuum and discrete maximum principle.
We employ models based on an \emph{equi-dimensional} representation of the fractures \cite{flemisch2018benchmarks,berre2020verification},
i.e., fractures are described as $d$-dimensional objects embedded in a $d$-dimensional porous matrix.
Equi-dimensional representations do not involve couplings between the discretizations of different dimensional objects.
Hence, they allow for a natural extension of classical results related to FE methods based on CG formulations.

\subsection{Preliminaries definitions}
\label{sub:matpre}

We denote an axis-aligned box domain in $\mathbb R^d$,
i.e., a rectangle for $d=2$ or a cuboid for $d=3$, by $\domain$,
where $d$ denotes the dimension of the problem.
The length along the $i$-axis is denoted by $L_i$.
We denote by $\partial \domain$ the boundary of $\domain$.
The outward unit vector to $\domain$ is denoted by $\mvecfun n$.
Points of $\domain$ are denoted by capital non-bold letters $X$ and their coordinates by $\mpoint {x}=(x_1,x_2, \ldots, x_d)$.

Over $\domain$, we introduce a set $\fracset$ of $d$-dimensional inclusions $f_i \subset \domain$ with $i=1,2\,\ldots \cardfracset$,
where $\cardfracset= | \fracset |$ is the number of inclusions.
We assume inclusions are non-axis-aligned boxes.
We are particularly interested in the case where inclusions have one dimension that is much smaller compared to the other(s).
In this case, we will refer to such inclusions as fractures and the smaller dimension will be referred to as aperture or thickness.

We define $\domain_f:=\bigcup_{i=1}^{\cardfracset}f_i$, i.e., the sets of the inclusions and $\domain_m := \domain \backslash \domain_f$, i.e., the set representing the matrix.
Both $\Omega_m$ and $\Omega_f$ may be non-connected sets.
The set $\Gamma:= \overline \Omega_m \cap \overline \Omega_f$ is the interface between the two subdomains.

We denote the fracture aperture by $\delta$, the permeability by $k$, and the porosity by $\phi$.
Permeability and porosity are assumed to be real functions attaining constant positive values over $\Omega_m$ and $\Omega_f$,
such as, for example,
\begin{equation}
k=\begin{cases} k_m, \quad \text{ in } \domain_m,\\
k_f, \,\,\quad \text{ in } \domain_f.\\
\end{cases}
\end{equation}

It is worth to point out that the assumptions adopted to introduce the reader to the topic could be relaxed:
$\domain$ has not to be a box domain;
fractures may have variable apertures and different shapes; the material properties can be heterogeneous. In particular, heterogeneous properties would require to replace scalar values by tensor quantities.

We suppose that $\partial \domain$ can be decomposed in two subsets $\dirset$ and $\neuset$
with $\partial \Omega=\cdirset \cup \cneuset $ and $\idirset \cap \ineuset = \emptyset$.
The set $\dirset$ denotes the subset of $\partial \domain$ where Dirichlet boundary conditions are imposed.
Notice that $\dirset$ may differ for flow and transport problems.

For a set $\randomset$, we denote by the set of square-integrable functions by $L^2(\randomset)$.
The scalar product and the norm in $L^2(\randomset)$ are $\| \, \cdot \, \|_{\randomset}$, and $(\,\cdot\,,\, \cdot \, )_{\randomset}$, respectively.
Given an unknown $q$ and Dirichlet boundary conditions $q_D$ imposed on $\dirset$, the set $\mS$ denotes a suitable function space where the solution of each problem is sought.
{The space $S$ will be specified in Section~\ref{sec:DMPAMR} for each of the problems}. 
Moreover, we introduce the following function spaces:
\begin{equation}
\mU=\{ v \in \mS \, : \, v=q_D \text{ on } \dirset\},
\end{equation}
and 
\begin{equation}
\mV=\{ v \in \mS \, : \, v=0 \text{ on } \dirset\}.
\end{equation}

Over $\domain$, we introduce a mesh $\mesh$ with $N_E$ elements and $N_N$ nodes.
Elements can be triangles or quadrilaterals for $d=2$ and tetrahedra or hexahedra for $d=3$.
We assume that all elements of a mesh are of the same type.
Elements are denoted by $E$ and
have vertices, edges, and, for $d=3$, faces.
The term side refers to edges for $d=2$ and to faces for $d=3$.
We assume, for now, that the mesh is conforming, i.e., 
the intersection of each pair of elements is empty, a vertex, an edge, or a face.

We let $\hh{\mS}$ denote the nodal interpolation space over $\mesh$, i.e.,
$$\hh{\mS}= \{ v \in C^0 (\Omega) \, : \, v |_E \in \mathbb L_1 \},$$
where $\mathbb L_1$ is the space of linear functions on $E$ for triangular and tetrahedral meshes (usually denoted by $\mathbb P_1$)
or the space of multi-linear functions on $E$ for quadrilateral or hexahedral meshes (usually denoted by $\mathbb Q_1$).

We let $\iset$ denote the set of all nodal indices $i = 1, 2, \ldots ,N_N$ and
we let $\isetdir$ be the subset corresponding to nodes located in $\dirset$, i.e.,
$$\isetdir = \{ i \, : \, \mpoint{x}_i \in \dirset \},$$
and $\isetrest = \iset \setminus \isetdir$.
We call $\mb{i}$ the Lagrangian basis function associated with the node $\mpoint{x}_i$.
Observe that the space $\hh{\mS}= \mathrm{span}\{ \mb{i}\}_{i\in \iset}$.

We denote the set of functions in $\hh{\mS}$ which satisfy the Dirichlet boundary condition by $\hh{\mU}$ and
the set of functions $\hh{\mS}$ which attain zero value on the Dirichlet boundary by $\hh{\mV}$.
Finally, we introduce $\hh{\mG} = \mathrm{span}\{\mb{i}\}_{i \in \isetdir}$, i.e.,
the span of the basis functions associated with the set of nodes in $\isetdir$.
We observe that the space $\hh{\mV}$ is the span of the basis functions associated with each node except the ones where Dirichlet conditions are imposed. Moreover, the following relations hold:
\begin{equation}
\label{eq:discretespacedefinition}
\begin{array}{lcll l}
\hh{\mU} &= & \mU \cap \hh{\mS},\\
\hh{\mV} &= & \mV \cap \hh{\mS} &= & \mathrm{span}\{ \mb{A}\}_{A\in \isetrest},\\
\hh{\mS} &= & \hh{\mV} \bigoplus \hh{\mG} &= & \mathrm{span}\{ \mb{A}\}_{A\in \iset},\\
\hh{\mU} &\subset& \hh{\mS}.
\end{array}
\end{equation}

\subsection{The flow problem}

For an incompressible single-phase fluid, the flow through a porous medium is described by
\begin{equation}
\label{eq:flowstrong}
\begin{array}{c c c l}
\nabla \cdot \mvecfun{u} &= & 0 & \quad \text{in } \domain,\\
\mvecfun{u} &= & -k \nabla p & \quad \text{in } \domain,\\
 p & = & g & \quad \text{on} \,\, \dirset,\\
 k \nabla p \cdot \mvecfun n& = & h & \quad \text{on} \,\, \neuset,
\end{array}
\end{equation}
where $\mvecfun{u}$ is the velocity and $p$ the pressure.
For the flow problem, the space $S$ coincides with the Sobolev space $H^1(\domain)$.
The weak formulation of problem \eqref{eq:flowstrong} reads:
\begin{equation}
\label{eq:flowweak}
\begin{array}{l}
\text{Find } p\in U \, \text{ such that }\\[2mm]
{\displaystyle d(p,v) = f(v) \quad \forall v \in V,}
\end{array}
\end{equation}
where
$$d(p,v)= \int_\domain k\, \nabla p \cdot \nabla q \,\text{d}V \quad \text{and} \quad f(v)= \int_{\neuset} h \, v \,\text{d}A.$$

\begin{mytheo}
The following a-priori estimates hold for problem~\eqref{eq:flowweak}:
\begin{align}
\nonumber
\min_{\partial \domain} p < p<\max_{\partial \domain}  p & \qquad \quad (\text{Maximum Principle}),\\[1mm] 
\nonumber
p|_{\partial \domain}\geq0 \implies p\geq 0 & \qquad \quad \text{(Positivity preservation)}.
\end{align}
\end{mytheo}

With the definitions introduced in Subsection \ref{sub:matpre},
the FE approximation to problem~\eqref{eq:flowstrong} reads
\begin{equation}
\label{eq:flowweakdiscrete}
\begin{array}{l}
\text{Find } \hh{p} \in \hh{U} \, \text{ such that }\\[2mm]
{\displaystyle d( \hh{p}, \hh{v}) = f( \hh{v}) \quad \forall \hh{v} \in \hh{V},}
\end{array}
\end{equation}
where
$$ d( \hh{p}, \hh{v}) = \int_\domain k\, \nabla p_h \cdot \nabla v_h \,\text{d}V \quad \text{and} \quad f(v)= \int_{\neuset} h \, v _h \,\text{d}A.$$

Hence, the flow problem admits the following algebraic representation
\begin{equation}
\label{eq:flowdiscrete_0}
\mathbf D\\
\mathbf p\\
=
\mathbf f,
\end{equation}
or, explicitly,

\begin{equation}
\label{eq:flowdiscrete_2}
\begin{pmatrix}
\mathbf D_{II} & \mathbf D_{ID} \\
\mathbf 0 & \mathbf I
\end{pmatrix}
\begin{pmatrix}
\mathbf p_{A} \\
\mathbf p_D 
\end{pmatrix}
=
\begin{pmatrix}
\mathbf f_{I} \\
\mathbf g_D 
\end{pmatrix},
\end{equation}
where
\begin{equation*}
\begin{array}{l c l c l}
\matent{D_{II}}{i}{j} & = &d(\mb{j},\mb{i}) & \quad & i,j\in \isetrest,\\
\matent{D_{ID}}{i}{j} & =& d(\mb{j},\mb{i}) & \quad & i\in \isetrest, \, j \in \isetdir,\\
\matent{f_{I}}{i}{} & =& f(\mb{i}) & \quad & i\in \isetrest.\\
\end{array}
\end{equation*}

\subsection{The transport problem}

The transport problem is coupled to the flow problem through the velocity $\mvecfun u$ as it enters both in the definition of the boundary conditions
and in the transport equation.
We let $\Gamma_{in} := \{ \mathbf x \in \partial \Omega \, : \, \mvecfun u \cdot \mvecfun n < 0\}$ denote the inflow boundary
and $\Gamma_{out} := \{ \mathbf x \in \partial \Omega \, : \, \mvecfun u \cdot \mvecfun n \geq 0\}$ denote the outflow boundary.
The unknown is the concentration $c=c(\mathbf x, t)$ of a given solute.
As we consider only Dirichlet problems, we have $\Gamma_{in} = \Gamma_{D}$.
The equi-dimensional model for the advective transport problem reads:

\begin{equation}
\label{eq:flowclassic}
\left\{
\begin{array}{r c l l}
{\displaystyle \phi \frac{\partial c}{\partial t} + \nabla \cdot ( \mvecfun u \, c ) }& = &0 &\quad \text{in}\,\, \Omega \times \mathscr I\\
c( \,\cdot \, , 0) &= & c_0 & \quad \text{in}\,\, \Omega\\
c & = & g & \quad \text{on}\,\, \Gamma_{in} \times \mathscr I
\end{array}
\right.
\end{equation}
Here, $\mathscr I=(0, T_{\text{fin}}]$ is the time interval, $c_0$ is the initial condition, and
$g$ is the inflow Dirichlet boundary condition.

Defining by $\Sigma$ the set of points where initial and boundary conditions are prescribed, i.e., 
$$\Sigma:=\{(\mathbf x,t) \, : \, \mathbf x \in \Gamma_{in} \lor t = 0\},$$
the following theorem holds.
\begin{mytheo}
\label{theo:this}
For problem \eqref{eq:flowclassic} the following a-priori estimates hold:
\begin{equation*}
\begin{array}{r c l l}
\nabla \cdot \mvecfun u=0&\implies & \min_{\Sigma}  c \leq c\leq \max_{ \Sigma} c& \quad(\text{Maximum Principle}), \\[1mm]
c|_{\Sigma}\geq0&\implies & c\geq 0 & \quad\text{(Positivity preservation)}.
\end{array}
\end{equation*}
\end{mytheo}

Theorem \ref{theo:this} ensures that if $c_0$ and $g$ are positive,
then the solution is positive all over the space-time domain $\Omega \times \mathscr I $.
Moreover, the concentration in the domain cannot be larger than the maximum value injected at the inflow boundary.

\subsubsection*{Finite element discretization of transport equation}

The function $c_h=c_h(\, \cdot \, , t) \in \hh{U}$ is the approximation of the concentration $c$ at time $t$
and $\mathbf c(t)$ is the time dependent array having as components the unknown coefficients $c_j$ with respect to the basis $\{ N_j \}$.
The spatial discretization of \eqref{eq:flowclassic} reads:
\begin{equation}
\label{eq:transportweakdiscrete}
\begin{array}{l}
\text{For all } t \in \mathscr I, \text{ find } c_h(\, \cdot \, , t) \in U_h \, \text{ such that } c( \,\cdot \, , 0) = c_0,\,\, \text{and} \\[2mm]
{\displaystyle m \left( \dfrac{ \partial c_h}{\partial t} , q_h \right ) + a (c_h , q_h)= 0 \quad \forall q_h \in V_h.}
\end{array}
\end{equation}
where we have used the following discrete bilinear forms:
\begin{equation}
\nonumber
\begin{array}{ l c l}
 m \left( \dfrac{ \partial c_h}{\partial t} , q_h \right ) &=& {\displaystyle \hphantom{-} \int_\Omega \phi \, \dfrac{ \partial c_h}{\partial t} \, q_h \, \text{d}V,}\\[3mm]
 \nonumber
 a (c_h , q_h) & = & {\displaystyle - \int_{\Omega} c_h\, \mvecfun u\cdot \nabla q_h \,\text{d}V + \int_{\Gamma_{out}} c_h \, q_h \, \mvecfun u \cdot \mvecfun n \,\text{d}A\, .}
\end{array}
\end{equation}
Problem \eqref{eq:transportweakdiscrete} admits the following algebraic representation
\begin{equation}
\label{eq:transport_discrete}
\mathbf M \frac{ \text{d} \mathbf c }{ \text{d} t} + \mathbf A \mathbf c = \mathbf 0,
\end{equation}
{where the matrices $\mathbf M$ and $\mathbf A$ are the scaled mass matrix and the discrete advection operator, respectively, and will be specified later.}

A classic theorem that provides sufficient conditions to ensure the maximum principle and positivity preservation for the semi-discrete transport problem~\eqref{eq:transport_discrete} is the following{~\cite{Kuzmin}}.
\begin{mytheo}
\label{th:DM}
\text{Suppose that}:
\begin{itemize}
\item $\matent{M}{i}{i} > 0$,
\item $\matent{M}{i}{j} = 0$ for $i\neq j$,
\item $\matent{A}{i}{j} \leq 0 $ for $i\neq j$.
\end{itemize}
Then the following a priori estimates hold for the coefficient $c_i$:
\begin{itemize}
\item[1a] The {semi-Discrete Maximum Principle} (DMP)  is satisfied:

$\quad\quad{\displaystyle \sum_{j}\matent{A}{i}{j}=0, \quad c_j\geq c_i, \quad \forall j \neq i \quad \implies \quad \dfrac{d c_i}{dt} \leq 0.}$
\item [1b] The {positivity preservation} is satisfied: $$c_j(0) \geq 0,\, \forall j\,\, \,\,\implies c_i(t)\geq0, \forall t>0.$$
\end{itemize}
\end{mytheo}

Employing an implicit Euler scheme,
one needs to solve a linear system at each time step.
The fully discrete counterpart of problem \eqref{eq:transport_discrete} is a sparse linear system of the form
\begin{equation}
\label{eq:transport_discrete_full}
\begin{pmatrix}
\mathbf B_{II} & \mathbf B_{ID} \\
\mathbf 0 & \mathbf I
\end{pmatrix}
\begin{pmatrix}
\mathbf c_{I}^{n+1} \\
\mathbf c_D^{n+1} 
\end{pmatrix}
=
\begin{pmatrix}
\mathbf M_{II} & \mathbf M_{ID} \\
\mathbf 0 & \mathbf 0
\end{pmatrix}
\begin{pmatrix}
\mathbf c_{I}^{n} \\
\mathbf c_D^{n} 
\end{pmatrix}
+
\begin{pmatrix}
\mathbf 0_{I} \\
\mathbf g_D 
\end{pmatrix},
\end{equation}
where
$$\matent{B}{i}{j} = \matent{M}{i}{j} + \Delta t \matent{A}{i}{j}.$$
with $\Delta t$ being the time-step size.
We point out that the backward Euler scheme is unconditionally stable with no restriction on the time-step size~\cite{nmdp}.

\subsection{Discrete maximum principle}

For FE approximations, a maximum principle for the discrete solutions, in general, does not descend from the properties of the continuous counterparts \eqref{eq:flowstrong} and \eqref{eq:flowclassic}.
In this section, we focus on the global \emph{discrete maximum principle}.
For a detailed discussion on local definition, we refer to~\cite{kuzmin2009explicit,Kuzmin}.

\begin{mydef} The solutions to \eqref{eq:flowdiscrete_2} or \eqref{eq:transport_discrete_full} satisfy the discrete maximum principle if
\begin{equation}
\label{eq:dmp}
\min_{j} g_j \leq u_i \leq \max_j g_j.
\end{equation}
\end{mydef}

\begin{mydef} The solutions to \eqref{eq:flowdiscrete_2} or \eqref{eq:transport_discrete_full} are said to be globally positivity-preserving if
\begin{equation}
\label{eq:pp}
\mathbf g \geq 0 \implies \mathbf u \geq 0.
\end{equation}
\end{mydef}
Typical proofs of \eqref{eq:dmp} and \eqref{eq:pp} are based on the theory of monotone matrices and in particular on M-matrices.

\begin{mydef} A regular matrix $\mathbf A$ is said to be monotone if $\mathbf A^{-1}\geq \mathbf 0$ or, equivalently, if
\label{def:3}
\begin{equation}
\mathbf A \mathbf x \geq \mathbf 0 \implies \mathbf x \geq \mathbf 0.
\end{equation}
\end{mydef}

\begin{mydef}
\label{def:4}
A monotone matrix $\mathbf A$ with $\matent{A}{i}{j} \leq 0$ for $i \not = j$ is called an M-matrix.
\end{mydef}
In general, the conditions expressed by Definitions~\ref{def:3} and~\ref{def:4} are quite complicated to verify.
The following theorem provides some mild conditions for a matrix to be an M-matrix~\cite{varga1966discrete}.
\begin{mytheo}
\label{theo:3}
A matrix $\mathbf A$ that
\begin{enumerate}
\item is diagonally dominant by rows and strictly diagonally dominant for at least one row,
\item satisfies $\matent{A}{i}{j}\leq 0$ for $i \not = j$,
\item satisfies $\matent{A}{i}{i} >0 $,
\end{enumerate}
is an M-matrix.
\end{mytheo}

\subsubsection{Discrete Diffusion Operators}

The violation of the DMP conditions in Theorem~\ref{theo:3} may be caused by the presence of positive extra-diagonal entries or negative diagonal entries.
{Negative diagonal entries can be found only in the discrete advection operator,
while positive extra-diagonal entries can be found both in the  the discrete advection operator and in the mass-matrix}.
A strategy employed to remove such entries is based on adding a suitable artificial diffusion operator $\mathbf S${~\cite{Kuzmin}}.

\begin{mydef}
\label{def:5}
A symmetric matrix $\mathbf S$ is called a discrete diffusion operator if 
\begin{itemize}
\item $\matent{S}{i}{j}=\matent{S}{j}{i} \leq 0$ for $i \not = j$,
\item $\sum_i \matent{S}{i}{j} = 0$.
\end{itemize}
\end{mydef}
Observe that for a discrete diffusion operator the diagonal entries are positive and the column sum is negative.

{Given a matrix $\mathbf Q$, we call $\mathbf S^Q$ the discrete diffusion operator  for which:}

\begin{itemize}
\item $\matent{S^Q}{i}{j}=-\,\max(0, \matent{Q}{i}{j}, \matent{Q}{j}{i}) \quad \text{for } i \neq j$,
\item $\matent{S^Q}{i}{i}=-\sum_{i \neq j} \matent{S^Q}{i}{j}$.
\end{itemize}

The matrix $\tilde {\mathbf Q} = \mathbf Q + \mathbf S^Q$ is by construction an M-matrix.

For a matrix $\mathbf Q$ arising from the FE discretization of a bilinear form,
the construction of $\mathbf S^Q$ can be performed during the assembly of the matrix $\mathbf Q$.

\subsection{Discrete maximum principle for flow problem}
For FE discretizations, trying to satisfy the conditions of Theorem~\ref{theo:3} leads, in general,
to geometric restrictions to the elements' shapes.
Actually, the imposition of such conditions for any element-wise contribution to the stiffness matrix can be
translated into specific constraints on the angles comprised between two different sides of the same element and on the ratio between the shortest and the largest lengths (i.e., the aspect ratio).
In particular, for two-dimensional problems, only angles smaller then $\pi/2$ are allowed in a mesh.
For quadrilateral meshes with rectangular elements, the elements have to be of a non-narrow type, i.e., with aspect ratio smaller than $\sqrt{2}$. 

If a mesh does not satisfy those geometrical constraints, one may replace the operator ${\mathbf D}$ with a corresponding stabilized counterpart $\tilde {\mathbf D} = \mathbf D + \mathbf S^D$.
The addition of an algebraic diffusion operator provides a stiffness matrix that satisfies the hypothesis of Theorem~\ref{theo:3}.

\subsection{Discrete maximum principle for the transport problem}

Continuous FE discretizations applied to transport problems are characterized by spurious undershoots and overshoots in the numerical solution
as the matrices arising do not satisfy the hypothesis of Theorem~\ref{theo:3}.
Several \emph{stabilization} strategies have been proposed
to make  $\mathbf M$ and $\mathbf A$ consistent with the hypotheses of the theorem.
In the next, we refer to the algebraic flux correction (AFC) scheme proposed by \cite{Kuzmin}.
The basic idea of this method is to employ an algebraic diffusion operator that allows ensuring the DMP and preserve of positivity of the solution.


The AFC is a predictor-corrector strategy which consists of two steps:
\begin{enumerate}
\item advance the solution in time by a low-order scheme that incorporates enough numerical diffusion to suppress undershoots
and overshoots;
\item add a correction to achieve higher accuracy without violating the DMP.
\end{enumerate}
In order to satisfy the hypotheses of Theorem~\ref{th:DM}, one needs to $1)$ approximate the consistent mass matrix $\mathbf M$ by its lumped counterpart $\mathbf{M}_L$ and
$2)$ eliminate the off-diagonal positive entries of the transport operator $\mathbf{A}$ by adding a discrete diffusion operator, i.e., replacing $\mathbf A$ with $\tilde {\mathbf A}= \mathbf{A} + \mathbf{S}^A$.

In the process of AFC, the discrete problem~\eqref{eq:transport_discrete_full} can be split into a \textit{good} diffusive part of the form
$$\mathbf{M}_L \mathbf{c}^{n+1} +\Delta t \,\tilde {\mathbf A}\,\mathbf{c}^{n+1}=\mathbf M_L \mathbf c^{n},$$
and a \textit{bad} anti-diffusive part given by
$$\mathbf{f}=(\mathbf{M}_L-\mathbf{M})\mathbf{c}^{n+1} + \Delta t \,\mathbf{S}^A\mathbf{c}^{n+1}- (\mathbf{M}_L - \mathbf{M})\mathbf{c}^{n}.$$

The \textit{good} diffusive counterpart is positivity-preserving but overly diffusive due to its linearity. Hence, an anti-diffusive correction is required to recover the accuracy of the original high-order discretization.

The bad anti-diffusive term admits the decomposition $[\,\mathbf{f}\,]_i=\sum_{j\neq i}\matent{F}{i}{j}$ where the components of the sum represent a numerical flux which attain a local conservation property,~i.e. $\matent{F}{i}{j}+\matent{F}{j}{i}=0$ and are defined as follows:
$$ \matent{F}{i}{j}=\matent{M}{i}{j}(\mathbf{c}^{n+1}_{i}-\mathbf{c}_{j}^{n+1})- \matent{M}{i}{j}(\mathbf{c}^{n}_{i}-\mathbf{c}_{j}^{n})-\Delta t \,\matent{S}{i}{j}^D(\mathbf{c}^{n+1}_{i}-\mathbf{c}_{j}^{n+1}).$$

The anti-diffusive flux contributions $\matent{F}{i}{j}$ are limited in such a way that overall equations becomes less diffusive
but spurious oscillations are still suppressed. To this aim, the terms $\matent{F}{i}{j}$ are multiplied by a solution-dependent correction factor
 $0\leq\alpha\leq 1$ and then added to the \textit{good} diffusive part of the problem:
$$\mathbf{M}_L \mathbf{c}^{n+1} + \Delta t\, \tilde{\mathbf{A}} \,\mathbf{c}^{n+1}=\mathbf M_L \mathbf c^{n} + \tilde{\mathbf{f}},$$
where $[\,\tilde{\mathbf{f}}\,]_i= \sum_{j\neq i} \matent{F}{i}{j}\matent{{\alpha}}{i}{j}$ is the modified correction vector.
If all the factors $\matent{\alpha}{i}{j}$ are equal to 1 the original high-order Galerkin discretization is recovered, while the lower order discretization is attainied for $\matent{\alpha}{i}{j}=0$. 
{A detailed description of the technique used to compute the coefficients  $\matent{\alpha}{i}{j}$ can be found in~\cite{Kuzmin}}.

\section{Discrete maximum principle\\ on non-conforming meshes}
\label{sec:DMPAMR}

In this section, we present the derivation of FE discretizations which satisfy the DMP on non-conforming meshes.
As the DMP properties of the matrices and the discrete diffusion operators depend on the basis employed,
we show how basis functions on non-conforming meshes are constructed from the shape functions usually employed in the assembly of FE matrices.
Finally, we present a detailed study of the assembly of two matrices.
For a diffusion problem, such as the flow problem,
we show that on non-conforming meshes with square elements,
the conditions for the DMP may not be satisfied and its discretization may require stabilization with a discrete diffusion operator to ensure the DMP.
For a transport problem, we show that a discrete diffusion operator constructed employing the elemental shape functions, which may be intuitively chosen in the assembly process,
and the one constructed employing the FE basis functions provide different stabilization matrices and that only the second one ensures the DMP.

\subsection{Conforming spaces on non-conforming meshes}

We now let the possibility of $\mesh$ be a non-conforming mesh, i.e.,
a mesh where a vertex of an element may belong to the interior of an edge or a face of another element.
Such a vertex is said to be hanging.
Non-hanging vertices are instead said to be regular.
We restrict ourselves to cases of $1$-irregular meshes,
i.e., the interior of an edge or a face can have at most one hanging vertex.
In particular, the hanging vertex corresponds to the center of an edge or of a face of another element.
Mesh of this kind are usually obtained when AMR is employed.

As we work with low-order FE method and Lagrangian bases,
nodes will be placed at the vertices and, hence, they inherit the same terminology, i.e. regular nodes and hanging nodes.
We denote the set of all nodes by $\isetall$ and its cardinality by $\cardall$.
while we denote the set of regular nodes by $\isetreg$ and its cardinality by $\cardreg$.
We define the set of non-conforming basis functions as
$$\basissetnc=\{ \mbnc{i} \, : \, \domain \to \mathbb R, \text{ such that } \mbnc{i} |_{E} \in \mathbb L_1 \text{ and } \mbnc{i} |_{E} (\mvertex{j}) = \delta_{ij} \},$$
and the set of conforming basis functions as
$$\basissetconf=\{ \mbconf{i} \, \in C^0(\domain), \text{ such that } \mbconf{i} |_{E} \in \mathbb L_1 \text{ and } \mbconf{i} |_{E} (\mvertex{j}) = \delta_{ij} \},$$
with $\delta_{ij}$ begin the Kronecker delta.
A function $\mbnc{i}$ is created by glueing together the elemental shape functions which attains value one at the node with coordinates $\mathbf x_i$.
We respectively define the non-conforming and conforming interpolation spaces as
$$  \mSnc=\mathrm{span} \,\basissetnc  \quad \text{and} \quad \quad \mS^h =\mathrm{span}\, \basissetconf.$$

The space $\mSnc$ contains also some discontinuous functions as degrees of freedom are associated with all the nodes of the mesh.
Its dimension is $\cardall$.
The space $\mS^h$ contains only continuous functions as degrees of freedom are associated only with regular nodes of the mesh.
Its dimension is $\cardreg$.
The function space $\mS^h$ is spanned by the Lagrangian basis functions $\basissetconf$,
but its definition coincides with the one introduced in~\eqref{eq:discretespacedefinition}.

\begin{figure}
\begin{center}
\includegraphics[scale=0.4]{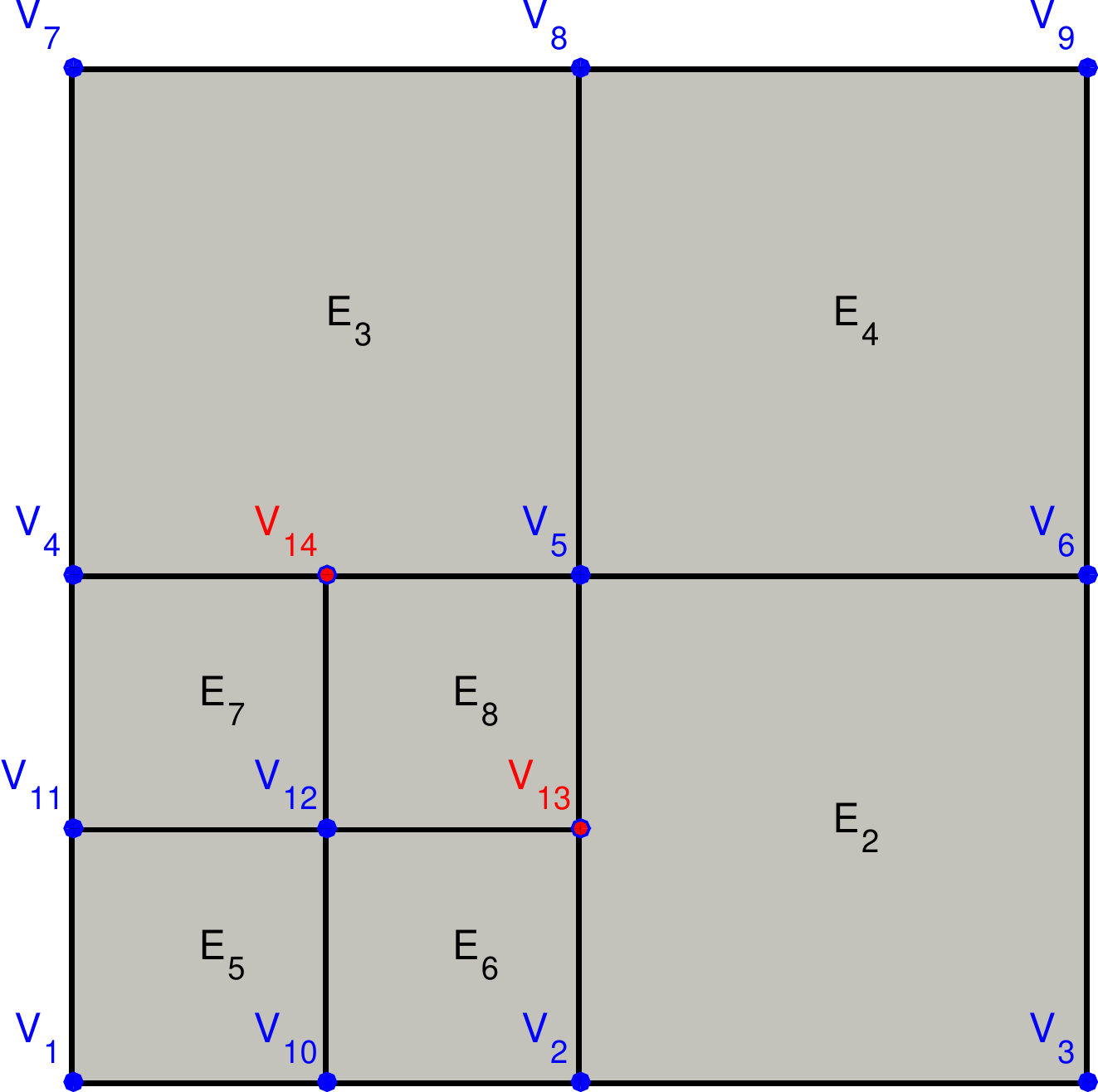}
\end{center}
\caption{Exemplary of a non conforming mesh with  $\isethang=\{13,14\}$ and $\isetreg=\{1,2,3,4,5,6,7,8,9,10,11,12\}$.}
\label{fig:adaptedmesh}
\end{figure}

The difference from the case of conforming meshes is that the support of Lagrangian basis functions $\basissetconf$ is not composed only of the elements that are adjacent to the node.
We explain this concept referring to the simple example of the mesh depicted in Figure~\ref{fig:adaptedmesh}.
Regular nodes are numbered from $1$ to $12$, while hanging nodes have indices $13$ and $14$.

\begin{figure}
\centering
\subfloat[$N^r_2$]{\includegraphics[width= 2.5in]{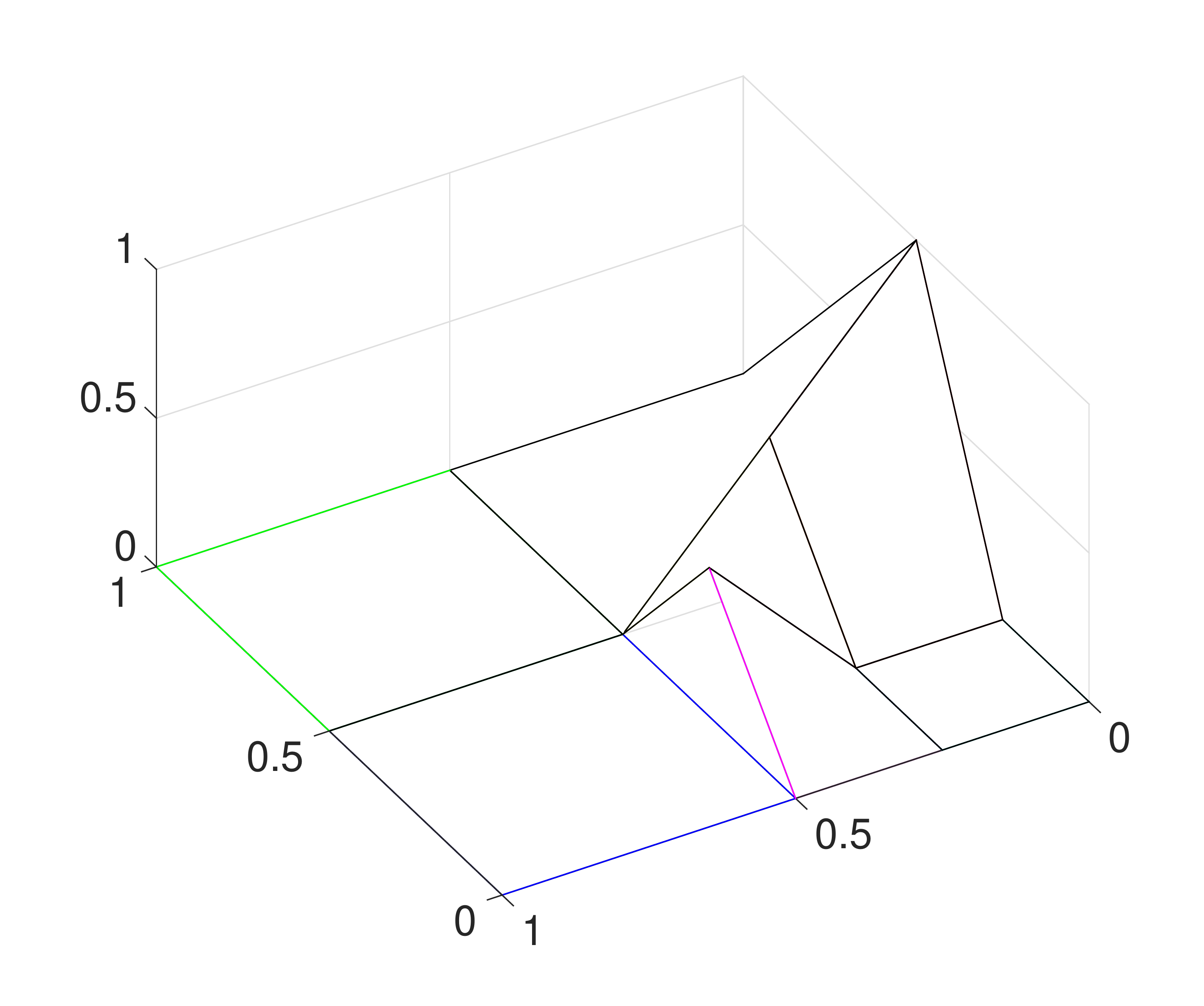}}
\subfloat[$N^r_5$]{\includegraphics[width= 2.5in]{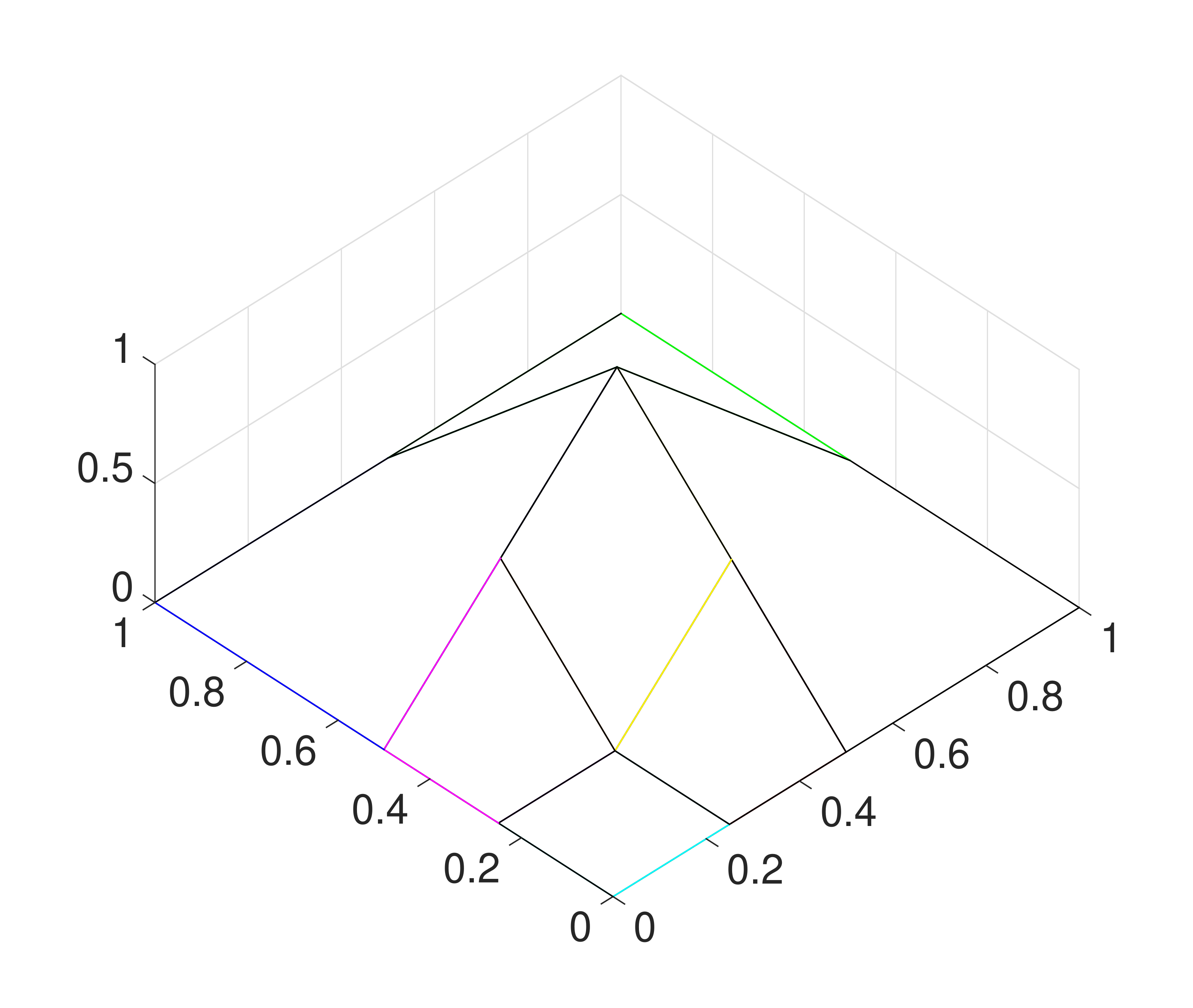}}
\caption{Two basis functions of the conforming space $V^r$ defined over the non-conforming mesh in Figure~\ref{fig:adaptedmesh}.}
\label{fig:twobasis}
\end{figure}
For the basis function $\mbconf{2}$ the support is formed by elements $E_6$ and $E_2$, which are adjacent to $\mvertex{2}$
but also by element $E_8$.
For the basis function $\mbconf{5}$, the support is formed by the adjacent elements $E_2$, $E_4$, $E_3$, and $E_8$ and by the elements $E_6$ and $E_7$.

The sets $\basissetconf(E)$ and $\basissetnc(E)$ denote the conforming and non-conforming basis function on the element $E$.
Observe that $\basissetnc(E)$ is composed of the elemental shape functions.
Observe that $\basissetnc(E)$ is composed of the shape functions defined on $E$.
On a single element, both these sets are a basis for $\mathbb L_1$.
We can introduce a change of basis from $\basissetnc(E)$ to $\basissetconf(E)$.
For example, for the element $E_6$, we have

\begin{equation}
\label{eq:changeofbasis}
\begin{array}{l c l}
\mbconf{10}&=& \mbnc{10},\\
\mbconf{2\phantom{0}}&=& \mbnc{2}+\dfrac{1}{2}\mbnc{13},\\
\mbconf{5\phantom{0}}&=&\dfrac{1}{2}\mbnc{13},\\
\mbconf{12}&=&\mbnc{12}.
\end{array}
\end{equation}

The basis functions of $\mS^h$ are constructed by removing the basis functions of $\mSnc$ at the hanging nodes and enriching the basis functions $\mbconf{a}$ with contributions from the removed ones.
Observe that $\mbconf{2}$ and $\mbconf{5}$ attain value $1/2$ at hanging nodes.

\subsubsection*{Restriction and prolongation operators}

We consider the evaluation of a linear functional $w$ on an element $E_2$, with respect to $\basissetconf(E_2)$ and $\basissetnc(E_2)$, i.e.,
\begin{equation}
\begin{array}{l c l l}
\nonumber
w_i^r & = & w(\mbconf{i}), & i \in \isetreg_{E_2}\\[2mm]
\nonumber
w_i^H & = & w(\mbnc{i}), & i \in \isetall_{E_2}.
\end{array}
\end{equation}
Given the change of basis \eqref{eq:changeofbasis}, we can define the elemental restriction operator $\mrestrictionoperator{E_2}$ as
\begin{align}
\nonumber
w^r_{10}&= \phantom{\dfrac{1}{2}} w^H_{10}\\
\nonumber
w^r_{2\phantom{0}}&= \phantom{\dfrac{1}{2}}  w^H_2 +\dfrac{1}{2}w^H_{13}\\
\nonumber
w^r_{5\phantom{0}}&=\dfrac{1}{2}w^H_{13}\\
w^r_{12}&= \phantom{\dfrac{1}{2}}  w^H_{12}.
\label{eq:HiertoNode}
\end{align}
With a little abuse of notation, i.e., denoting the algebraic representation of the restriction operator as the restriction operator itself,
equation~\eqref{eq:HiertoNode} admits the following algebraic representation $\mathbf w^r_{E_2} = \mrestrictionmatrix{E_2} \mathbf w^H_{E_2}$:
\begin{equation}
\label{eq:trasformazione_1}
\begin{pmatrix}
\nonumber
w^r_{10}\\[1mm]
w^r_{2\phantom{0}}\\[1mm]
w^r_{5\phantom{0}} \\[1mm]
w^r_{12}
\end{pmatrix}=
\begin{pmatrix}1&0&0&0\\[0.5mm]
\nonumber
0&1&1/2&0\\[0.5mm]
0&0&1/2&0\\[0.5mm]
0&0&0&1
\end{pmatrix}
\begin{pmatrix}w^H_{10}\\[1mm] w^H_{2\phantom{0}} \\[1mm] w^H_{13} \\[1mm] w^H_{12}
\end{pmatrix}
\end{equation}
The elemental prolongation operator $\mprolongationoperator{E_2}$ is the adjoint of $\mrestrictionoperator{E_2}$.
For elements that do not have hanging nodes, the elemental prolongation and restriction operators are the identity.
Putting together all the elemental operators we can construct the global restriction and prolongation operators,
denoted by $\mrestrictionoperator{}$ and $\mprolongationoperator{}$, respectively.
Their algebraic representations are still denoted with the same symbols
$\mrestrictionmatrix{} \in \mathbb R^{n^r \times n^H}$ and \mbox{$\mprolongationmatrix{} \in \mathbb R^{n^H \times n^r}$.}

\subsection{DMP conditions on non-conforming meshes}

We refer to an abstract variational problem:
\begin{equation}
\label{eq:abstractweak}
\begin{array}{l}
\text{Find } p\in \mU  \, \text{ such that }\\[2mm]
{\displaystyle a(p,q)=f(q) \quad \forall{q} \in \mV},
\end{array}
\end{equation}
for which we assume that a (continuous) maximum principle holds.
The solution of the problem above is equivalent to the solution of the linear system
\begin{equation}
\label{eq:examplelocal}
\mathbf{A}^r \mathbf{u}^r  = \mathbf f^r,
\end{equation}
where
\begin{equation}
\begin{array}{ l c l l}
\nonumber
\matent{A^r}{i}{j} & = & a( \mbconf{j} ,\mbconf{i}), &  \quad \text{for} \,\, i,j \in \isetreg,\\[2mm]
\nonumber
\matent{f^r}{i}{} &= & f(\mbconf{i}), & \quad \text{for} \,\, i \in \isetreg.
\end{array}
\end{equation}
In the process of assembling FE matrices, only elemental shape functions are available.
Hence, it is possibile to easily assemble
\begin{equation}
\begin{array}{ l c l l}
\nonumber
\matent{A^H}{i}{j} & = & a( \mbnc{j} ,\mbnc{i}), &  \quad \text{for} \,\, i,j \in \isetall,\\[2mm]
\nonumber
\matent{f^H}{i}{} &= & f(\mbnc{i}), & \quad \text{for} \,\, i \in \isetall.
\end{array}
\end{equation}
The stiffness matrix and right-hand-side associated with the conforming discretization space can be obtained by applying the following transformations
\begin{equation}
\label{eq:intA}
\mathbf{A}^r= \mrestrictionmatrix{} \mathbf{A}^H \mprolongationmatrix{}, \quad  \mathbf{f}^r=\mrestrictionmatrix{} \mathbf{f}^H.
\end{equation}
Once the solution $\mathbf p^r$ is found, it may be interpolated on the non-conforming space $\mSnc$ applying
$$
\mathbf{p}^H=\mprolongationmatrix{} \mathbf{p}^r.
$$
In order to ensure the DMP for problem~\eqref{eq:abstractweak},
the conditions of the Theorem~\ref{theo:3} has to be verified for the matrix in \eqref{eq:examplelocal}.
Even if the matrix $\mathbf A^H$ would fulfill such conditions,
the operation~\eqref{eq:intA} may introduce positive extra-diagonal entries.

\subsection{Finite element assembly on non-conforming meshes}

The FE assembly of the problem~\eqref{eq:examplelocal} is realized on each element $E$ constructing the matrix $\mathbf A^H_E$ and the vector $\mathbf f^H_E$.
The restriction operations~\eqref{eq:intA} are hence realized per element, resulting in
\begin{equation}
\nonumber
\begin{array}{l c c}
\mathbf{A}_E^r &= &\mrestrictionmatrix{E}\mathbf{A}_E^H\mprolongationmatrix{E},\\[2mm]
\mathbf{f}_E^r &= & \mrestrictionmatrix{E}\mathbf{f}_E^H.
\end{array}
\end{equation}
The operation above, may introduce positive extra-diagonal entries in $\mathbf{A}_E^r$, even if the matrix $\mathbf{A}_E^H$ has none.
Hence, one needs to check if  the assumptions of Theorem~\ref{theo:3} are satisfied for the matrices  $\mathbf{A}_E^r$ and in case add the elemental contribution of the algebraic diffusion operator.
In Figure~\ref{algo:1}, we report the steps to assemble a stiffness matrix that ensures the DMP.
In order to ensure the DMP, we check that all local stiffness matrices have positive diagonal entries, and negative extra-diagonal entries.
In this way, the stiffness matrix, before setting the Dirichlet boundary conditions, will have the same properties,
and, finally, the imposition of Dirichlet boundary conditions makes $\mathbf A^r$ an M-matrix.

\begin{figure}[ht!]
\begin{algorithm}[H]
  \SetAlgoLined
  \KwData{Bilinear form $a(\,\cdot \,,\,\cdot \,)$, linear functional $F(\,\cdot \,)$, mesh $\mesh$}
  \KwResult{The associated FE matrix $\mathbf A^r$ and right-hand-side $\mathbf f^r$, with $\mathbf A^r$ which satisfies the DMP condition.}
  \For{$E \leftarrow \mesh$}
  {
  Assemble $\mathbf A^H_E$ and $\mathbf f^H_E$\;
 $\mathbf{A}_E^r =\mrestrictionmatrix{E}\mathbf{A}_E^H\mprolongationmatrix{E}$ and $\mathbf{f}_E^r=\mrestrictionmatrix{E}\mathbf{f}_E^H$\; 
\If{$\mathbf{A}_E^r$ has positive extra-diagonal entries}
{
Compute $\mathbf{S}_E^{A^r_E}$\;
$\mathbf{A}_E^r \leftarrow \mathbf{A}_E^r+\mathbf{S}_E^{A^r_E}$\;
}
Modify rows of $\mathbf A^r$ and $\mathbf f^r$ to impose Dirichlet boundary conditions\;
  }
  \caption{Algorithm to assemble problem~\eqref{eq:abstractweak} on a non-conforming mesh and ensuring the DMP.}
  \label{algo:1}
\end{algorithm}
\end{figure}

\subsubsection{Failure of Discrete Maximum Principle conditions}

In order to justify the use of the Algorithm~\ref{algo:1}, we present the assembly of two matrices associated with a flow problem and a transport operator.
We consider a square element, whose nodes are numbered as follow $\mvertex{1}=(0,0)$, $\mvertex{2}=(1,0)$, $\mvertex{3}=(0,1)$, $\mvertex{4}=(1,1)$.
The node $\mvertex{1}=(0,0)$ is a hanging node and the two regular nodes connected to it are $\mvertex{0}=(-1,0)$ and $\mvertex{2}$.
The restriction operator that goes from a space where the index set is $J^H_E=(1,2,3,4)$ to a space where the index set is $J^r_E=(0,2,3,4)$ is
\begin{equation}
\nonumber
\mrestrictionoperator{E}=
\begin{pmatrix}1/2 &0&0&0\\[0.5mm]
1/2&1&0&0\\[0.5mm]
0&0&1&0\\[0.5mm]
0&0&0&1
\end{pmatrix}.
\end{equation}

\paragraph{Diffusion operator}

We consider the weak form associated
\begin{equation}
\nonumber
\int_E \nabla \mbnc{i} \cdot \nabla \mbnc{i} \mbox{d}\Omega,
\end{equation}
and the associated local stiffness matrix
\begin{equation}
\label{eq:localtrans}
\mathbf A_E^H=
\begin{pmatrix}
\phantom{-}2/3&-1/6&-1/6& -1/3\\[0.5mm]
-1/6&\phantom{-}2/3&-1/3&-1/6\\[0.5mm]
-1/6&-1/3&\phantom{-}2/3&-1/6\\[0.5mm]
-1/3&-1/6&-1/6&\phantom{-}2/3
\end{pmatrix}.
\end{equation}
The matrix $\mathbf A_E^H$ has positive diagonal entries and negative extra-diagonal entries, as it usually happens for diffusion operators.
After applying the transformation~\eqref{eq:localtrans},
we obtain a local stiffness matrix evaluated with respect to the regular basis that reads
\begin{equation*}
\mathbf A_E^r=
\begin{pmatrix}
\phantom{-}2/3\phantom{0}     &       \phantom{-}1/12     &     -1/3       &    -5/12    \\[0.5mm]
\phantom{-}1/12     &      \phantom{-}1/6\phantom{0}    &       -1/6       &    -1/12    \\[0.5mm]
-1/3\phantom{0}      &     -1/6\phantom{0}       &     \phantom{-}2/3     &      -1/6\phantom{0}     \\[0.5mm]
-5/12    &      -1/12    &      -1/6       &     \phantom{-}2/3\phantom{0} 
\end{pmatrix}.
\end{equation*}
Hence, in order to have a local stiffness matrix that ensures the DMP conditions, the following algebraic diffusion operator has be added to $\mathbf A_E^r$.
\begin{equation*}
\mathbf S_E^{A_E^r}=
\begin{pmatrix}
\phantom{-}1/12     &       -1/12     &    0      &   0    \\[0.5mm]
-1/12     &      \phantom{-}1/12    &       0    &    0  \\[0.5mm]
0      &     0      &     0     &      0    \\[0.5mm]
0 &      0 						&     0	 &   0
\end{pmatrix}.
\end{equation*}

\paragraph{Transport operator}

We consider now the weak form associated to the transport operator
\begin{equation}
\nonumber
\int_E \mbnc{j} \dfrac{\partial \mbnc{i}}{\partial x_1}
\end{equation}
which corresponds to having chosen a velocity vector of the form $\mvecfun{u}=(1,0)$.
The elemental stiffness matrix assembled on the non-conforming basis functions and the associated algebraic diffusion operator are

\begin{equation*}
\mathbf A_E^H=\frac{1}{12}
\begin{pmatrix}
-2      &       -2      &       -1     &        -1        \\[0.5mm]
\phantom{-}2         &    \phantom{-} 2        &   \phantom{-}   1    &    \phantom{-}      1      \\[0.5mm]  
  -1         &    -1        &     -2        &     -2        \\[0.5mm]
       \phantom{-}1     &         \phantom{-}1   &         \phantom{-}  2     &       \phantom{-}  2 
  \end{pmatrix},
\quad
\mathbf S_E^{ A_E^H}=\frac{1}{12}
\begin{pmatrix}
 \phantom{-}3 & -2 &  \phantom{-}0 & -1 \\[0.5mm]
-2 &  \phantom{-}4 & -1 & -1 \\[0.5mm]
 \phantom{-}0 & -1 &  \phantom{-}3 & -2 \\[0.5mm]
-1 & -1 & -2 &  \phantom{-}4 
\end{pmatrix}.
\end{equation*}
The matrix $\mathbf A_E^H$ has negative diagonal entries and positive extra-diagonal entries.
Moreover, it is not diagonally dominant.
On the other hand, the matrix $\mathbf A_E^H + \mathbf S_E^H$ respects the conditions of Theorem~\ref{theo:3},
but the computation of $ \mrestrictionmatrix{E} (\mathbf A_E^H + \mathbf S_E^H) \mprolongationmatrix{E} $ gives an elemental stiffness matrix which does not respect such conditions.
On the other hand, the computation of $\mathbf A_E^r = \mrestrictionmatrix{E} \mathbf A_E^H \mprolongationmatrix{E} $ and the corresponding algebraic diffusion operator are

\begin{equation*}
\mathbf A_E^r=\frac{1}{24}
\begin{pmatrix}
-1 & -3 & -1 & -1 \\[0.5mm]
 \phantom{-}1 &  \phantom{-}3 &  \phantom{-}1 &  \phantom{-}1 \\[0.5mm]
-1 & -3 & -4 & -4 \\[0.5mm]
 \phantom{-}1 &  \phantom{-}3 &  \phantom{-}4 &  \phantom{-}4 
\end{pmatrix},
\quad
\mathbf S_E^{A_E^r}=\frac{1}{24}
\begin{pmatrix}
 \phantom{-}2 & -1 &  \phantom{-}0 & -1 \\[0.5mm]
-1 &  \phantom{-}5 & -1 & -3 \\[0.5mm]
 \phantom{-}0 & -1 &  \phantom{-}5 & -4 \\[0.5mm]
-1 & -3 & -4 &  \phantom{-}8 
\end{pmatrix}.
\end{equation*}
The sum of $\mathbf A_E^r$ and $\mathbf S_E^{A_E^r}$ gives an M-matrix.

\section{Conservation properties of continuous finite element methods}
\label{sec:conservation}
In this section, we derive the global and local conservation properties of FE discretizations on non-conforming meshes. 
We refer to a diffusion problem, such as the flow problem, but they can be equally proven for all the cases considered in~\cite{hughes}.
The idea at the base of the proofs is that the basis functions form a partition of unity, i.e., the sum of the evaluation of all the basis functions at a point $\mathbf x$ is one.
Equivalently, we can observe that the function one belongs to the space of test functions or that the sum of the rows of the stiffness matrix is zero.

While the proof of the global conservation properties is identical to the one for conforming meshes,
the main difference in the proof of the local conservation properties is in the definition of a suitable space of test functions.
The use of this definition allows us to show that the FE method is conservative across any interface that separates two subsets $\domain$,
without the requirement that such an interface coincides with the sides of some elements of the mesh.

\paragraph{Global conservation properties}

The function spaces defined in equation~\eqref{eq:discretespacedefinition} are spanned by the functions $\mbconf{i}$.
In order to show the global conservation, we introduce the flux $\mflux{\domain}{D}$ along the Dirichlet boundary $\Gamma_D$, and restate problem~\eqref{eq:flowweakdiscrete} as follows:
\begin{equation}
\label{eq:auxflow_global}
\begin{array}{l}
\text{Find}\,\, (\hh{u} , \mflux{\domain}{D}) \in \hh{U}\times \hh{G}\,\,\text{such that}\\[2mm]
{\displaystyle d( \hh{p}, \hh{v}) +(\mflux{\domain}{D} ,\hh{v})_{D} = f( \hh{v})  \quad\forall \hh{v}\in \hh{\mS} }.
\end{array}
\end{equation}
Setting $\hh{v}=1$ into the above equation, we obtain
$$ \int_{\Gamma_D}\mflux{\domain}{D}\,\text{d}{\Gamma} + \int_{\Gamma_N} h  \,\text{d}\Gamma=0.$$
This shows that that $\mflux{\domain}{D}$ is the flux associated to the Dirichlet boundary conditions that makes the method globally conservative.

Formulation~\eqref{eq:auxflow_global} can be split into two sub-problems:\\
$\text{Find}\,\, (\hh{u} , \mflux{\domain}{D}) \in \hh{U}\times \hh{G}\,\,\text{such that}$
\begin{align}
\label{eq:auxflow_global_split_a}
d( \hh{p}, \hh{v}) - f( \hh{v})&=0 & \forall \hh{v} \in \hh{V},\\
\label{eq:auxflow_global_split_b}
(\mflux{\domain}{D} ,\hh{v})_{D} &= d( \hh{p}, \hh{v}) - f( \hh{v}) & \forall \hh{v}\in \hh{G},
\end{align}
which admit the following algebraic representation:
\begin{equation}
\label{eq:flowdiscretehughes}
\begin{pmatrix}
\mathbf D_{II} & 0 \\
-\mathbf D_{ID} & \mathbf B_{DD}
\end{pmatrix}
\begin{pmatrix}
\mathbf p_{I} \\
\mathbf h_D 
\end{pmatrix}
=
\begin{pmatrix}
\mathbf f_{I} -  \mathbf D_{ID}  \mathbf g_D \\
-\mathbf f_D 
\end{pmatrix},
\end{equation}
In the linear system above, we have
\begin{equation*}
\begin{array}{l c l  c  l}
\matent{B_{DD}}{i}{j} & = & (\mbconf{j},\mbconf{i})_{\Gamma_D} &  \quad & i,j\in \isetdir,\\
\matent{D_{II}}{i}{j} & = &d(\mbconf{j},\mbconf{i}) &  \quad & i,j\in \isetrest,\\
\matent{D_{ID}}{i}{j} & =& d(\mbconf{j},\mbconf{i}) &  \quad & i\in \isetrest, \, j \in \isetdir,\\
\matent{f_{I}}{i}{} & =& f(\mbconf{i}) &  \quad & i\in \isetrest.\\
\end{array}
\end{equation*} 
The matrix $\mathbf B_{DD}$ is the boundary mass on $\Gamma_D$.
From \eqref{eq:auxflow_global_split_b} or, equivalently, from the second row of \eqref{eq:flowdiscretehughes},
we notice that the computation of the boundary flux associated to $\Gamma_D$ requires the inversion of the boundary mass matrix $\mathbf B_{DD}$,
which can be thought as a \textit{post-processing} evaluation.



\paragraph{Local conservation properties}

Let $\domain_1 \subset \domain$ and $\domain_2 \subset \domain$ be two non-overlapping subdomains such that $\overline{\domain}_1\cup \overline{\domain}_2= \overline{\domain}$.
We denote the boundary of the subdomain $\domain_1$ by $\Gamma_1$,
the boundary of the subdomain $\domain_2$ by $\Gamma_2$,
and the interface between the $\domain_1$ and $\domain_2$ (i.e., $\Gamma_1 \cap \Gamma_2$) by $L$ .
The subdomain $\domain_1$ could be interior to $\domain$
or intersect the boundary of $\domain$ in either a finite number of isolated points or a set with a finite measure.

For the sake of simplicity, we assume that $\Gamma_1$ only overlaps with the Neumann boundary $\Gamma_N$ and denote by $\Gamma_1^N=\Gamma_1\cap\Gamma_N$ their intersection such that $\Gamma_1=L \cup \Gamma_1^N$. In the same way, the boundary $\Gamma_2\setminus L$ is split into the Neumann and the Dirichlet portions and defined as $\Gamma_2=L \cup \Gamma_2^N \cup \Gamma_D$.

We refer to $\isetreg_L$ as the index set of the conforming basis functions whose support is non-null on the interface ${L}$, i.e.,
$$\isetreg_L = \{ i \, : \, | \text{supp} (\mbconf{i}) \cap L | >0\}.$$
and introduce the function space  $G_{L}^h =  \mathrm{span}  \{ \mbconf{i} \}_{\{i \in \isetreg_L \}}$. We also define the bilinear form and the linear functional restricted to the subdomain $\domain_a$ for $a=\{ 1,2\}$ as follows:
$$d_a({p},{q}) =  \int_{\domain_a} k \, \nabla p \cdot \nabla q \,\text{d}\Omega \quad \text{and} \quad f_a(q )=\int_{\Gamma_a^N} h \, {q} \,\text{d}\Gamma,$$
with $a=\{ 1,2\}$.

In order to derive the local conservation properties,
we denote the approximated flux across the boundary $L$ with respect to $\domain_1$ by $\mflux{\domain_1}{L} $.

The boundary flux  $\mflux{\domain_1}{L} $ is the solution of following formulation:
\begin{equation}
\begin{array}{l}
\text{Find} \,(\hh{p},\mflux{\domain_1}{L})\in U^h\times G_{L}^h \, \text{ such that }\\[2mm]
\label{eq:bflow_local_1}
(\hh{v}, \mflux{\domain_1}{L})_{L}=d_1(\hh{p}, \hh{v}) - f_1(\hh{v})\quad\forall \hh{v}\in S^h,
\end{array}
\end{equation}
which can be split into two sub-problems:
\begin{align}
\nonumber
\text{Find} \,(\hh{p},\mflux{\domain_1}{L})\in U^h\times G_{L}^h \, \text{ such that }\\[2mm]
\nonumber
d_1(\hh{p}, \hh{v}) - f_1(\hh{v})&=0 & \forall \hh{v}\in \mS^h- G_{L}^h,\\
\nonumber
(\hh{v}, \mflux{\domain_1}{L})_{L} & = d_1(\hh{p}, \hh{v}) - f_1(\hh{v}) & \forall \hh{v}\in  G_{L}^h.\quad\quad
\end{align}

In the same way, one can define the flux for the complementary domain $\mflux{\domain_2}{L}$ across the boundary $L$ as the solution of the problem: 
\begin{equation}
\begin{array}{l}
\text{Find} \,(\hh{p},\mflux{\domain_2}{L})\in U^h\times G_{L}^h \, \text{ such that }\\[2mm]
\label{eq:bflow_local_2}
(\hh{v},\mflux{\domain_2}{L})_{L}=d_{2}(\hh{p}, \hh{v}) - f_{2}(\hh{v})-(\hh{v},q_D(\Omega))_{\Gamma_D} \quad \forall \hh{v}\in \mS^h,
\end{array}
\end{equation}

We establish the local conservation by first proving that the interface fluxes, $\mflux{\domain_1}{L}$ and $\mflux{\domain_2}{L} $ are conservative with respect to corresponding subdomains $\domain_1$ and $\domain_2$. 
To this aim we select any $\hh{v}\in\mS^h$ such that $\hh{v}_{\domain_a}=1$ with $a\in\{1,2\}$, and obtain:
\begin{equation*}
\int_{L}\mflux{\domain_1}{L}\,\text{d}\Gamma+ \int_{\Gamma_1^N} h \,\text{d}\Gamma=0
\end{equation*}
and
\begin{equation*}
\int_{L}\mflux{\domain_2}{L}\,\text{d}\Gamma+ \int_{\Gamma_2^N} h \,\text{d}\Gamma+\int_{\Gamma_D}q_D^h(\Omega)\,\text{d}\Gamma=0.
\end{equation*}

It only remains to prove the uniqueness of the boundary flux over $L$. Summing together equation~\eqref{eq:bflow_local_1} and equation \eqref{eq:bflow_local_2} and restricting $\hh{v}$ to $G_{L}^h\subset \mS^h$ leads to:
\begin{equation*}
(\hh{v},\mflux{\domain_1}{L}+ \mflux{\domain_2}{L})_{L}=0
\end{equation*}

which states that the two fluxes equilibrate on $L$ and completes the proof.

{\bf Remark 1}  If the vertices of the  interface  $L$ correspond to the vertices the mesh $\mathcal{T}$, the local conservation is attained point-wise:
$$\mflux{\domain_1}{L}+\mflux{\domain_2}{L}=0,$$
and the coefficients of the boundary fluxes can be obtained by computing the inverse of a boundary mass matrix restricted to $L$.

{\bf Remark 2} As the sum of the rows of a discrete diffusion operator is zero, due to the properties listed in Definition~\ref{def:5}, the global and local conservation properties hold also for the stabilized formulations.

\section{Numerical examples}
\label{sec:numerical_results}
In this section, we analyze the DMP and conservation properties of FE discretizations on non-conforming meshes.
The applications will be the coupled flow and transport problems in fractured porous media based on an equi-dimensional representation of fractures.
For the flow problem, we focus both on the conservation properties of the computed fluxes and the positivity of the solution.
For the transport problem, the main focus is to show the positivity of the solution.
Before discussing in detail the numerical results, we briefly describe the AMR strategy employed to generate the non-conforming adapted meshes for domains with fractures.

We consider three numerical examples.
The first example is a regular fracture network~\cite{flemisch2018benchmarks}.
Being the fractures of this example axis-aligned,
we can employ a structured mesh that resolves the interfaces between the fractures and the embedding matrix.
The second example is a two-dimensional adaptation of the benchmark proposed in~\cite{BenchPaper}.
The third example describes a realistic fracture network~\cite{flemisch2018benchmarks}.

We employ the first and the second benchmarks to estimate the conservative boundary fluxes, and the third test case to show the key role that stabilization plays to obtain solutions that satisfy the DMP.
In order to evaluate the accuracy of the boundary fluxes $q_{L}(\domain_i)$ with $i\in\{1,2\}$, we compute the relative $L^2$ error:
 \begin{equation}
 \label{eq:error}
 e_q=\dfrac{\|q_{L}(\domain_i)-q_{R}(\domain_j)\|_{L}}{\|q_{R}(\domain_i)\|_{L}},
 \end{equation} 
where $q_{R}(\domain_i)$ denotes the reference flux, usually obtained on a very fine mesh.
Moreover, we evaluate the total flux $Q_{L}(\domain_j)$ along the boundary $L$ with respect to the sub-domains $\domain_1$ and $\domain_2$ as follows:
 \begin{equation}
 \label{eq:tot_flux}
 Q_{L}(\domain_i)=\int_{L}q_{L}(\domain_i) d\Gamma.
 \end{equation} 

All numerical tests have been implemented in \lq\lq Parrot\rq\rq, an application implemented in the FE framework MOOSE.
Parrot allows for the simulation of several problems in heterogeneous materials, such as flow and transport, Biot's equations in the time and time-frequency domain.
The solution of the linear system arising from the FE discretization is performed using the parallel direct solver MUMPS.

\subsection{Adaptive mesh refinement for\\ heterogenous media}
In this work, we employ the strategy proposed in~\cite{Holliger2000} to create adapted meshes for fracture networks.
The idea of the method is to start from an initial uniform mesh.
Then, given the distribution of embedded fractures $\mathcal F$,
one may iterate over all elements and mark for refinement the ones that have non-empty intersections with at least one fracture.
The uniform refinement of such elements is a step of the AMR procedure.
At each step, we also mark for refinements the elements close to the refined one to keep the meshes $1$-irregular.
Iterating the refinement step several times allows creating sequences of meshes, which become progressively finer close to the interface between the fractures and the matrix,
where the discretization error is usually larger.
As interfaces are not explicitly resolved, material properties may be non-constant over the mesh elements close to the interfaces.
Hence, in the process of assembling the stiffness matrix, the material properties are evaluated at the quadrature points.

We point out that the employed mesh generation is a flexible procedure that allows us to obtain meshes for any distribution of the fractures and with some desired features.
For example, one may refine only elements that have a non-empty overlap with the boundary of fractures or a region around the fracture which does not necessarily coincide with its aperture.

\subsection{Resolved and unresolved meshes}

Meshes generated with the AMR strategy do not resolve the interfaces between matrix and fractures. Hence, we will refer to them as \emph{unresolved}, or \emph{adapted}, meshes.
This kind of meshes is in contrast to the \emph{resolved} meshes, i.e., meshes that explicitly resolve the interfaces.
For resolved meshes, elements close to the interface have edges that lie on the interfaces.

Discretization methods based on resolved meshes are usually referred to as fitted methods,
while methods based on unresolved methods are referred to as unfitted methods.
Despite the lack of regularity of interface problems, fitted methods, in general, provide optimal convergence rates compared to the polynomial order.
On the other hand, standard unfitted FE methods do not provide optimal convergence rates but allow for larger flexibility in the meshing of complicated domains and interfaces.
For the elements on which material properties are not constant,
the standard convergence proofs used for FE are not applicable.


We refer to the resolved meshes as $\mathcal R^{\text{be},\text{fe}}$, where $\text{be}$ is the number of vertical subdivisions in the mesh
and $\text{fe}$ is the number of element crossing the fracture.
We point out that large values of $\text{fe}$ with smaller values of $\text{be}$ give meshes with elongated elements in the fracture,
which may affect the accuracy of the solution.

Unresolved meshes are denoted by $\mathcal U^{\text{be}}_{\text{amr}}$,
where $\text{amr}$ is the number of AMR steps
applied to the uniform background mesh consisting of $\text{be}$ initial subdivisions.

For each example, we report the mesh characteristics, together with the number of elements $n^E$, and the number of elements per fracture aperture $w$.
and for the first two examples, for which a resolved mesh can be created, we perform a comparison of the solutions between the resolved and the unresolved meshes. 

\subsection{Regular fracture network}

The geometry of this test consists in a regular fracture network embedded in square domain with $L_1=L_2 = 1\,[\text{m}]$.
Material properties are listed in Table~\ref{tab:param_3}.
A flow problem on this geometry has been originally presented in~\cite{geiger2013novel}.
A comparison between different solution methods for both hybrid- and equi-dimensional representations of the fracture have been presented in~\cite{flemisch2018benchmarks}.
A hybrid-dimensional representation for the transport has been proposed in~\cite{odsaeter}.

We consider several resolutions for the resolved and the unresolved meshes.
For the resolved meshes, we choose $\text{be}\in \{ 80, 160, 320, 640, 1280, 2560, 4000, 5120 \}$ and $\text{fe} \in \{ 1,2,4,8,16 \}$.
In order to study the accuracy of the boundary fluxes,
we use a reference solution computed on a resolved mesh with $5120$ elements in the background and $64$ elements in an area that is the double of the fracture aperture.
This means that $32$ elements are inside the fracture and $16$ elements per side are close to each fracture side.
We remind that the total number of elements is $(\text{be}+3\text{fe})^2$,
while the total number of nodes is $(\text{be}+3\text{fe}+1)^2$.
The characteristics of the unresolved meshes are reported in Table~\ref{tab:meshcharacteristics1}.  
Pressure and concentration profiles are analized along the segments $\linea$, $(x_1=0.5$ [m]$)$, $\lineb$, $(x_2=0.7$ [m]$)$, and $\linec$, $(x_2=0.5$ [m]$)$, and $\lined$, $(x_2=0.75$ [m]$)$ as reported in Figure~\ref{fig:unresolved_regular_p}.

\subsubsection{Flow problem and local conservative flux calculations}

\begin{figure}[hbt!]
\centering
\includegraphics[width=9cm]{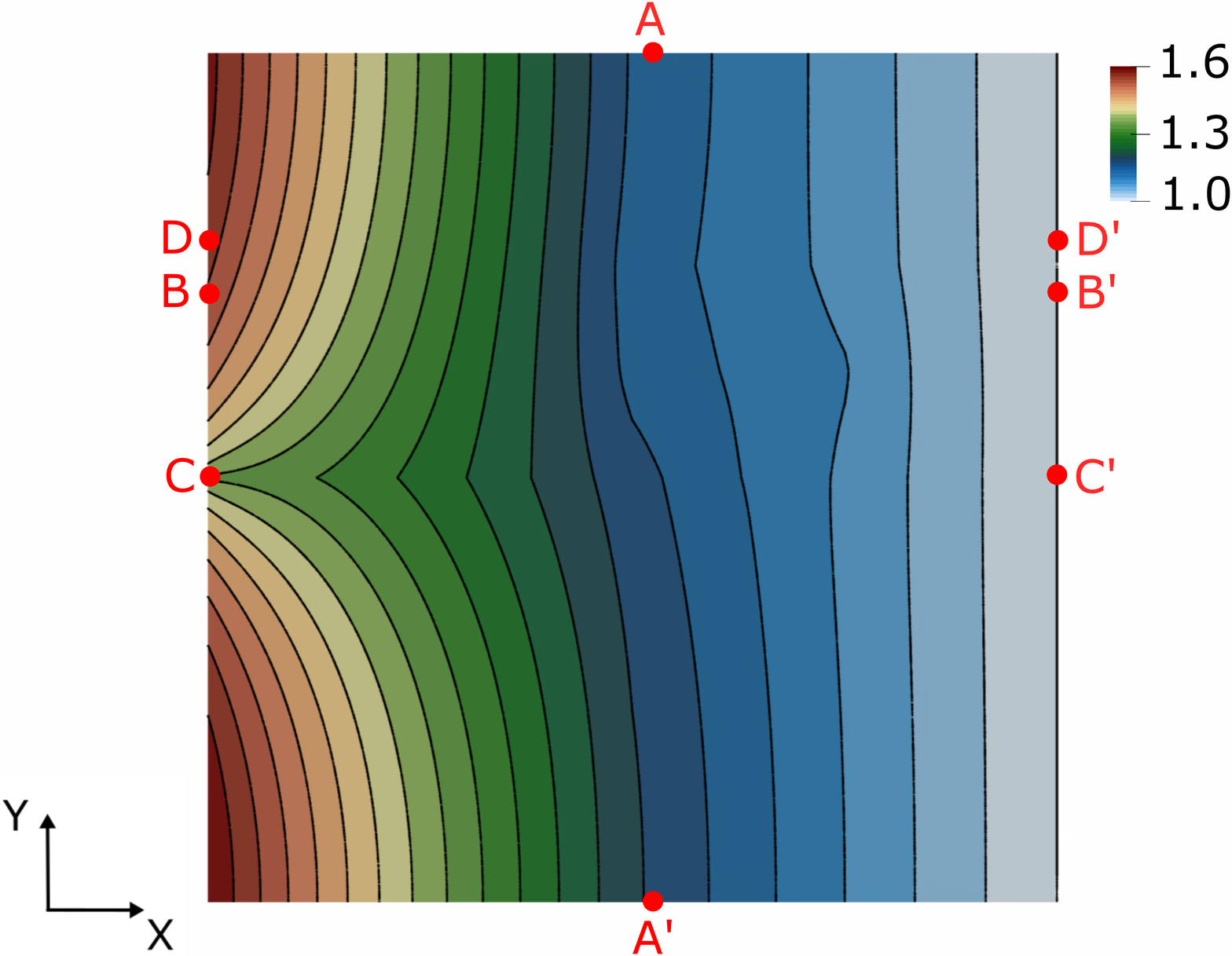}
\caption{Pressure distribution for the test case  $\mathcal{U}^{320}_{8}$. We also report the extrema of the three segments along which the relevant properties are computed.}
\label{fig:unresolved_regular_p}
\end{figure}

For the flow problem, $\Gamma_D$ is the right side at $x_1=L_1$,
while $\Gamma_N$ is composed of the other three sides.
We set $g=1$ [m] and $h=-1$ [m$\times$s$^{-1}$] on the left side of  $\Gamma_N$.
On the top and bottom sides, homogeneous Neumann conditions are imposed.
Figure~\ref{fig:unresolved_regular_p} shows the spatial distribution of the pressure computed on the mesh $\mathcal{U}^{320}_{8}$.

\begin{table}
\centering
\caption{Material properties employed for the benchmark  \lq\lq Regular fracture network\rq\rq.\label{tab:param_3}}
 \begin{tabular}{ | l | l | l | l |}
\hline
 Property & Symbol & Value & Unit \\
   \hline\hline
 Fracture  thickness & $\delta$ & $1 \times 10^{-4}$ & m\\
  \hline
   Matrix permeability & $k_m$ & $1 \times 10^{4}$ & m$\times$s$^{-1}$\\
   \hline
  Fracture permeability & $k_f$ & $1$ & m$\times$s$^{-1}$ \\
   \hline
   Matrix porosity & $\phi_m$ & $1$ & \\
   \hline
  Fracture porosity & $\phi_f$ & $1$ & \\
  \hline
  \end{tabular}
\end{table}

\begin{figure}[hbt!]
\centering
\includegraphics[width=12cm]{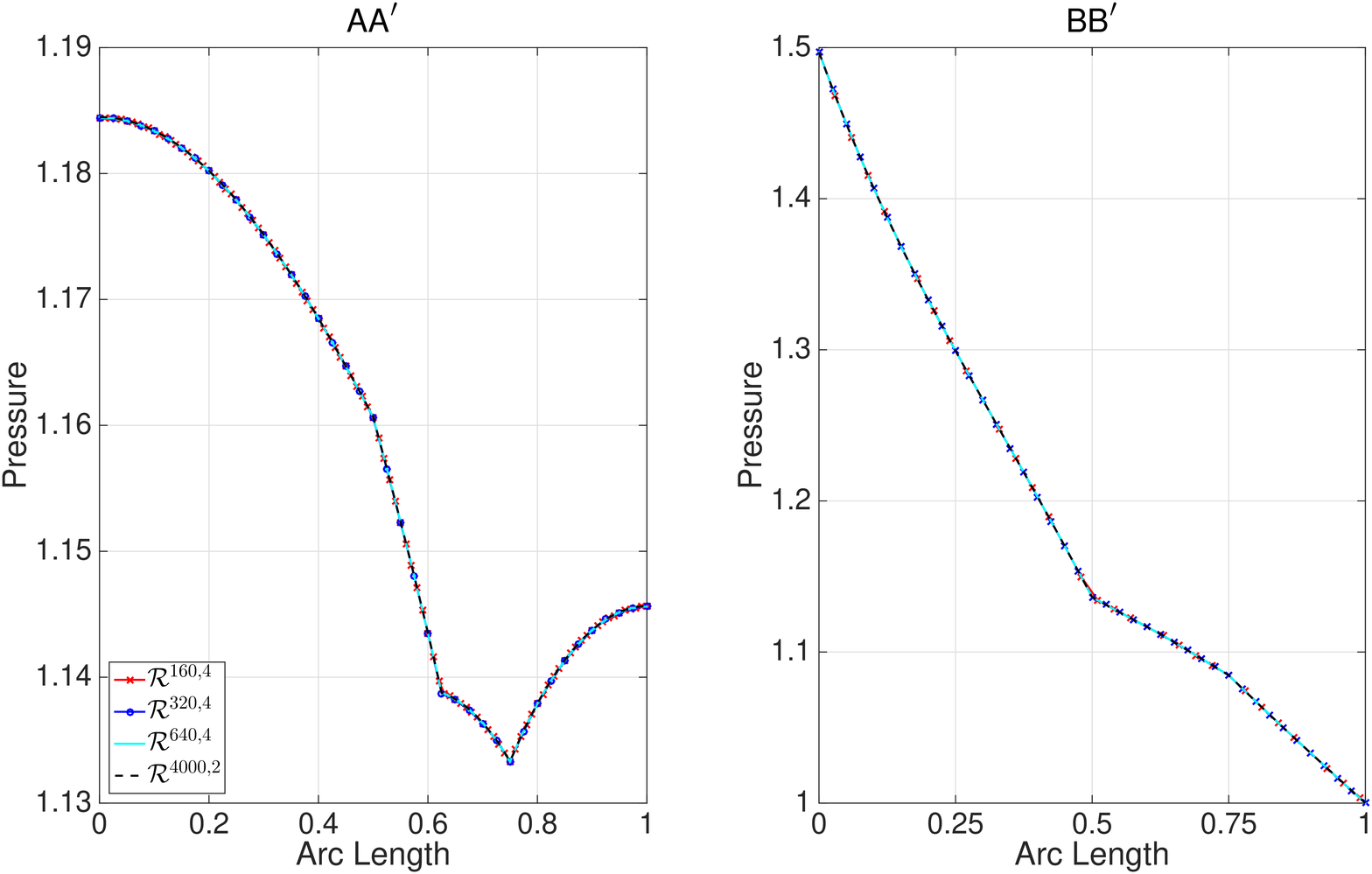}
\caption{Comparison of pressure values along two lines for the resolved test cases.}
\label{fig:regularresolved}
\end{figure}

\begin{figure}[hbt!]
\centering
\includegraphics[width=12cm]{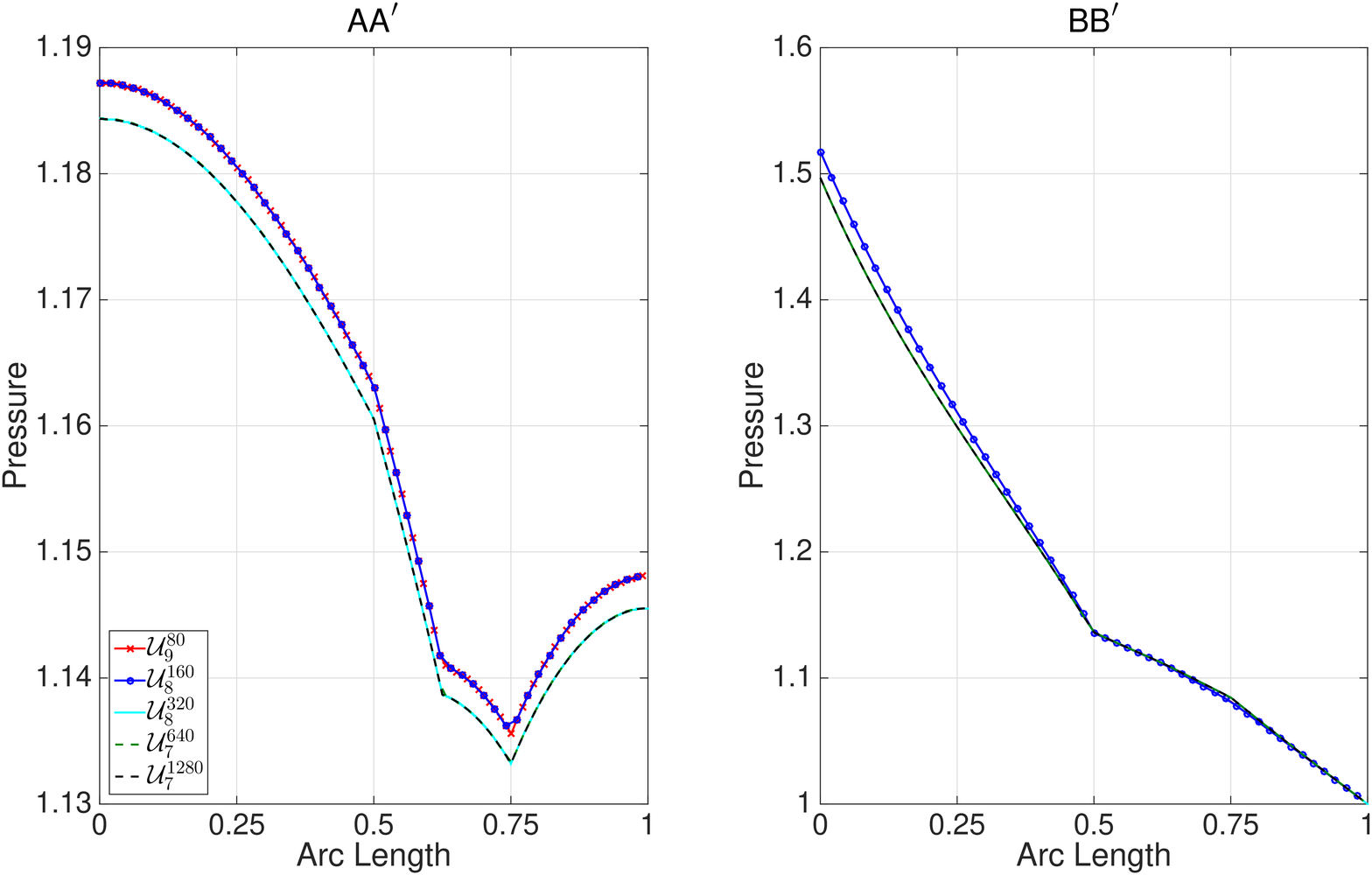}
\caption{Comparison of pressure values along two lines for the unresolved test cases.}
\label{fig:regularunresolved}
\end{figure}

\begin{figure}[hbt!]
\centering
\includegraphics[width=12cm]{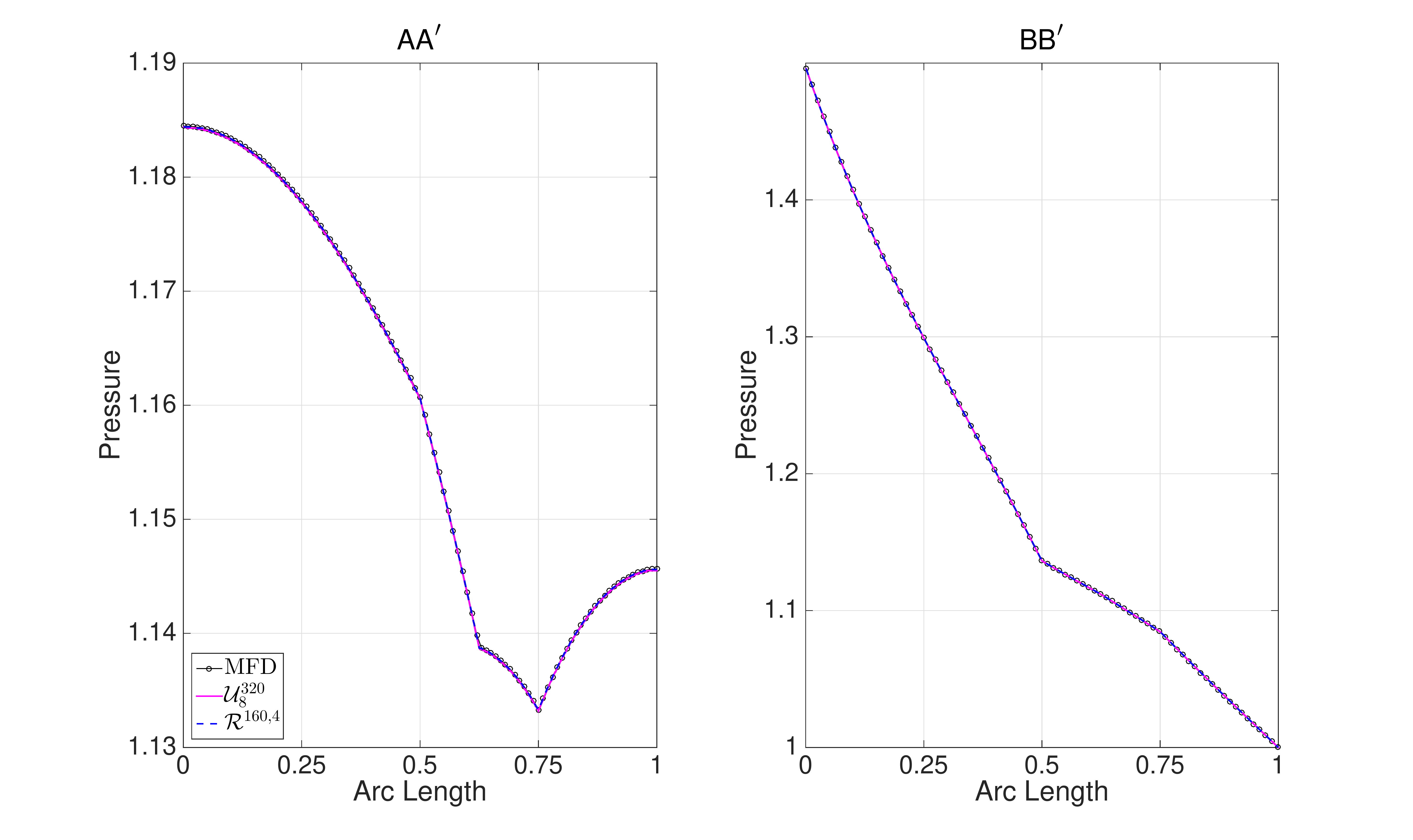}
\caption{Comparison of pressure values along two lines between resolved and unresolved test cases. The MFD solution computed in~\cite{flemisch2018benchmarks} is reported in black.}
\label{fig:regularComparison}
\end{figure}

We compare the pressure on some resolved and unresolved meshes computed along with the segments $\linea$ and $\lineb$ and depict all the results in Figures~\ref{fig:regularresolved} and~\ref{fig:regularunresolved}, respectively.
We observe that for resolved meshes no visible differences are present along both the lines. On the other hand, the accuracy of the solutions on unresolved meshes is mostly affected by $w$,
the number of elements per fracture width.
For this case, solutions obtained with meshes with $w>8$ present identical pressure distributions.
In Figure~\ref{fig:regularComparison}, we report a comparison between the solutions computed on the unresolved and the resolved meshes.
We observe that the two solutions coincide with one of the solutions presented in~\cite{flemisch2018benchmarks},
which has been computed with a mimetic finite difference (MFD) method on a very fine mesh,
consisting of $1\,175\,056$ elements and $10$ elements per each fracture in its normal direction.

To show the conservation properties of the FE on adapted meshes,
we consider two subdivisions of $\domain$.
The first one is the previously-introduced $\domain_m$ and $\domain_f$, for which we compute the total flux along with the interface $\Gamma$.
As a second one, we consider a subdivision in $\domain_1$ and $\domain_2$, which are divided by the segment $L:=\lineb$.
In Table~\ref{tab:flux_mesh_2},
we report the value of the total fluxes $Q_{L}(S)$ with $S \in \{ \domain_1,\domain_2 , \domain_m , \domain_f\}$ obtained by using Equation~\eqref{eq:tot_flux} for some of the resolved and unresolved meshes.
As already proved in Section~\ref{sec:conservation}, all the numerical results confirm that the FE method on adapted meshes attains local conservation properties. Indeed, we see that the value of the flux computed on one side of an interface is matched by the negative value on the other side of the domain.
 The total fluxes computed across the interface $\Gamma$ satisfy the same property.
 
\begin{table}[ht]
\centering
\caption{Total fluxes, Equation~\ref{eq:tot_flux}, for the resolved and unresolved test cases.}
\begin{tabular}{| l | r | r | r | r| }
 \hline
       Mesh & $Q_{\Gamma}(\domain_m)$ & $Q_{\Gamma}(\domain_f)$ & $Q_{\lineb}(\Omega_1)$ & $Q_{\lineb}(\Omega_2)$ \\
 \hline
  \hline
$\mathcal R^{320,4}$&  0.66219 & -0.66219 & -0.11775 & 0.11775\\
$\mathcal R^{640,4}$& 0.66219 & -0.66219 & -0.11775 & 0.11775\\
$\mathcal R^{1280,4}$& 0.66218&-0.66218& -0.11775 & 0.11775\\
$\mathcal R^{4000,2}$& 0.66218 & -0.66218& -0.11775 & 0.11775\\
 \hline
$\mathcal U^{80}_{9}$& 0.65692& -0.65692 & -0.11704 & 0.11704\\
$\mathcal U^{160}_{8}$ &0.656915 &-0.656915 &-0.11691 & 0.11691 \\
$\mathcal U^{320}_{8}$&0.662382 &-0.662382 & -0.11803 & 0.11803\\
$\mathcal U^{640}_{7} $& 0.662381& -0.662381 & -0.11799 & 0.11799\\
$\mathcal U^{1280}_{7}$  & 0.662415&-0.662415 & -0.11794 & 0.11794\\
 \hline
\end{tabular}
\label{tab:flux_mesh_2}
\end{table}

\begin{figure}
\begin{center}
\subfloat[Boundary flux across the region of the boundary $\lineb$ intersecting fractures.]{\includegraphics[width= 5in]{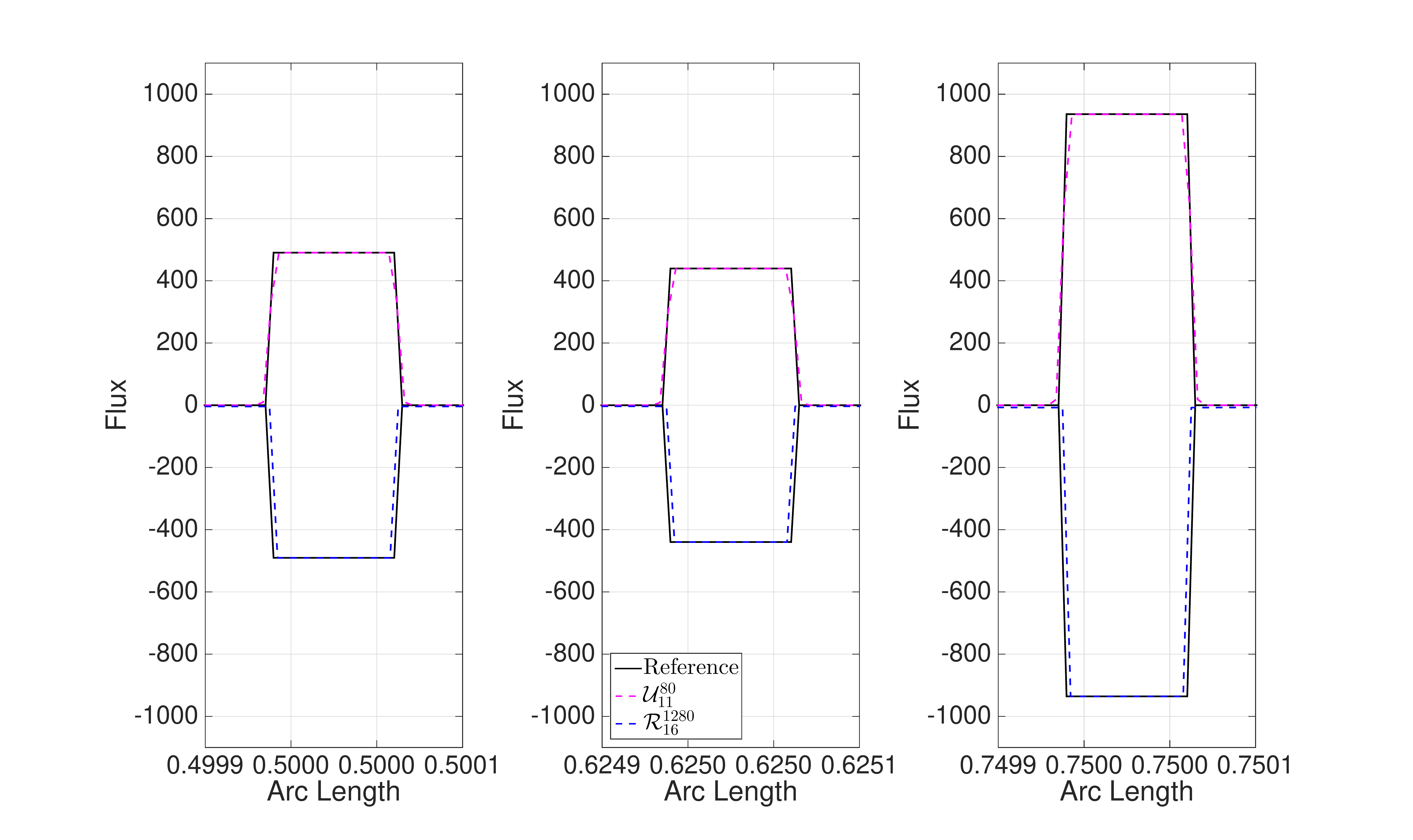}}\\
\subfloat[Boundary flux across the portion of the boundary $\lineb$ lying in the embedding matrix.]{\includegraphics[width= 5in]{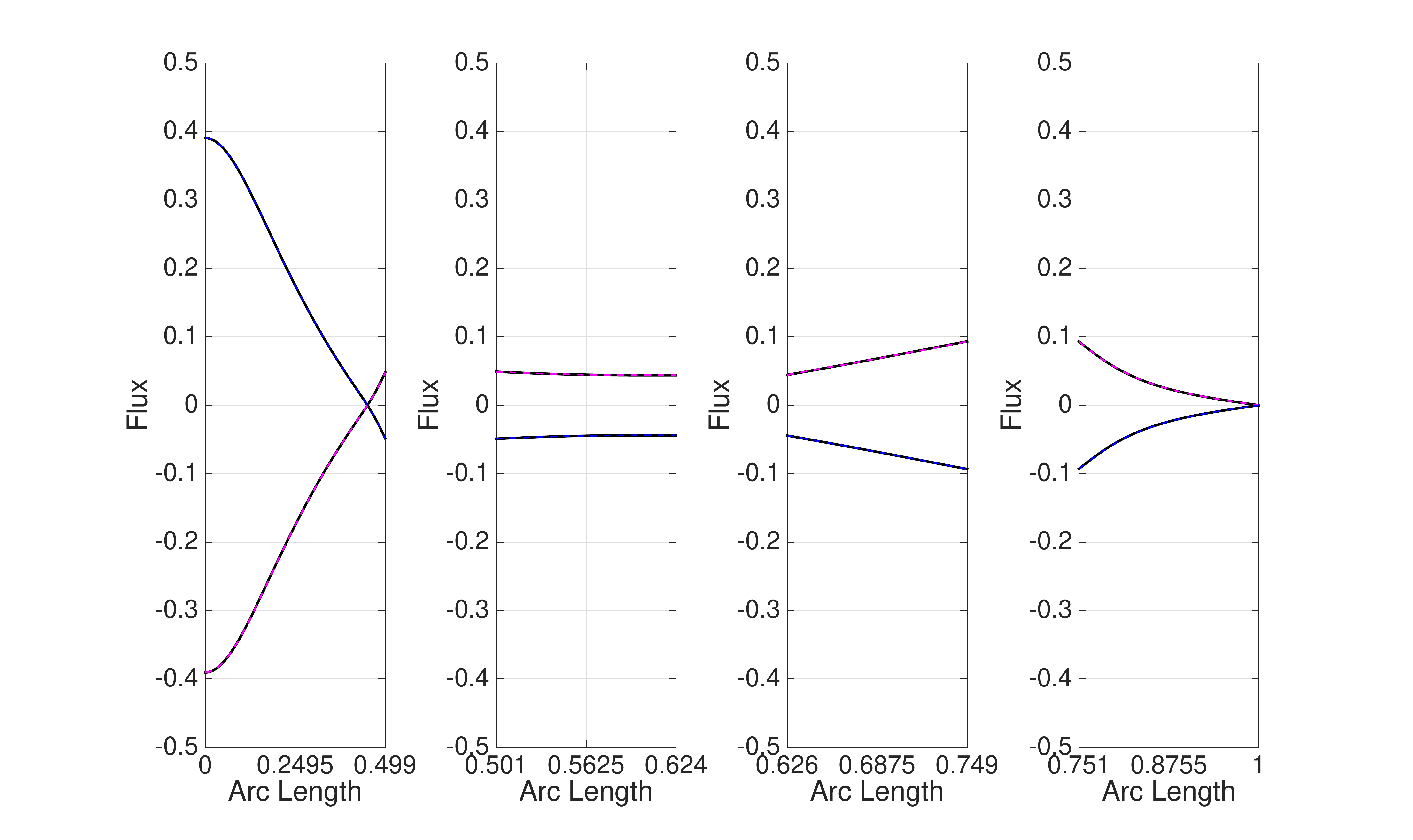}}
\caption{Boundary fluxes across boundary $\lineb$. We report $\mflux{\domain_2}{\lineb}$ for the unresolved mesh, and $\mflux{\domain_1}{\lineb}$ the resolved test cases. The black line is the reference solution.}
\label{fig:comparisonfluxed}
\end{center}
\end{figure}

In Figure~\ref{fig:comparisonfluxed}, we compare the fluxes computed on the unresolved mesh $\mathcal U^{80}_{11}$ and the resolved mesh $\mathcal R^{1280,16}$
against the fluxes computed for the reference solution.
Such mesh resolutions have been computed as they are the coarsest ones for each case that can correctly reproduce the flux.
For reproducing the flux correctly,
unresolved meshes with a coarse background are sufficient but they need a large number of AMR steps, requiring $w>16$.
On the other hand, fluxes computed on the resolved mesh require a much finer background.
Such results suggest that resolved meshes are not enough to correctly reproduce the flux,
but they also require a certain number of elements close to the interfaces.
This fact is also highlighted by the errors reported in Table~\ref{tab:convergenceResolved} and Table~\ref{tab:convergenceNonstab}.
We observe that the error is mostly dominated by the value of $w$ for unresolved discretizations.
For resolved meshes, although the error decreases more regularly, it remains larger even if a larger number of background and fracture elements are employed.

\begin{table}[h!]
\caption{$L^2$ relative error,  Equation~\ref{eq:error}, for the boundary flux $q_{\lineb}(\domain_1)$ computed on the resolved meshes.}
\begin{tabular}{| c | c | c | c | c | c | c | c | }
\hline
   \diagbox{fe:}{be:}
  & 80 & 160 & 320 & 640 & 1280 & 2560 & 5120 \\ 
\hline
\hline
1 & 0.99439  & 0.98905  & 0.97863  & 0.95874  & 0.92264  & 0.86367  & 0.78035  \\
2 & 0.56599  & 0.56301  & 0.56064  & 0.55603  & 0.54743  & 0.53292  & 0.51144  \\
4 & 0.38873  & 0.38697  & 0.38609  & 0.38437  & 0.3811  & 0.37542  & 0.3667  \\
8 & 0.25898  & 0.2578  & 0.25745  & 0.25679  & 0.25551  & 0.25325  & 0.24969  \\
16 & 0.16424  & 0.16339  & 0.16323  & 0.16294  & 0.16239  & 0.16138  & 0.15975  \\
\hline
\end{tabular}
\label{tab:convergenceResolved}
\end{table}

\begin{table}[h!]
\begin{tabular}{| c | c | c | c | c | c | c | }
\hline
   \diagbox{amr:}{be:}
  & 80 & 160 & 320 & 640 & 1280 & 2560 \\ 
\hline
\hline
2 & & & & & & 0.50383  \\
3 & & & & & 0.50383  & 0.3782  \\
4 & & & & 0.50383  & 0.3782  & 0.24934  \\
5 & & & 0.50383  & 0.3782  & 0.24934  & 0.1466  \\
6 & & 0.50383  & 0.3782  & 0.24934  & 0.1466  & 0.061532  \\
7 & 0.50384  & 0.3782  & 0.24934  & 0.1466  & 0.061532  & \\
8 & 0.3782  & 0.24934  & 0.1466  & 0.061532  & & \\
9 & 0.24934  & 0.1466  & 0.061533  & & & \\
10 & 0.1466  & 0.061533  & & & & \\
11 & 0.061535  & & & & & \\
\hline
\end{tabular}
\caption{$L^2$ relative error,  Equation~\ref{eq:error}, for the boundary flux $q_{\lineb}(\domain_1)$ computed on the unresolved meshes.}
\label{tab:convergenceNonstab}
\end{table}

\subsubsection{Transport problem and DMP}

\begin{figure}[hbt!]
\centering
\includegraphics[width=12cm]{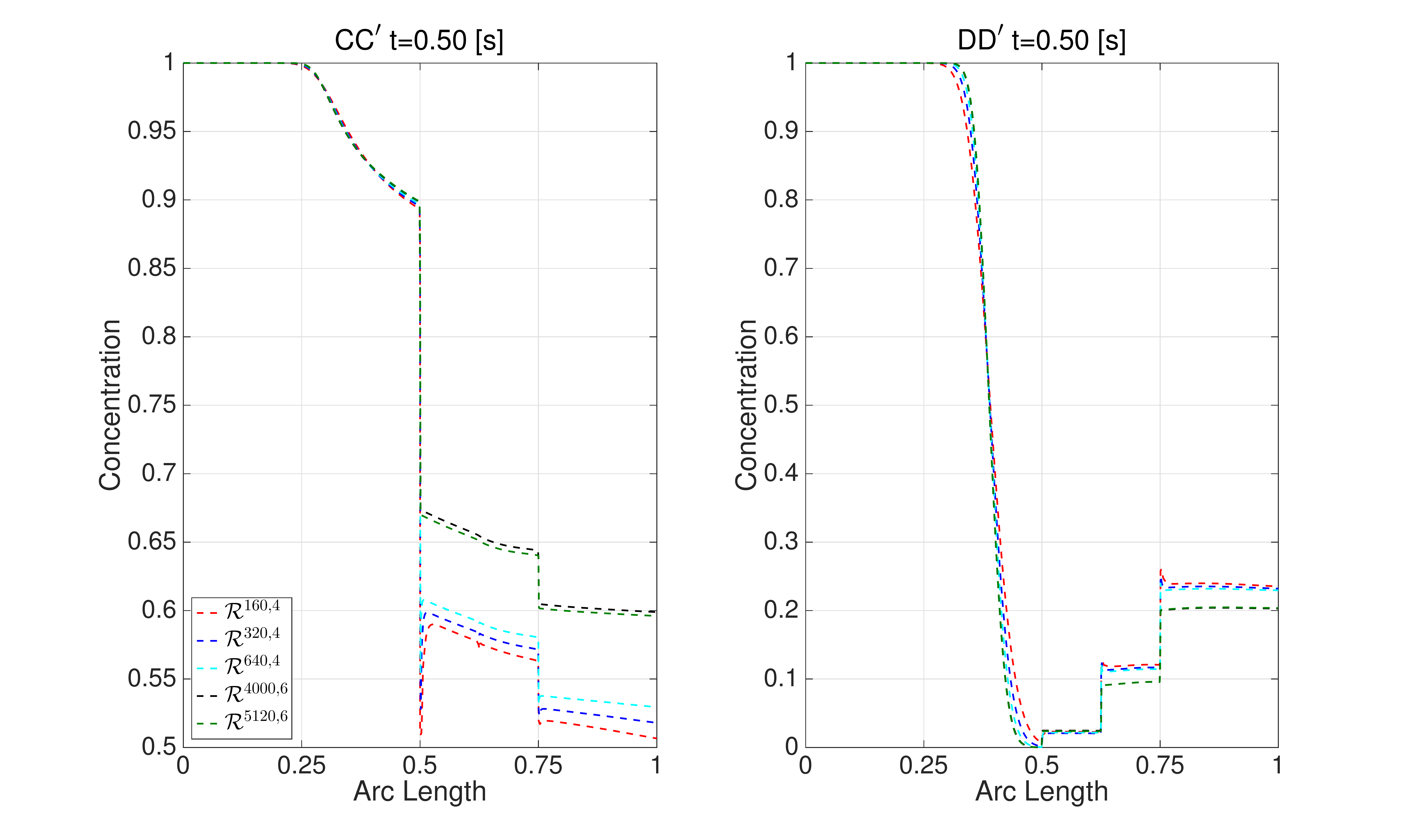} 
\caption{Concentration values along two lines for the resolved test cases.}
\label{fig:regular_c_res}
\end{figure}

\begin{figure}[hbt!]
\centering
\includegraphics[width=12cm]{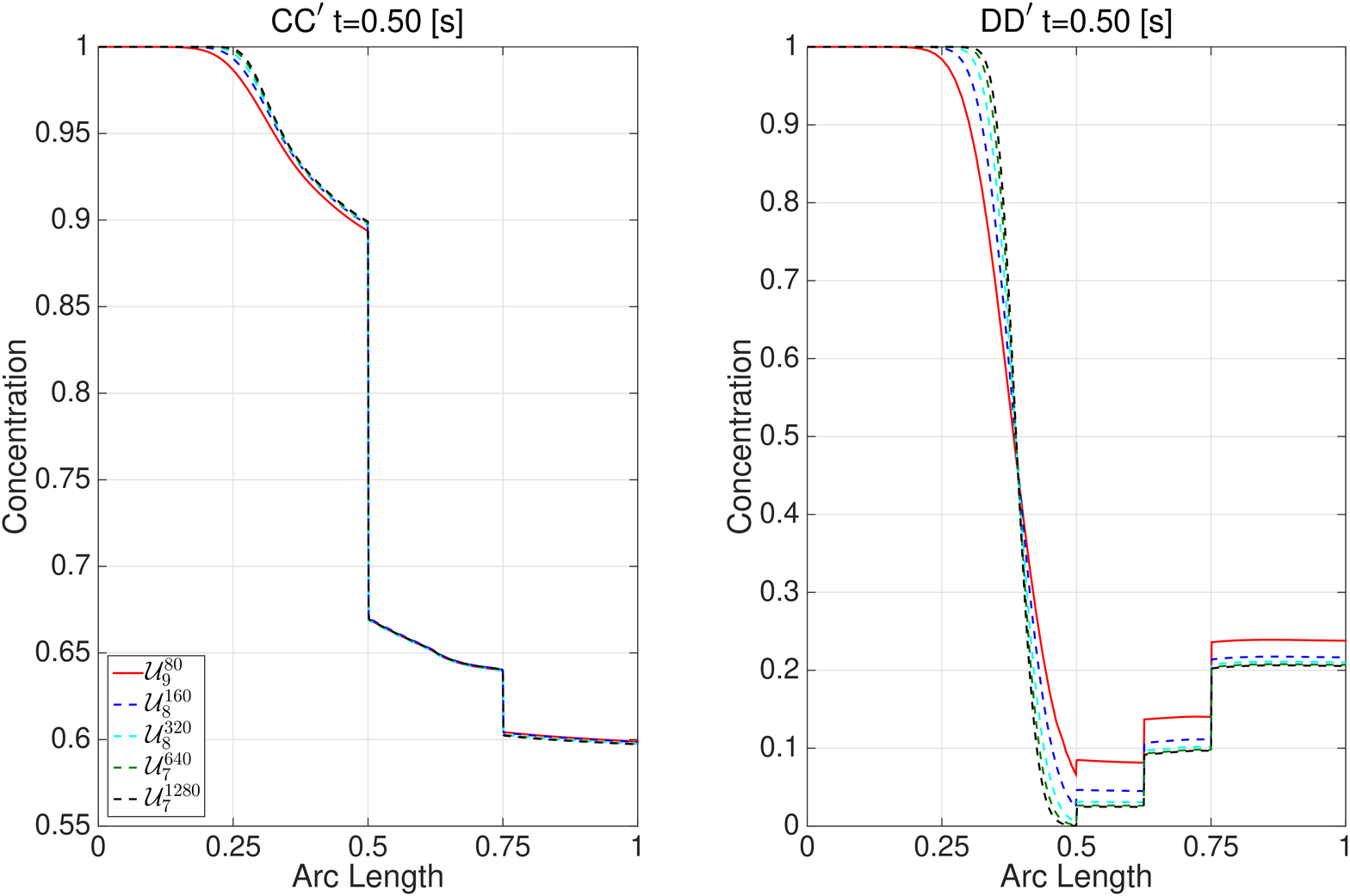}
\caption{Concentration values along two lines for the unresolved test cases.}
\label{fig:regular_c_unres}
\end{figure}

\begin{figure}[hbt!]
\centering
\includegraphics[width=13cm]{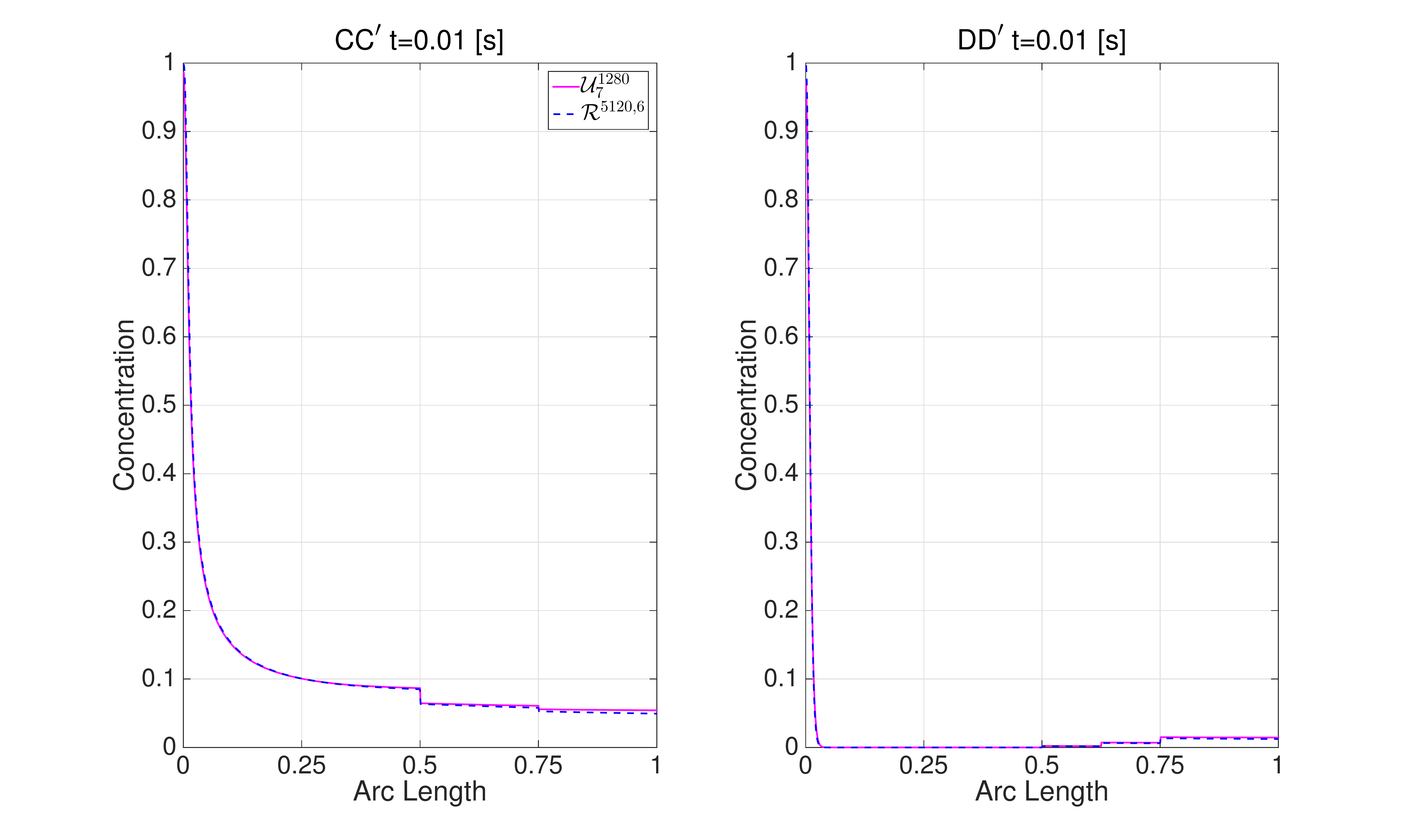}
\caption{Comparison of concentration values along two lines between resolved and unresolved test cases at time $t=0.01\,\text[s]$.}
\label{fig:regular_c_res_unres_1}
\end{figure}

\begin{figure}[hbt!]
\centering
\includegraphics[width=13cm]{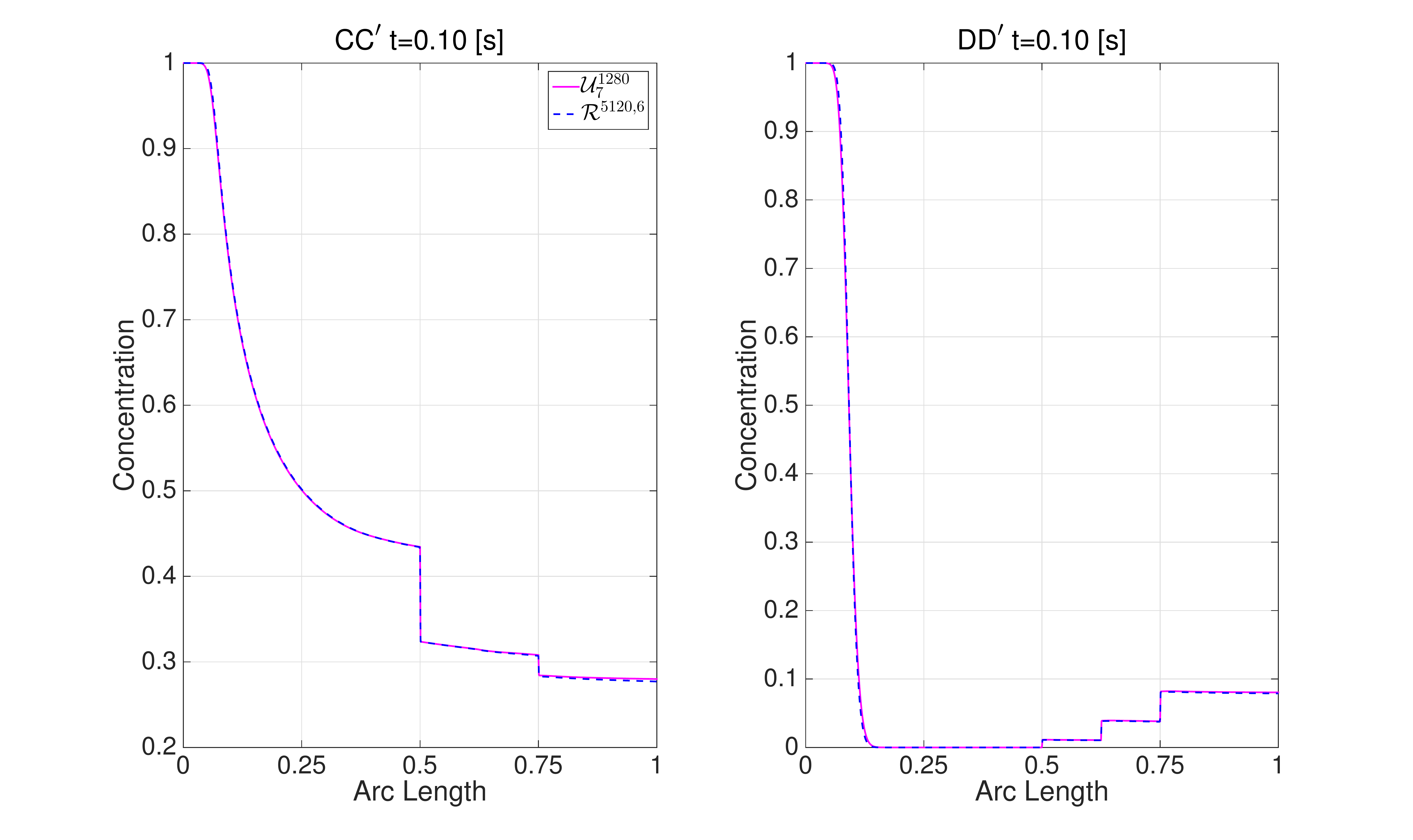}
\caption{Comparison of concentration values along two lines between resolved and unresolved test cases at time $t=0.1\,\text[s]$.}
\label{fig:regular_c_res_unres_2}
\end{figure}

\begin{figure}[hbt!]
\centering
\includegraphics[width=13cm]{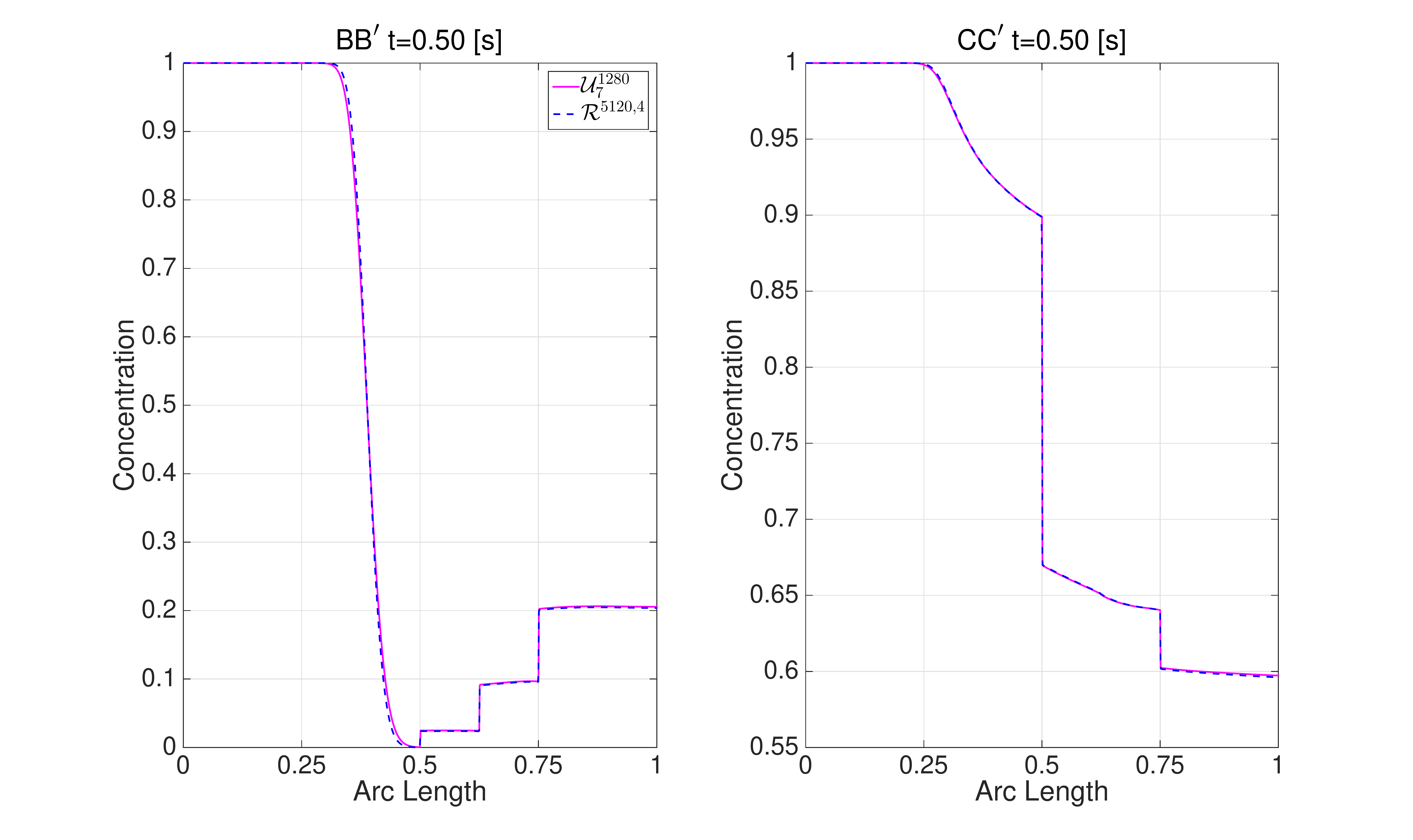}
\caption{Comparison of concentration values along two lines between resolved and unresolved test cases at time $t=0.5\,\text[s]$.}
\label{fig:regular_c_res_unres_3}
\end{figure}

\begin{figure}[hbt!]
\centering
\includegraphics[width=12.5cm]{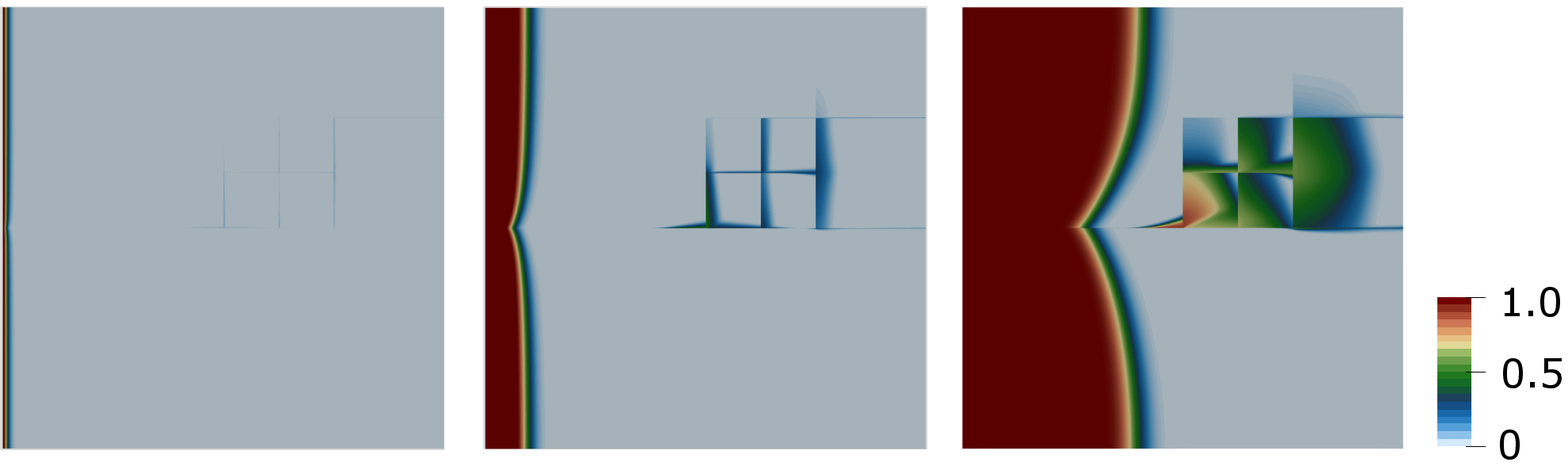} 
\caption{Concentration distribution at time  $t_1=0.01\,\text{[s]}$,  $t_2=0.10\,\text{[s]}$ , and $t=T_{\rm{fin}}=0.50\,\text{[s]}$ for the test case $\mathcal{U}_{7}^{1280}$.}
\label{fig:regular_DMP}
\end{figure}

For the transport problem, $\Gamma_D$, which coincides with $\Gamma_{in}$, is the left edge of $\domain$.
The boundary condition is set as $g=1.0$ [m$^{-3}$].
We consider an overall simulation time equal to $T_{\rm{fin}}=0.50~\rm{[s]}$, and set the time step equal to $\Delta t=0.025~\rm{[s]}$. 

In Figures~\ref{fig:regular_c_res} and~\ref{fig:regular_c_unres}, we report concentration profiles along the segments $\linec$ and $\lined$ for some resolved and unresolved meshes, respectively.
While no relevant differences are found for the unresolved meshes,
the numerical results computed on resolved meshes are very sensitive to the number of subdivisions of the background mesh $be$.
In particular, similar solutions are obtained for meshes with at least $be=4000$,
while all the meshes with lower resolution produce results with notable differences, especially along the segment $\linec$.
The results obtained with the unresolved meshes highlight again the importance of not only resolving the interfaces but also to refine the region close to them.
We point out that the high aspect ratio of the elements located inside the fracture aperture may affect the stability of the numerical method and reduce the accuracy of the approximated solutions.
Indeed, the numerical diffusion added in the direction of the longest side causes the unphysical undershoots visible where the segment $\linec$ crosses by the vertical fracture at $x=0.5$.

In Figures~\ref{fig:regular_c_res_unres_1},~\ref{fig:regular_c_res_unres_2}, and~\ref{fig:regular_c_res_unres_3},
we compare one of the finest resolved test cases with the numerical solution computed on the unresolved mesh with the largest number of element layers along the fracture aperture.
We analyze the concentration along lines $\linec$ and $\lined$ at different time instants: $t=0.01$, $t=0.1$ and $t=T_{\rm{fin}}=0.50$ [s] and observe similar solutions for both the types of meshes.
The good match found between the unresolved and the resolved results confirms that while the solution computed on the unresolved mesh is mostly affected by the number AMR,
the resolved test cases require a much finer background to improve the aspect ratio of the elements located inside the fracture.

In Figure~\ref{fig:regular_DMP}, we show the distribution of the concentration at the same times selected for the profiles depicted in Figure~\ref{fig:regular_c_res_unres_3}.
We observe that concentrations are negligible inside the embedding matrix and that the physical quantity is mainly transported inside the high permeable fracture network.


\subsection{Single 2D}

For the second example, we set $L_1=L_2 = 100\,[\text{m}]$.
Such a square domain is crossed by an oblique fracture with a thickness $\delta=0.01\,[\text{m}]$.
The subdomain $\domain_m$ consists of two subdomains, $\domain_{m1}$ and $\domain_{m2}$,
with different permeabilities, $k_{m1}$ and $k_{m2}$, and porosities ${\phi_{m1}}$ and ${\phi_{m2}}$.
We report details about the geometrical set up in Table~\ref{tab:geo_1}, where we denote by $x_i^{s}$ and $x_i^{e}$ the coordinates of starting and ending points, respectively.
The material properties are summarized in Table~\ref{tab:param_1}.

\begin{table}[ht]
\caption{Geometrical sizes employed in the problem \lq\lq Single fracture network\rq\rq. Here $s$ and $e$ refer to the starting and the end position, respectively.}
\centering
\begin{tabular}{ |p{1cm}||p{1cm}|p{1cm}|p{1cm}| |p{1cm}|  }
 \hline
 \multicolumn{5}{|c|}{Rock Matrix} \\
 \hline
      & {$x_{1}^s$} [m]& $x_{1}^e$ [m]& $x_{2}^s$ [m]& $ x_{2}^s$ [m] \\
 \hline
     \hline
$\mathrm\domain_{m1}$ &-50 & 50 & -50 & -40    \\
$\mathrm\domain_{m2}$  &-50 & 50 & -40 & -50  \\
$\mathrm{\Gamma_{in}}$  &0 & 0 & 40 & 50  \\
$\mathrm{\Gamma_{out}}$  &50 & 50 & 0 & 10  \\
$\mathrm\Omega_{A}$ &-50 & 0 & -50 & 50    \\
$\mathrm\Omega_{B}$  & 0 & 50 & -50 & 50  \\
 \hline
     \hline
 \multicolumn{5}{|c|}{Fracture} \\
  \hline
    \hline
     & {$x_{1}^s$} [m]& $x_{1}^e$ [m]& $x_{2}^s$ [m]& $ x_{2}^s$ [m] \\
 \hline
 $\domain_{f}$& -50   & 50 & -30 & 30\\
 \hline
     \hline
  \multicolumn{5}{|c|}{Segments} \\
  \hline
    \hline
      & {$x_{1}^s$} [m]& $x_{1}^e$ [m]& $x_{2}^s$ [m]& $ x_{2}^s$ [m] \\
 \hline
 $\linea$& -50   & 50 & 50 & -50\\
 $\lineb$& -50   & 50 & 30 & -30\\
 \hline
\end{tabular}
\label{tab:geo_1}
\end{table}

\begin{table}[ht]
\centering
\caption{Material properties employed in the problem \lq\lq Single fracture network\rq\rq.\label{tab:param_1}}
 \begin{tabular}{ | l | l | l | l |}
\hline
 Property & Symbol & Value & Unit \\
 \hline
  \hline
 Fracture thickness & $\delta$ & 0.01 & m\\
Matrix conductivity ($\domain_{m1}$) & $k_{m1}$    &  $10^{-5}$ & $\text{m}\times \text{s}$ \\
Matrix conductivity ($\domain_{m2}$)& $k_{m2}$    &  $10^{-6}$ &  $\text{m}\times \text{s}$\\
Fracture conductivity & $ k_{f}$     &  $10^{-1}$  &  $\text{m}\times \text{s}$\\
Matrix porosity ($\domain_{m1}$) &${\phi_{m1}}$    & $0.2 $ & \\
Matrix porosity ($\domain_{m2}$)& ${\phi_{m2}}$   & $0.25$ & \\
Fracture porosity & ${\phi_{f}}$      & $0.4$ &\\
 \hline
\end{tabular}
\end{table}

\begin{figure}[hbt!]
\centering
\subfloat[Resolved]{\includegraphics[width=5.98cm]{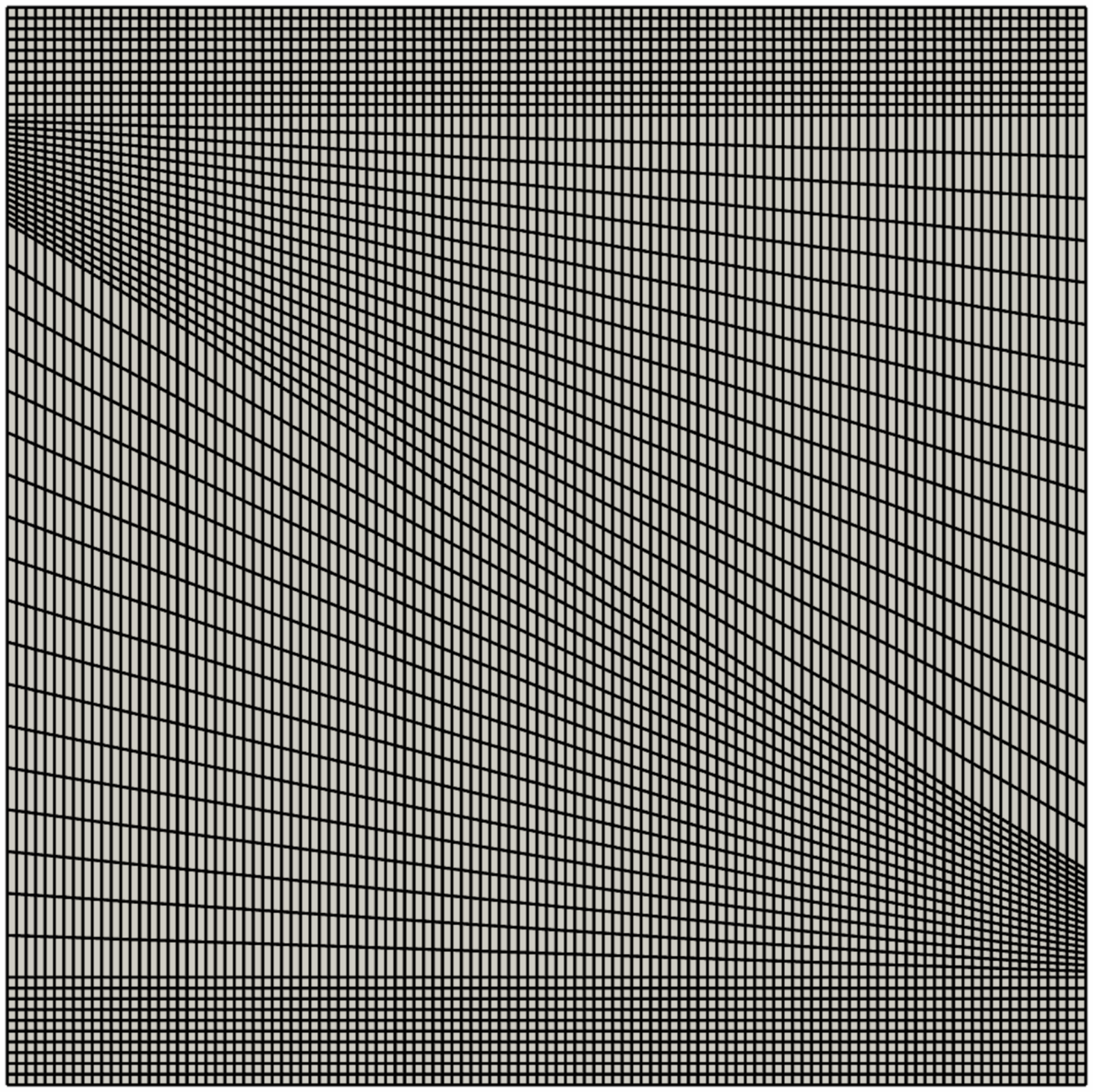}}
\subfloat[Unresolved]{\includegraphics[width=6.02cm]{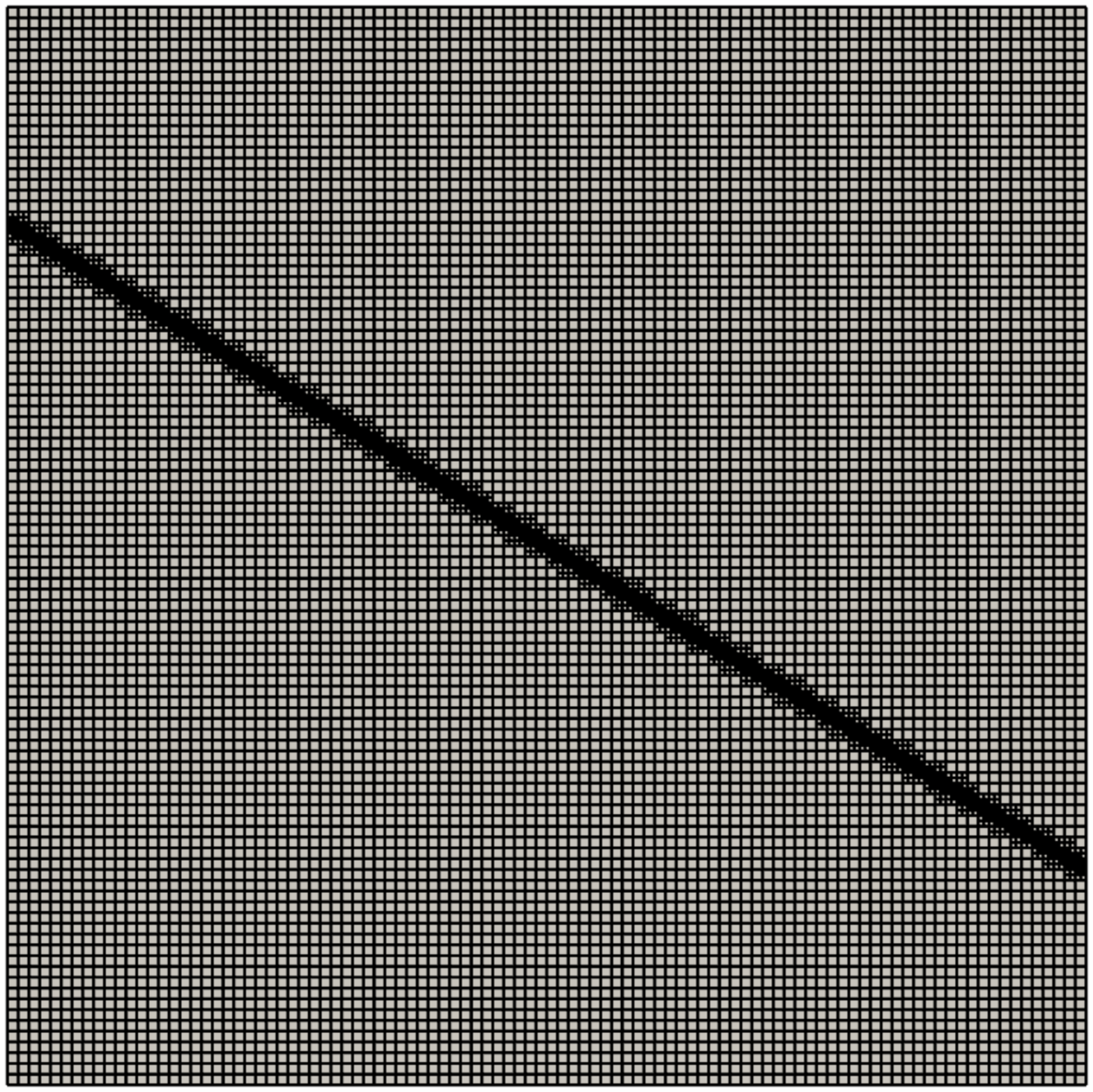}}
\caption{(a) Example of a resolved mesh ($\mathcal{R}^{114,2}$). (b) Example of an unresolved mesh ($\mathcal{U}^{100,8}$).}
\label{fig:single_unresolved_mesh}
\end{figure}

We compare numerical results computed on several resolved and unresolved meshes.
The characteristics of the considered meshes are reported in Table~\ref{tab:mesh_single_res} and Table~\ref{tab:mesh_single_unres}, respectively.
In Figure~\ref{fig:single_unresolved_mesh}, we report the resolved mesh $\mathcal{R}^{114,2}$ and the unresolved mesh $\mathcal{T}^{100,8}$. 
As the fracture is oblique, mesh elements of the resolved mesh are skewed,
while the unresolved meshes allow to keep the same aspect ratio for all the elements.
We analyze pressure and concentrations along the segments $\linea$ and $\lineb$ depicted in Figure~\ref{fig:single_p_unresolved_setup}, and whose coordinates are reported in Table~\ref{tab:geo_1}. For the flow problem, we estimate the approximated flux across the line $L:=\linec$ located at $x_1=0$ [m].

\begin{table}[ht]
\caption{Total fluxes, Equation~\ref{eq:tot_flux}, for the resolved and unresolved test cases. All the fluxes values have been multiplied by a scale factor $1.0e^{6}$.}
\centering
\begin{tabular}{ |p{1cm}|p{2cm}|p{2cm}|  }
 \hline
      Mesh & $Q_{\linec} (\domain_1)$& $Q_{\linec} (\domain_2)$\\
 \hline
   \hline
$\mathcal{R}^{114,2}$&$1.86675$  & $-1.86675$\\ 
$\mathcal{R}^{228,4}$&$1.85977$  & $-1.85977$  \\
$\mathcal{R}^{456,8}$&$1.85632$  & $-1.85632$ \\
$\mathcal{R}^{921,16}$&$1.85457$  & $-1.85457$\\
 \hline
 $\mathcal{U}^{100}_8$&$1.87065$         & $-1.87065$ \\
 $\mathcal{U}^{200}_6$ & $1.86277$ & $-1.86277$\\
 $\mathcal{U}^{200}_7$ & $1.86277$ & $-1.86277$\\
 $\mathcal{U}^{200}_8$& $1.86174$ & $-1.86174$\\
 $\mathcal{U}^{400}_6$ & $1.85779$ & $-1.85779$\\
 $\mathcal{U}^{400}_7$ & $1.85747$ & $-1.85747$  \\
 $\mathcal{U}^{400}_8$&$1.85534$ & $-1.85534$ \\
 $\mathcal{U}^{800}_6$&$1.85534$ & $-1.85534$ \\
 $\mathcal{U}^{800}_7$ & $1.85517$ & $-1.85517$ \\
  $\mathcal{U}^{800}_8$&$ 1.85505$ & $-1.85505$ \\
  $\mathcal{U}^{1600}_6$& $ 1.85411$ & $-1.85411$ \\ 
 $\mathcal{U}^{1600}_7$& $ 1.85399$ & $-1.85399$ \\  

 \hline
\end{tabular}
\label{tab:flux_single_res_unres}
\end{table}

\subsubsection{Flow problem and local conservation of the flux}

\begin{figure}[hbt!]
\centering
\includegraphics[width=8cm ]{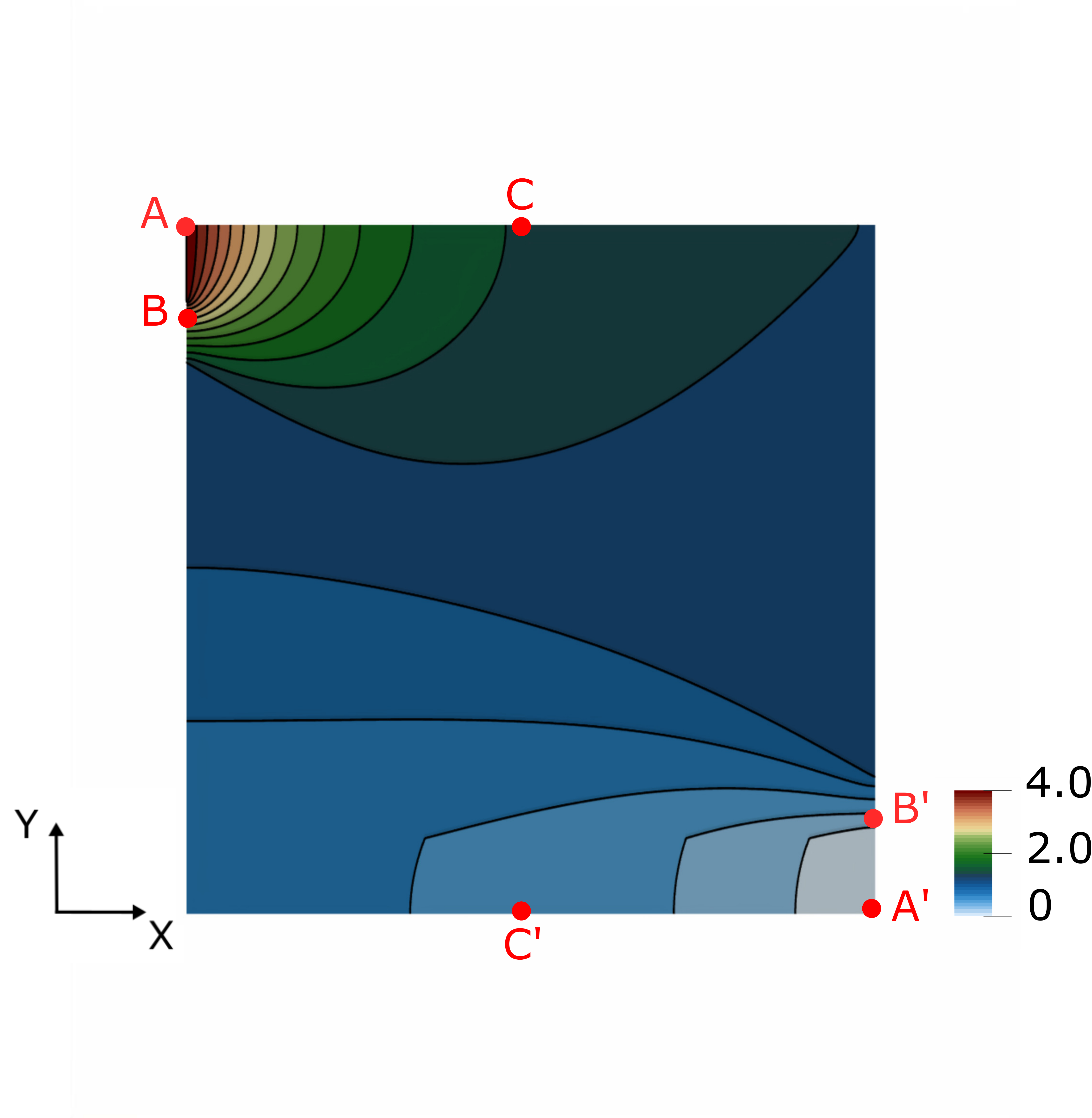}
\caption{Pressure distribution for the unresolved test case $\mathcal{U}_8^{320}$. We also report the extrema of the three segments along which the relevant properties are computed.}
\label{fig:single_p_unresolved_setup}
\end{figure}

For the flow problem, the Dirichlet boundary $\Gamma_D$ consists of two segments $\Gamma_{in}$ and $\Gamma_{out}$, as reported in Table~\ref{tab:geo_1},
while $\Gamma_N$ is composed of the other remaining sides.
We set $g_{in}=4$ [m] and $g_{out}=1$ [m] on $\Gamma_{in}$ and $\Gamma_{out}$, respectively.
Homogeneous Neumann conditions are imposed on $\Gamma_N$. Figure~\ref{fig:single_p_unresolved_setup} shows the spatial distribution of the pressure obtained for the test case $\mathcal{U}^{400}_6$.

\begin{figure}[hbt!]
\centering
\includegraphics[width=12cm]{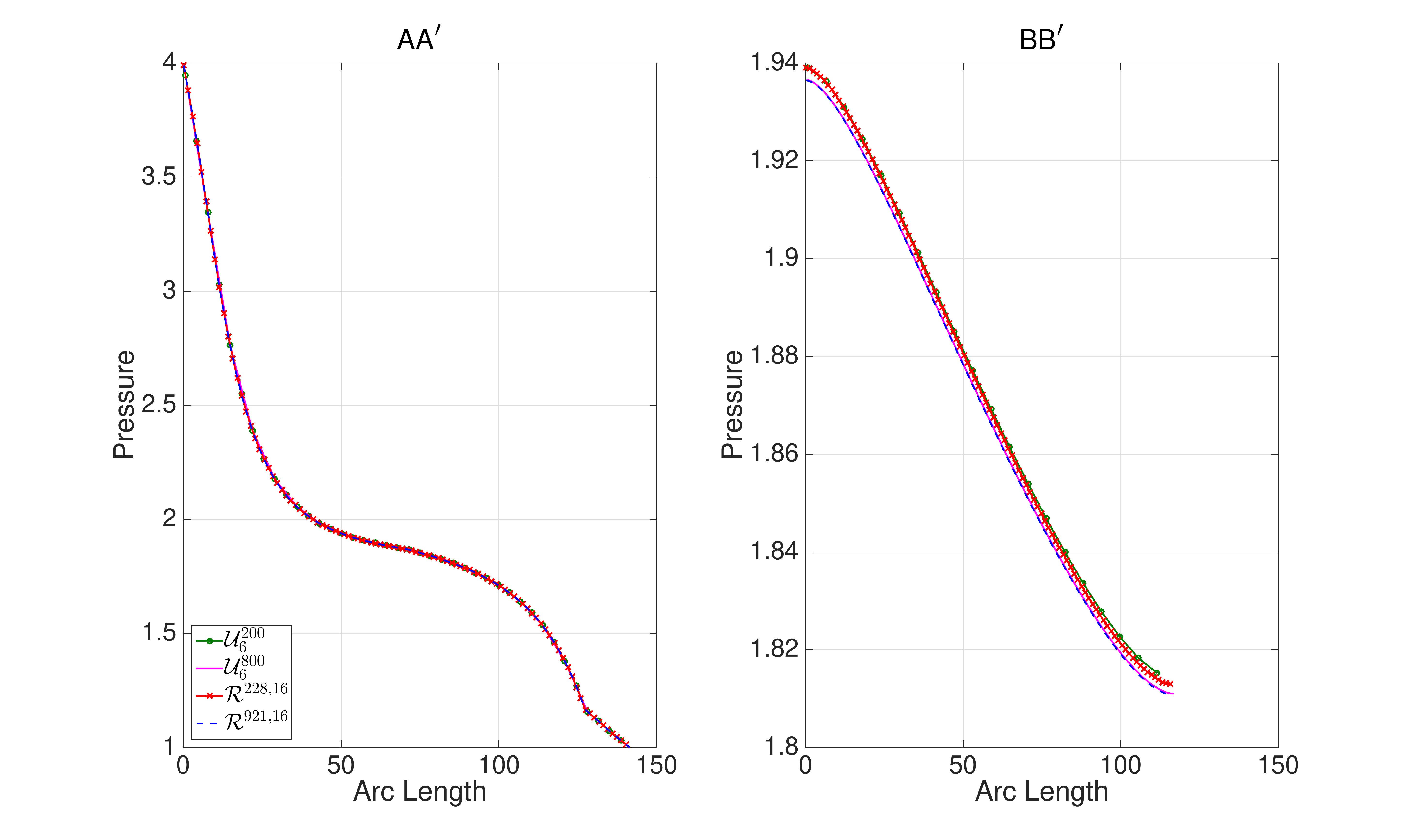}
\caption{Comparison of pressure values along two lines between resolved and unresolved test cases.}
\label{fig:figure_p_res_unres_bench_1}
\end{figure}

In Figure~\ref{fig:figure_p_res_unres_bench_1}, we compare the pressure computed along the lines $\linea$ and $\lineb$ on some resolved and unresolved meshes. In particular, we compare coarser mesh solutions with results computed on meshes fine enough to guarantee mesh-independence.
We observe that solutions present a similar behavior along the segment $\linea$, while the results computed along the line $\lineb$ are more sensitive to the mesh resolution. 
As already observed for the benchmark \lq\lq Regular fracture network\rq\rq, the resolved meshes require a certain number of elements in the background, while the solution computed on the unresolved meshes is mainly affected by the number of elements close to the matrix-fracture interface.

\begin{figure}[hbt!]
\centering
\includegraphics[width=12cm]{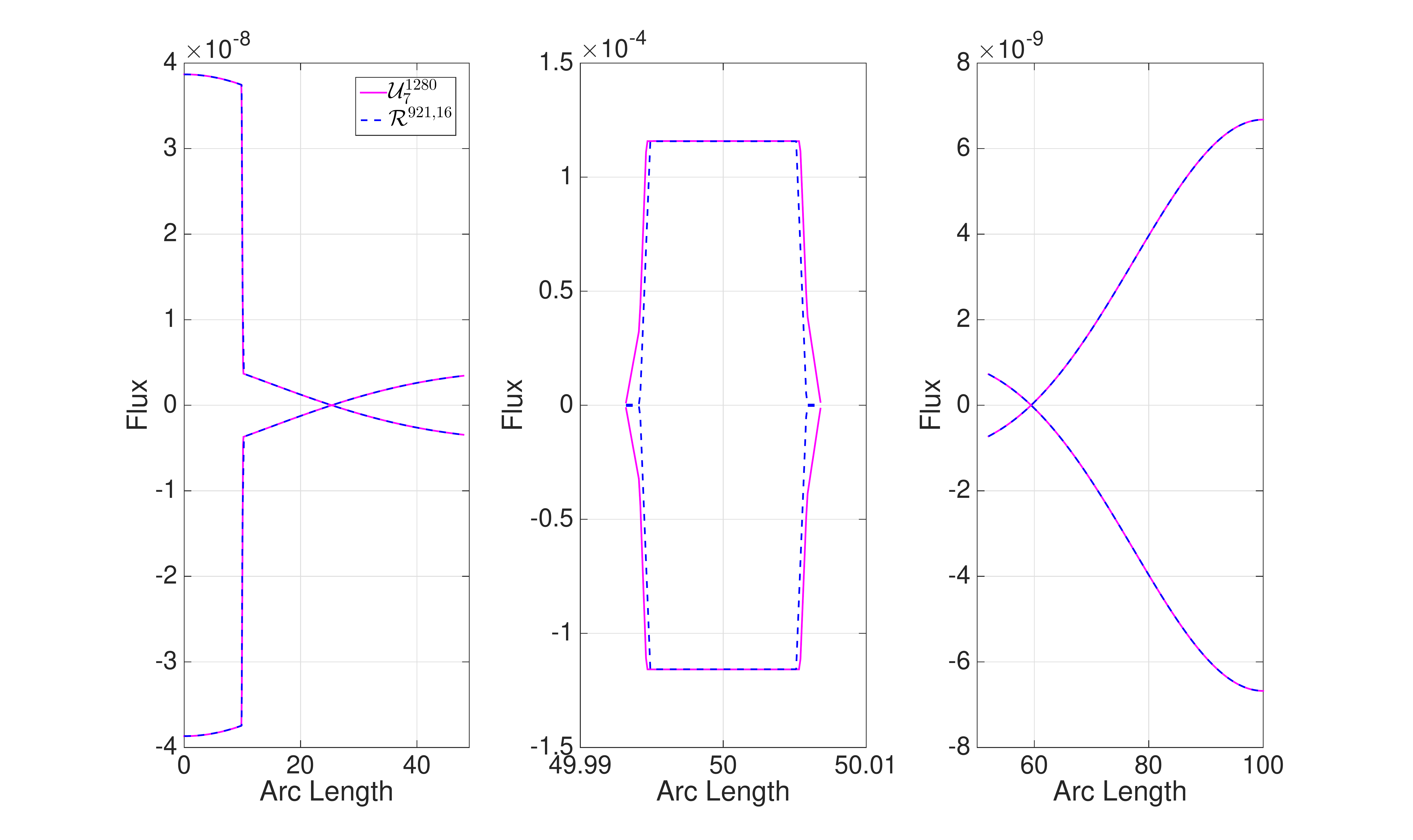}
\caption{Boundary fluxes, $q_{\linec}(\Omega_1)$  and $q_{\linec}(\Omega_2)$, across boundary $\linec$ for the test cases $\mathcal{R}^{921,16}$ and $\mathcal{U}^{1600}_7$.}
\label{fig:single_flux}
\end{figure}

\begin{figure}[hbt!]
\centering
\includegraphics[width=8cm]{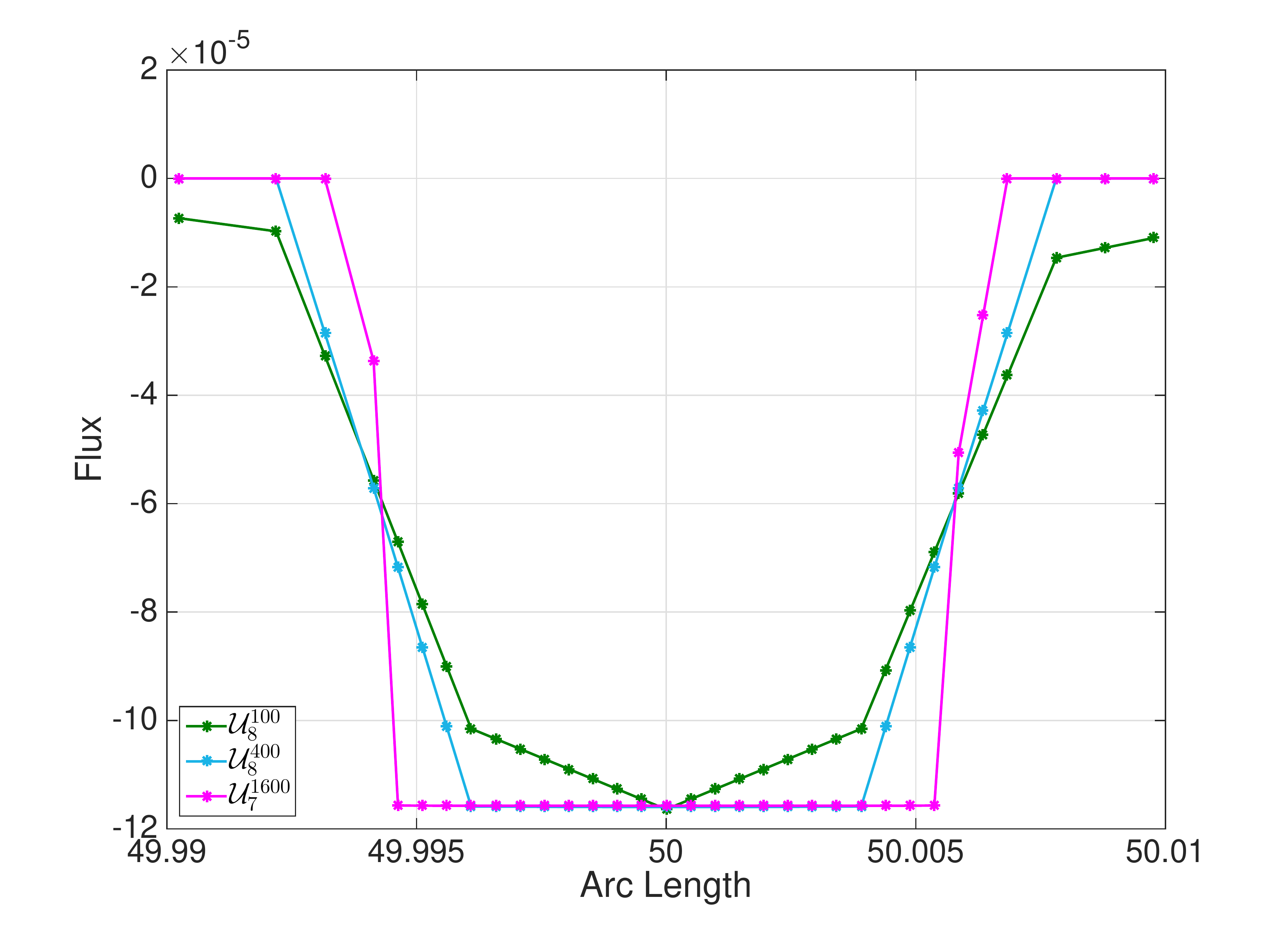}
\caption{Boundary flux $q_{\linec}(\Omega_2)$ across the region of the boundary $\linec$ intersecting the fracture for the test cases $\mathcal{U}^{100,8}$, $\mathcal{U}^{200,8}$ and $\mathcal{U}^{1600}_7$.}
\label{fig:single_flux_zoom_in}
\end{figure}

In Figure~\ref{fig:single_flux}, we report the flux along the segment $\linec$ computed on the meshes $\mathcal{R}^{921,16}$ and $\mathcal{U}^{1600}_7$.
Results for both types of meshes present the same profile but with opposite sign, showing that the two fluxes are matched point-wise and not only in an integral form.
On the other hand, differences between the results obtained with the two types of mesh are visible close to the two sides of the fracture.
The relative position between the segment $\linec$ and the oblique fracture require a large number of either mesh subdivisions,
$be$, or adaptive refinement steps, AMR, to reproduce the steep profile computed for the resolved test cases. 
The results reported in Figure~\ref{fig:single_flux_zoom_in}) for some unresolved meshes confirm this observation.
Indeed, by progressively increasing the resolution in the fracture, the flux $q_{\linec}(\domain_2)$ in the central part of the segment $\linec$ tends to correctly reproduce the steep characteristics of the flux.

For the unresolved test cases, we estimate the error by adopting the flux computed on $\mathcal{U}^{1600}_7$ as a reference solution.
The results reported in Table~\ref{tab:flux_error_single} show that, for a fixed number of AMR steps,
the error $e_q$ progressively reduces by increasing the number of mesh subdivisions $be$.
In particular, the estimated converge rate presents an order larger than $1.5$.
Similar behavior can be observed fixing the number of mesh subdivisions and analyzing the error behavior with respect to the number of AMR steps.

In Table~\ref{tab:flux_single_res_unres} we report the value of the total fluxes $Q_{L}(S)$ with $S \in \{ \domain_1,\domain_2\}$ computed for meshes listed in Tables~\ref{tab:mesh_single_res} and~\ref{tab:mesh_single_unres}.
Again, all the numerical results confirm that the local conservation is satisfied for both resolved and unresolved meshes. 

\begin{table}[ht]
\caption{$L_2$ relative error,  Equation~\ref{eq:error}, for the boundary flux $q_{\linec}(\domain_1)$ computed on the unresolved meshes.}
\begin{tabular}{ |p{1.8cm}|p{3cm}|p{3cm}|p{3cm}| }
  \hline
   \diagbox{be:}{amr:}
   & 6 & 7 & 8\\
 \hline
 \hline
 $100$  & - -& - -& 0.004862\\
 $200$  & 0.009018  &  0.004894 & 0.003498\\
 $400$  & 0.004911  & 0.003508 &  0.002235\\
 $800$   &0.003512  &  0.002238 &  0.001351\\
 $1600$ & 0.002239 & - - & - -\\
 \hline
\end{tabular}
\label{tab:flux_error_single}
\end{table}

\subsubsection*{Transport problem and DMP}

\begin{figure}[hbt!]
\centering
\includegraphics[width=12cm]{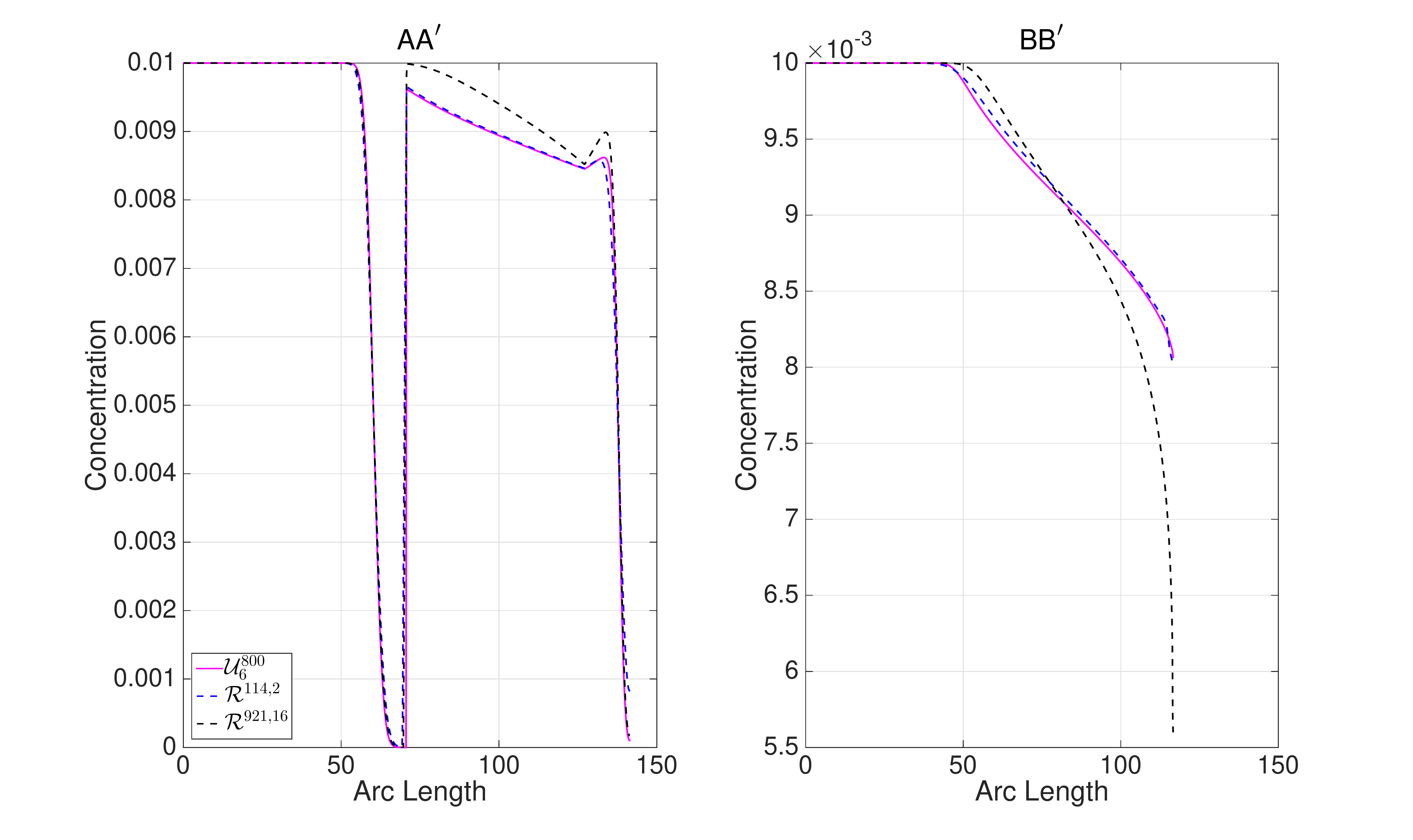}
\caption{Comparison of concentration values along two lines between two selected resolved and unresolved test cases at time $t=T=10^{9}\,\text{[s]}$.}
\label{fig:single_CM}
\end{figure}

\begin{figure}[hbt!]
\centering
\includegraphics[width=13cm]{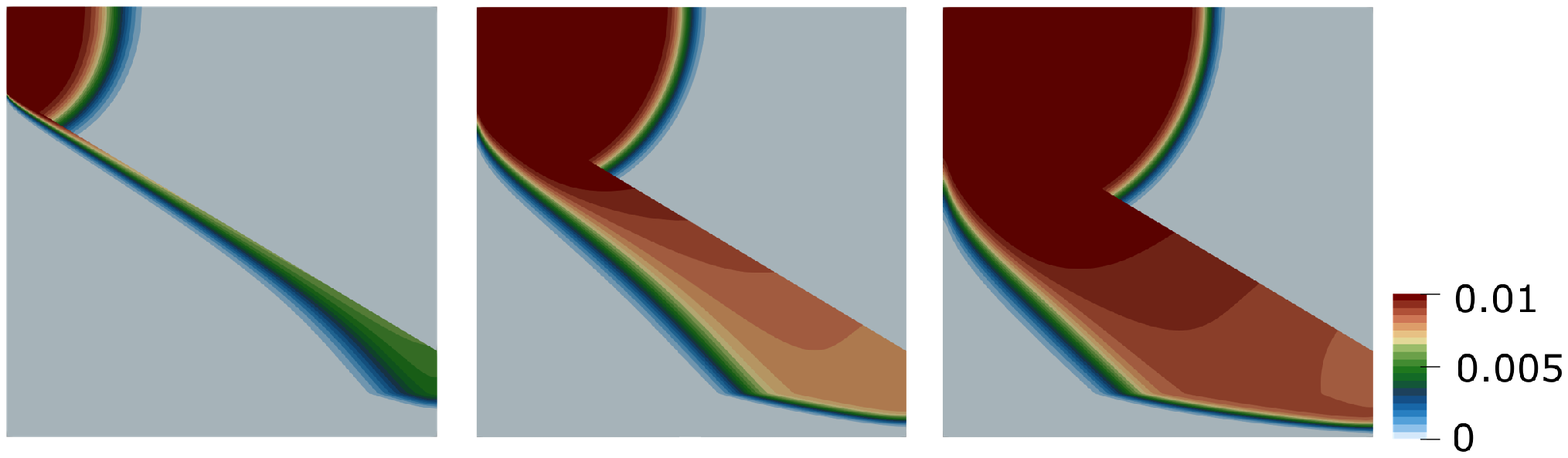}
\caption{Concentration distribution at time  $t_1=10^{8}\,[s]$,  $t_2=5\,10\,^{8}[s]$ , and $t=T=\,10^{9}\,[s]$ for the test case $\mathcal{U}^{800}_6$.}
\label{fig:unresolved_single_c}
\end{figure}

For the transport problem, $\Gamma_D$ coincides with $\Gamma_{in}$, on which we impose $g= 0.01$ [m$^{-3}$].
The rest of the boundary is $\Gamma_{out}$.
We set $T_{\rm{fin}}=10^{9}$ [s] and $\Delta t=10^2$ [s].
We compare the concentrations obtained on some resolved and unresolved meshes by analyzing their value at the final time $t=T_{\rm{fin}}$ along the lines $\linea$ and $\lineb$.
The results are reported in Figure~\ref{fig:single_CM}.

Although for the flow problem results coincide between the two types of meshes,
the concentration profiles computed on $\mathcal{R}^{921,16}$ and $\mathcal{U}^{800}_6$ present different behaviours close to the point $\text{A}^{\prime}$ inside $\domain_m$ and close to the point $\text{B}^{\prime}$ inside $\domain_f$.
Such differences are related to the skewness of the elements of the resolved meshes which introduce a larger stabilization to preserve the positivity of the solution but
reduce its accuracy.
Moreover, the aspect ratio of the bilinear elements deteriorates when increasing the number of subdivision $be$.
This observation is supported by the good agreement found between the concentrations computed on the unresolved meshes and the one computed on $\mathcal{R}^{114,2}$.

In Figure~\ref{fig:unresolved_single_c}, we show the distribution of the concentration $c$ at the time instants $t=10^{8}$, $t=5\,10\,^{8}$ , and $t=T_{\rm{fin}}=\,10^{9}$ [s]. 
Here, the values of the concentrations are always comprised between $0$ and $0.01$ for the entire simulated time, and verify the ability of the AFC scheme to ensure the DMP and to preserve the positivity of the solution.


\subsection{Realistic fracture network}

This example has been presented in~\cite{flemisch2018benchmarks}, where a comparison between different discretizations based on a one-dimensional representation of the fractures has been performed.
For this numerical example,
we have $L_1=700$ [m] and $L_2=600$ [m] and the set $\mathcal F$ is composed of $63$ fractures, ranging from isolated fractures up to groups of tens of fractures each.
The aperture of the fractures is $10^{-2}$ [m].
Indeed, such a geometry represents an outcrop in the Sotra island, near Bergen in Norway.
The matrix permeability is $k_m=10^{-14}$ [m$\times$s$^{-1}$], while the fracture permeability is to $k_f=10^{-8}$ m$\times$s$^{-1}$.
We consider no-flow boundary condition on top and bottom, pressure $1013250$ [m] on the left, and pressure $0$ [m] on the right of the boundary of the domain.
Due to the high geometrical complexity of the fracture network, meshes which resolve the interface $\Gamma$ are difficult to be generated.
On the other hand, AMR method enables a fully automatic generation of meshes for domains that contain a large number of highly complex heterogeneities.

As we have discussed in Section~2 and we will show in the numerical results,
the DMP may not hold for FE discretizations of diffusion problems on adapted meshes.
We will show that the proposed stabilization instead ensures the DMP and provides a suitable discretization method.

The characteristics of the considered meshes are reported in Table~\ref{tab:meshcharacteristicsrealsitc}.
The superscript $\text{be}$ refers to the number of subdivisions along the horizontal direction.
The number of intervals along the vertical direction is $6\text{be}/7$.

\begin{table}
\centering
\caption{Maximum and the minimum values for the pressure for the benchmark   \lq\lq Realistic fracture networks\rq\rq.}
\begin{tabular}{| l | r | r | r | }
\hline
Mesh & $\min \hh{p}$ & $\max \hh{p}$\\
\hline
\hline
$\mathcal U^{7}_{7}$&-3255.8181&1013250\\
\hline
$\mathcal U^{7}_{8}$&-2.7606e-14&1014210.2276\\
\hline
$\mathcal U^{7}_{9}$&-100.4391&1013250\\
\hline
$\mathcal U^{7}_{10}$&-10.6924&1013250\\
\hline
$\mathcal U^{14}_{6}$&-3255.8181&1013250\\
\hline
$\mathcal U^{14}_{7}$&-2.7606e-14&1014210.2276\\
\hline
$\mathcal U^{14}_{8}$&-100.4391&1013250\\
\hline
$\mathcal U^{14}_{9}$&-10.6924&1013250\\
\hline
$\mathcal U^{28}_{5}$&-3255.8451&1013250\\
\hline
$\mathcal U^{28}_{6}$&-2.7607e-14&1014210.112\\
\hline
$\mathcal U^{28}_{7}$&-100.4418&1013250\\
\hline
$\mathcal U^{28}_{8}$&-10.6928&1013250\\
\hline
\end{tabular}
\label{tab:pressure_realistic}
\end{table}

\begin{figure}
\centering
\subfloat[Step 0]{\includegraphics[width= 2in]{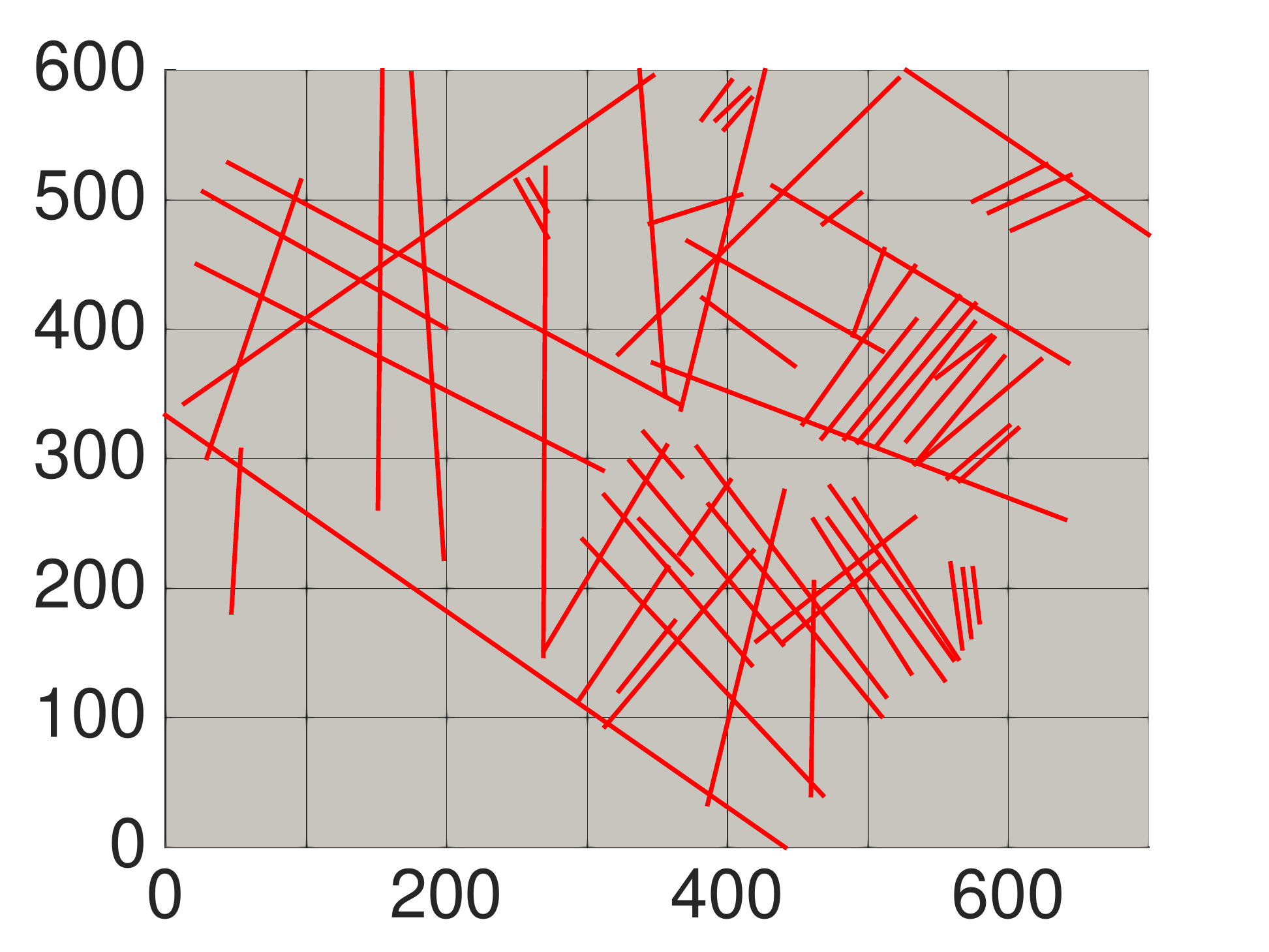}}
\subfloat[Step 1]{\includegraphics[width= 2in]{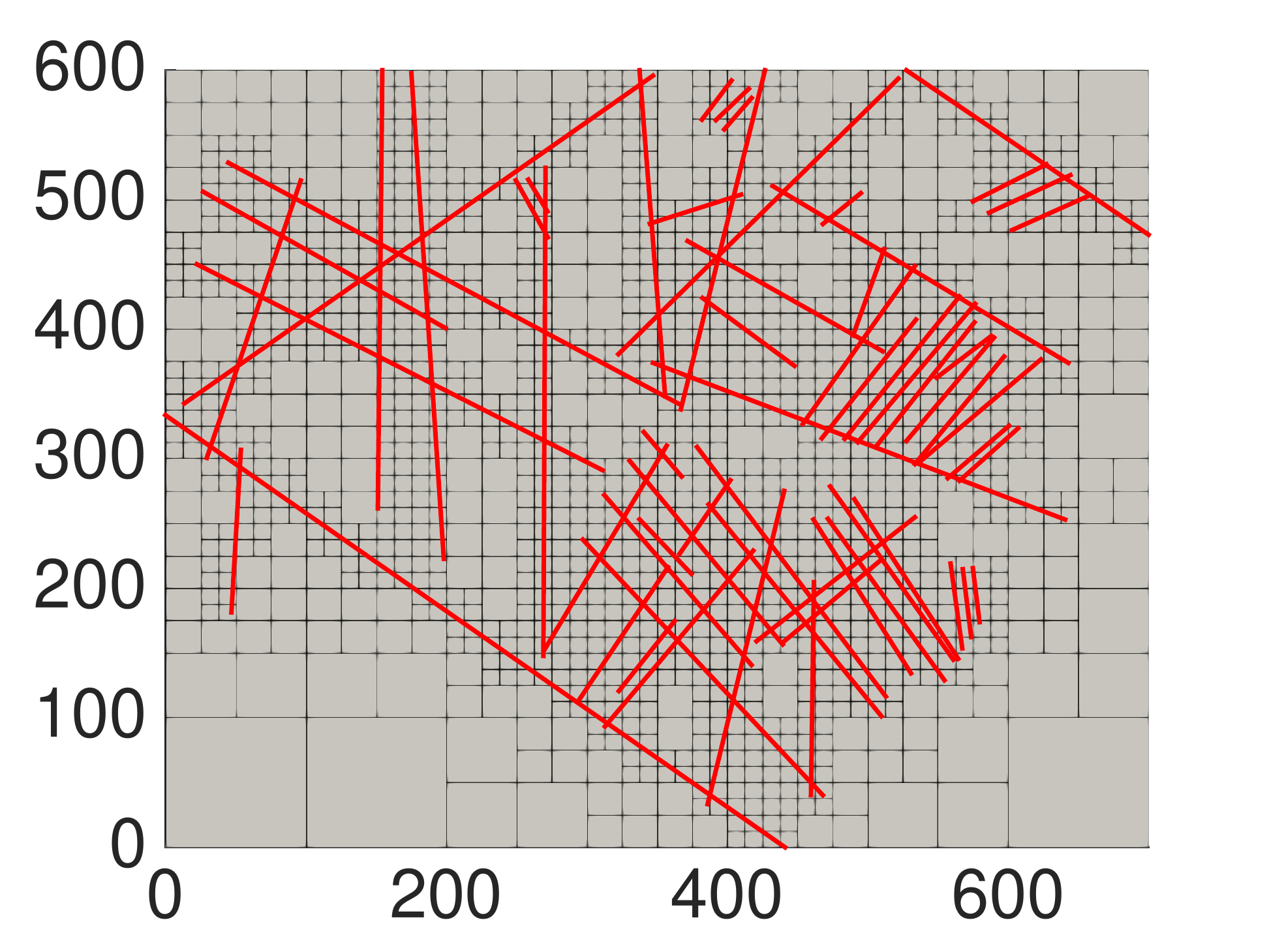}}\\
\subfloat[Step 2]{\includegraphics[width= 2in]{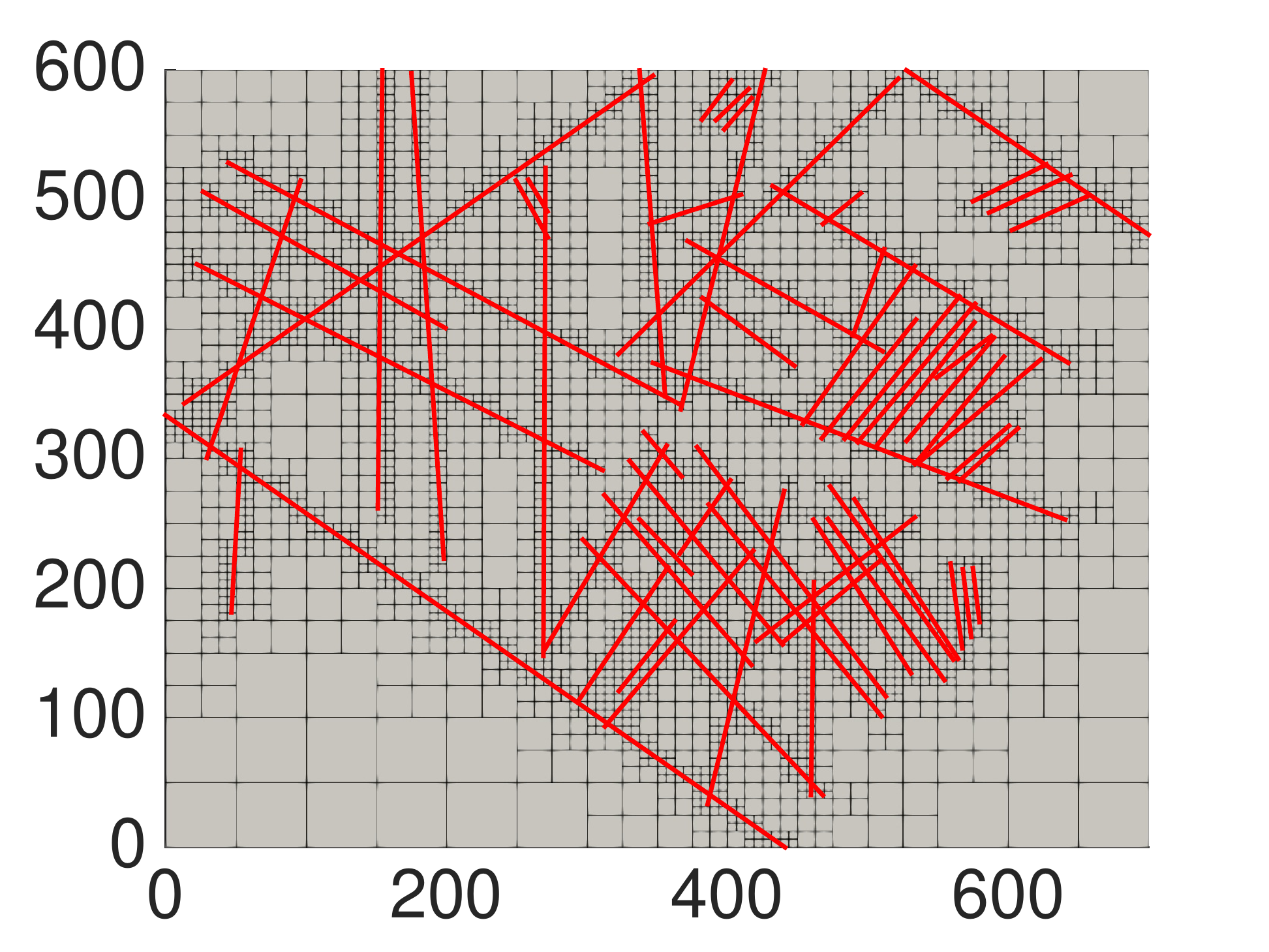}}
\subfloat[Step 4]{\includegraphics[width= 2in]{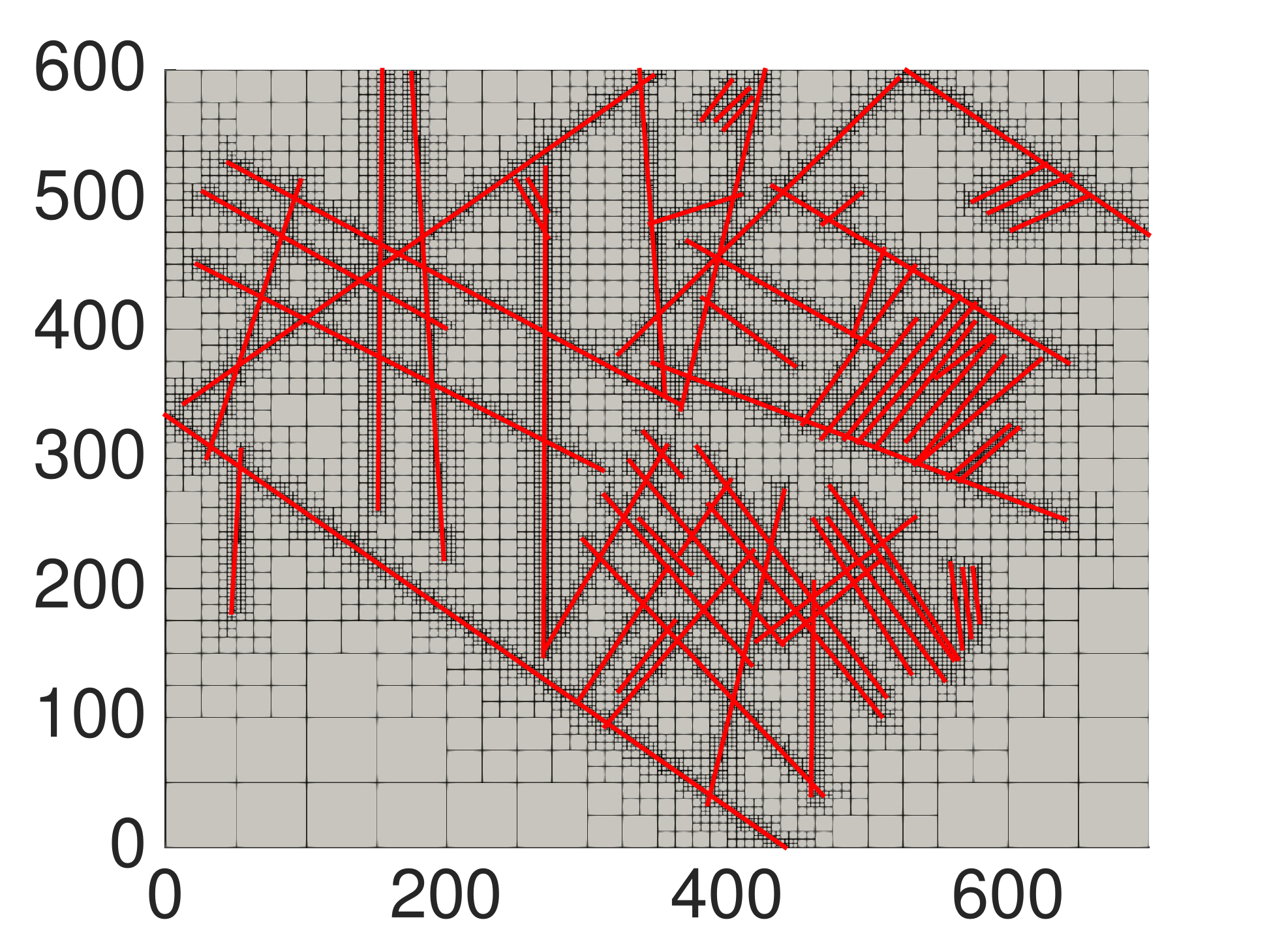}}
\caption{Sequences of meshes generated by iterating several times the AMR procedure.}
\label{fig:realisticadaptive}
\end{figure}

\begin{figure}
\centering
\subfloat[Unresolved Non Stabilized]{\includegraphics[width=3.0in]{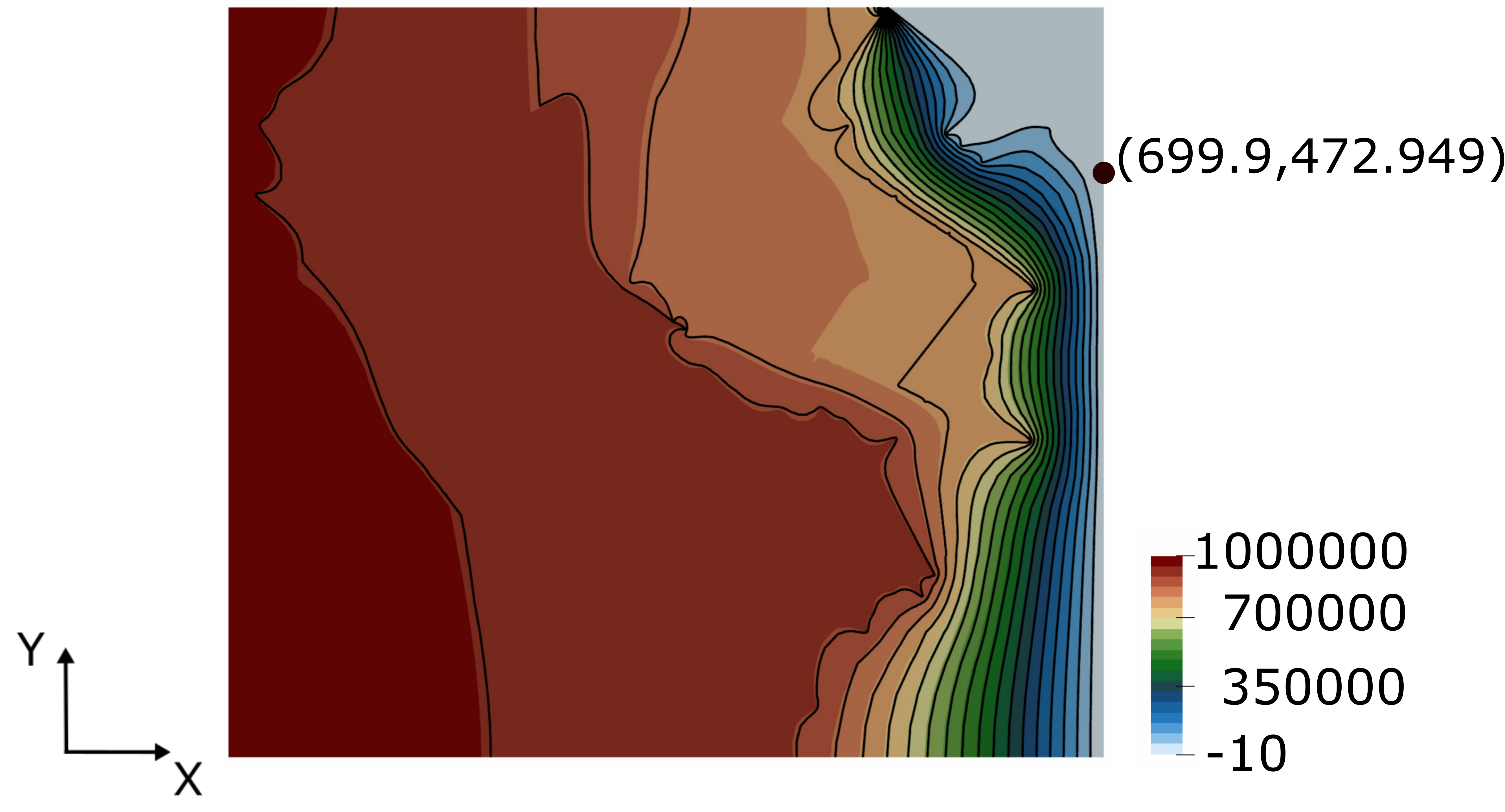}}
\subfloat[Unresolved Stabilized]{\includegraphics[width=2.44in]{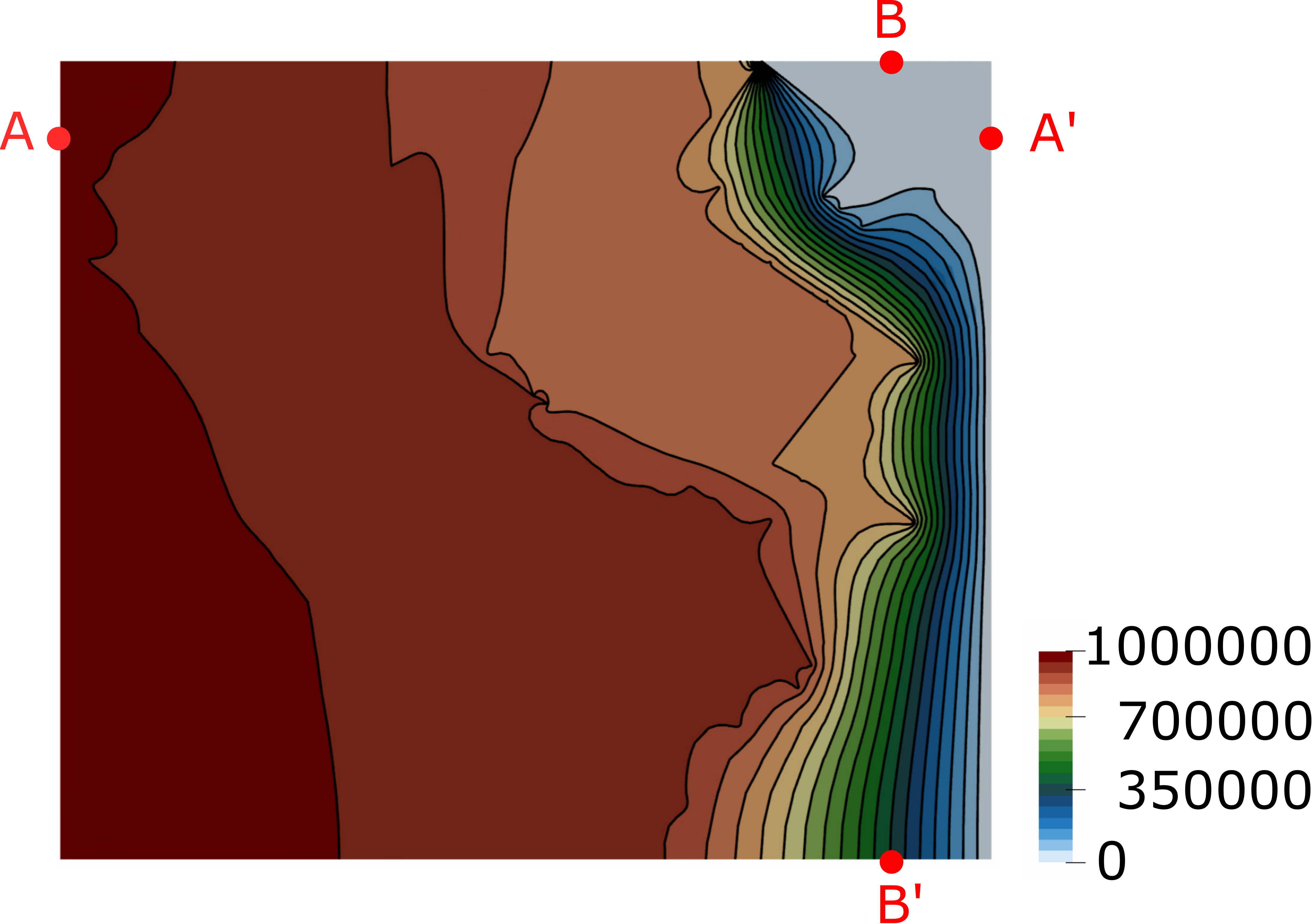}}
\caption{Pressure distribution for the unresolved test case $\mathcal U^{28}_8$. We also report the extrema of the three segments along which the relevant properties are computed.}
\label{fig:realisticsolution}
\end{figure}

In Table~\ref{tab:pressure_realistic}, we report the minimum and the maximum
for a non-stabilized FE discretization implemented on the corresponding grid.
We observe that for all the cases, the DMP is violated,
as the solution provided either a minimum smaller than $0$ or a maximum larger than the value imposed on the left side of the domain.
On the other hand, the stabilized discretization provides a solution that attains its maximum and minimum value at the boundary.
In Figure~\ref{fig:realisticsolution}, we report the solutions obtained with the stabilized and non-stabilized methods on the mesh $\mathcal U^{28}_8$.
On the non-stabilized solution (a), we superimpose the contour lines of the stabilized solution (b)  to highlight the differences between the two solutions.
We observe that the largest differences are found in the central part of the domain.

In Table~\ref{tab:realisticerrors}, we report the relative $L^2$ errors (i.e., replacing $q_L(\domain_i)$ with $p^h$ in Equation~\ref{eq:error}) for the two methods with respect to a reference solution computed on a mesh $\mathcal U^{56}_{11}$ with a non-stabilized method.
The magnitude of the errors computed with the stabilized and non-stabilized methods are comparable,
with the stabilized one providing, in general, a bit larger errors but ensuring the DMP.
The addition of the algebraic stabilization is in general non-consistent,
unless non-linear stabilization methods are considered~\cite{KUZMIN20093448}.

\begin{figure}
{\includegraphics[width= 12cm]{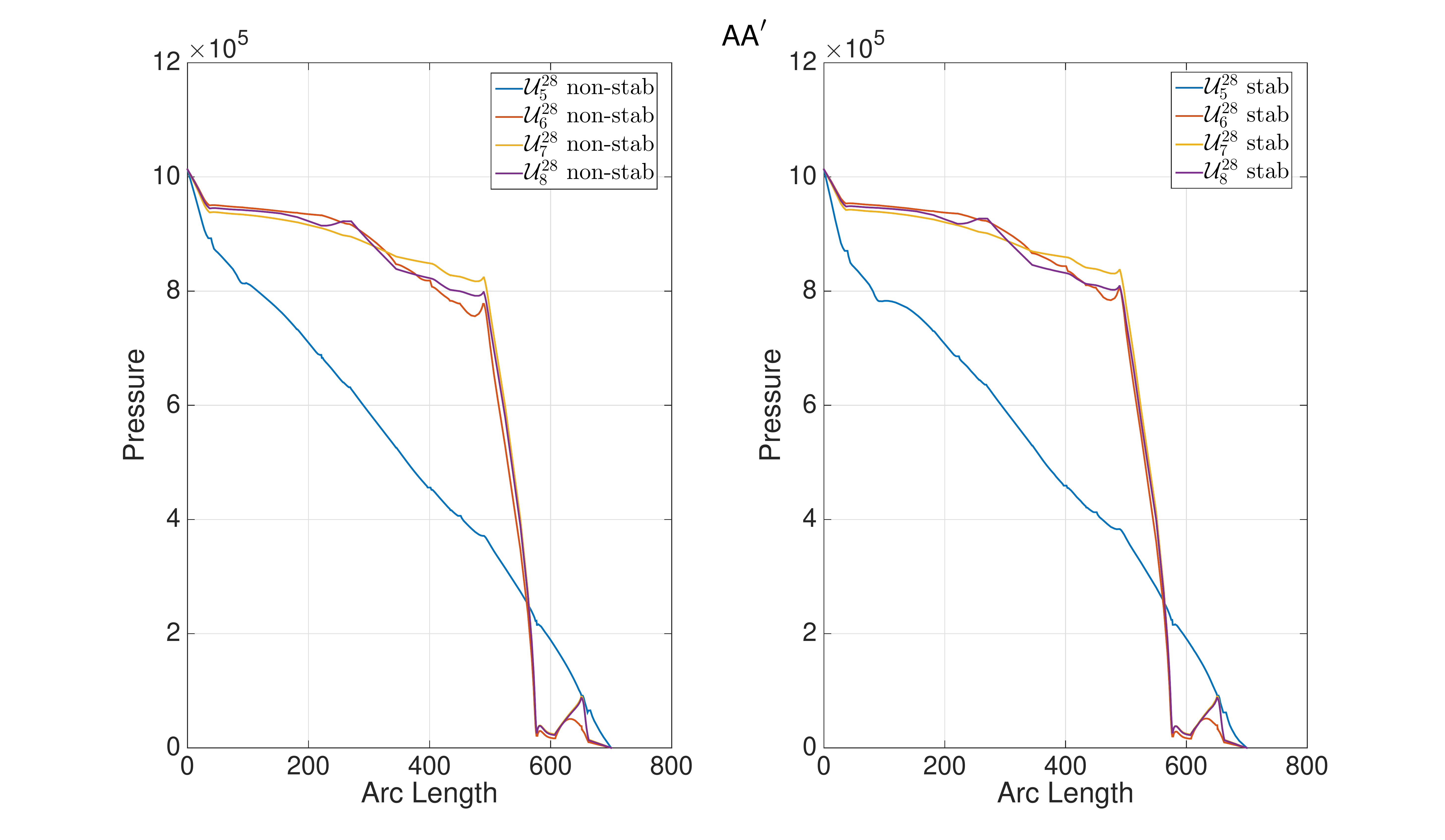}}\\
{\includegraphics[width= 12cm]{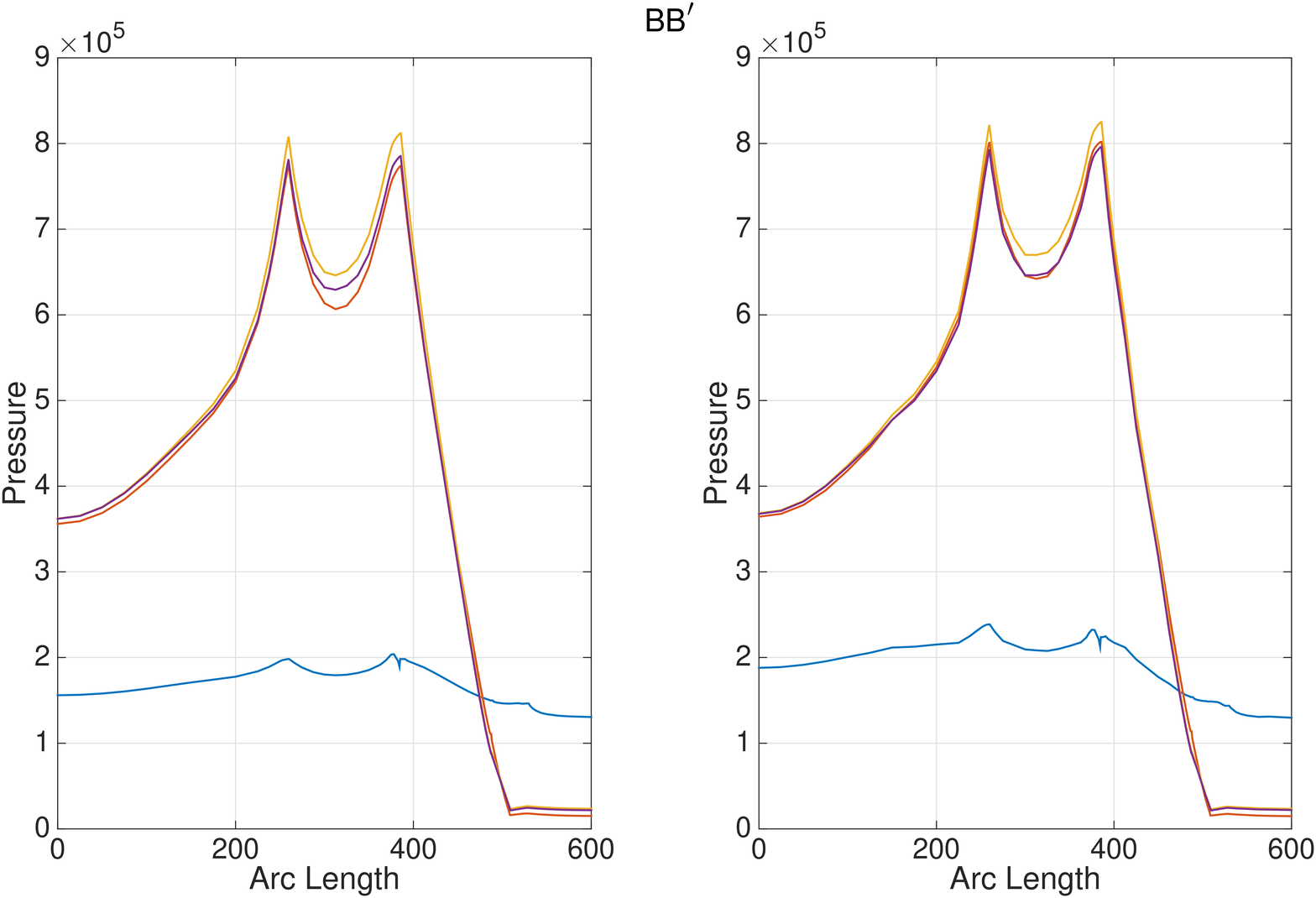}}
\caption{Comparison of pressure values along two lines between stabilized and non stabilized unresolved test cases. }
\label{fig:realisticalongline}
\end{figure}

We report in Figure~\ref{fig:realisticalongline} the pressure solution along two different lines: $\linea$, $y = 500$ [m], and $\lineb$, $x = 625$ [m], for the meshes with $\text{be}=28$.
We note that except for the case with $\text{amr}=5$, the solutions behave similarly,
showing again that the behavior of the solutions is mostly dominated by the number of elements in the fracture.

\begin{figure}[hbt!]
\centering
\includegraphics[width=12cm]{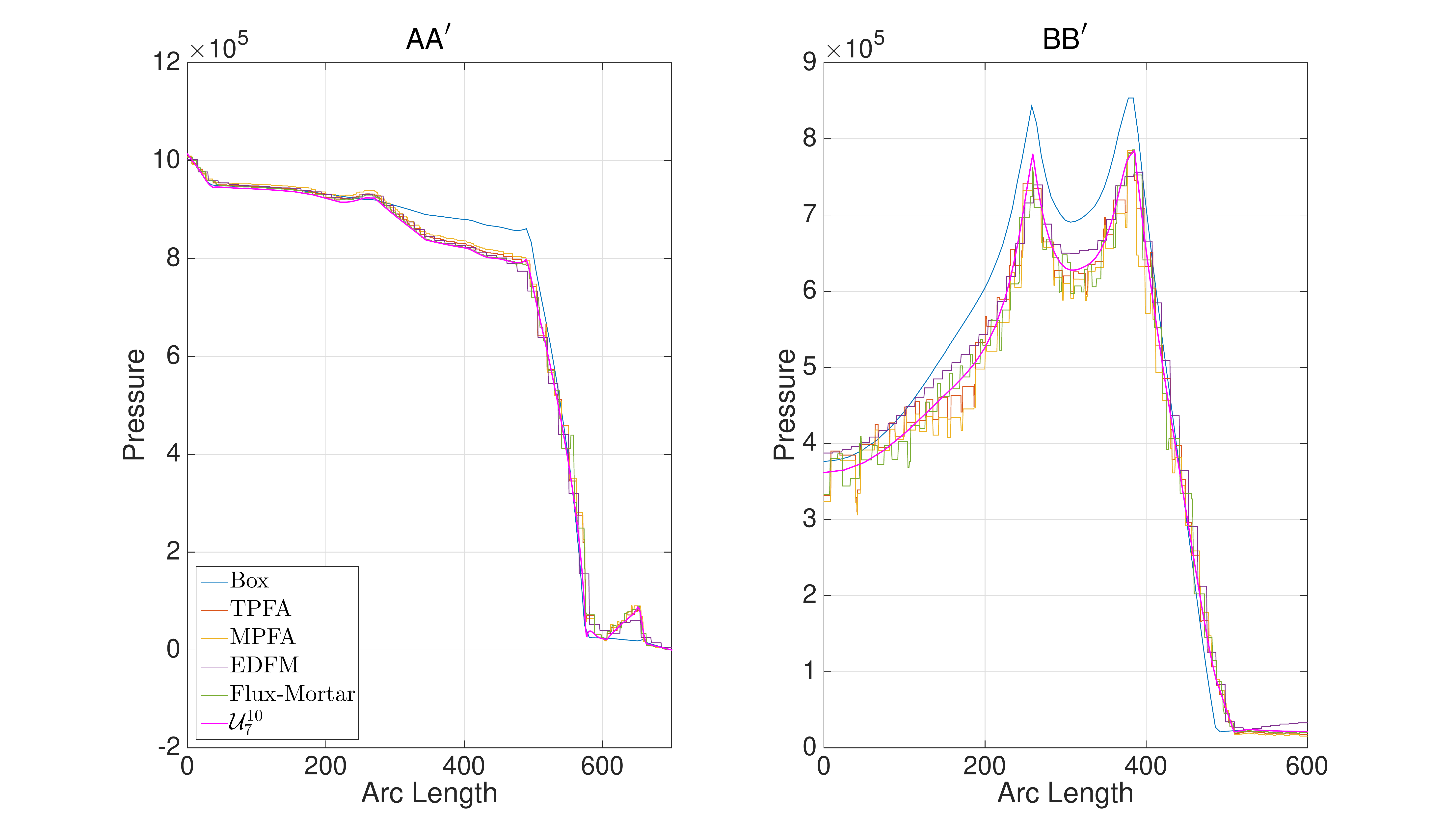}
\caption{Comparison of pressure values along two lines between the unresolved test case and the numerical solutions based on the hybrid-dimensional approaches proposed in~\cite{flemisch2018benchmarks}.}
\label{fig:rcomparisonwiththeir}
\end{figure}

Finally, in Figure~\ref{fig:rcomparisonwiththeir}, we report a comparison between our stabilized method and the methods based on the hybrid-dimensional representation of the domain
presented in~\cite{flemisch2018benchmarks}.
Our stabilized method based on AMR provides a possible alternative to hybrid methods for computing solutions with a two-dimensional representation of the fractures.
As we see, all the characteristics of the solutions along these two lines are correctly described.

\begin{table}
\caption{$L_2$ relative error,  Equation~\ref{eq:error}, for the boundary flux $q_{\lineb}(\domain_1)$ computed on the resolved meshes.}
\centering
\begin{tabular}{| c | c | c | c |}
\hline
&  \multicolumn{3}{ c | }{Non-stabilized}  \\
\hline
   \diagbox{amr}{be}  & 7 & 14 & 28 \\ 
\hline
\hline
5 & & & 0.31745 \\
6 & & 0.31745  & 0.024124  \\
7 & 0.31745  & 0.024133  & 0.020829  \\
8 & 0.024133  & 0.020821  & 0.0016197  \\
9 & 0.020821  & 0.0016153  &  \\
10 & 0.0016153  & & \\
\hline
&\multicolumn{3}{ c | }{Stabilized}\\
\hline
   \diagbox{amr}{be}  & 7 & 14 & 28 \\ 
\hline
\hline
5 & & &0.31604\\
6 & & 0.31635  & 0.021175 \\
7 & 0.31635  & 0.020986  & 0.027017 \\
8 & 0.020986  & 0.026737  & 0.0078999\\
9&0.026737  & 0.0074733  &\\
10& 0.0074733  & & \\
\hline
\end{tabular}
\label{tab:realisticerrors}
\end{table}

\newpage

\section{Conclusions}

We have presented a continuous finite element method for the numerical simulation of the coupled flow and transport problems in fractured porous media
based on an equi-dimensional representation of the fractures and an adaptive mesh refinement.
Adaptive mesh refinement allows to automatically and efficiently adapt an initially uniform mesh to any kind of heterogeneities.
The resulting non-conforming mesh does not resolve the matrix-fracture interfaces which can be approximated with the desired accuracy.

For the proposed finite element discretization on non-conforming meshes,
we have extended standard results on the discrete maximum principle, proposed a stabilization strategy based on algebraic diffusion operators, and shown that global and local conservation properties hold.
The stabilization strategy is based on the algebraic flux correction method which allows us to ensure the discrete maximum principle and to remove the unphysical oscillations in the solution.

We have presented a numerical comparison with a standard finite element method based conforming meshes which resolve the interfaces.
Such numerical examples highlight the main advantages in terms of accuracy and solution properties deriving from the use of meshes with non-skewed elements.
In particular, the use of adapted meshes allows to avoid large aspect ratios for the elements inside the fracture aperture, which represents the main drawback for equi-dimensional models.
On the other hand, the numerical results illustrate that the number of elements inside and around the fracture has a significant influence on the accuracy of the solution.
As the last contribution, we have shown that the proposed discretization method provides accurate results also for realistic cases with complex fracture networks.

%
%

{It is worth to point out that the proposed approach allows describing the realistic equidimensional  structure of  fracture networks without requiring their explicitly representation.

In the future we plan to employ such a discretization method  to validate the approaches based on a lower-dimensional representation of fractures,
which are usually employed to perform simulation in fractured media.}

The study of the accuracy of the proposed method for realistic three-dimensional cases should also be performed.
As the proposed method does not involve the coupling of different discretizations for the matrix and the background,
the linear systems to be solved are characterized by positive definite M-matrices and, hence, are suitable for the efficient solution of the problem with multigrid methods.
This is another fundamental characteristic of the proposed method.
In this context, the resulting hierarchy of meshes could also be considered in the construction of the levels of a multigrid solver.
Moreover, we marked for refinement the elements which presented a non-empty overlap with at least one fractures.
This idea comes on the \text{a-priori} information that the discretization error in the simulation of heterogenous media is usually localized close to the interfaces.
Obviously, the development of suitable error estimators would also provide a strategy to control such numerical errors.

\section*{Acknowledgments}
M. G. C. Nestola  acknowledges gratefully the support the PASC Project FASTER.
Marco Favino acknowledges gratefully the support of the Swiss National Science Foundation (SNSF) through the grant PZ00P2$\_$180112.
All methods and routines used in this study are implemented within the open-source software library \emph{Parrot}:~\url{https://github.com/favinom/parrot/}.
\emph{Parrot}'s lead developers are Maria GC Nestola (nestom@usi.ch) and Marco Favino (marco.favino@unil.ch).



\section{Appendix A. Mesh Characteristics}
Mesh characteristics of the resolved and the unresolved meshes employed in Sec,~\ref{sec:numerical_results}.

\begin{table}[ht!]
\caption{Resolved meshes for the benchmark  \lq\lq Single fracture network\rq\rq.}
\centering
\begin{tabular}{| l | r | r | r | }
 \hline
Mesh & $n^E$ & $ n^r$ & $w$\\
 \hline
   \hline
$\mathcal{R}^{114,2}$ &6\,612&6\,785 & 2\\ 
$\mathcal{R}^{228,4}$ &26\,448&27\,793 & 4  \\
$\mathcal{R}^{456,8}$ & 105\,792 & 106\,481 & 8\\
$\mathcal{R}^{921,16}$  & 423\,168 & 424\,545 &16\\
 \hline
\end{tabular} 
\label{tab:mesh_single_res}
\end{table}

\begin{table}
\caption{Unresolved meshes for the benchmark \lq\lq Regular fracture network\rq\rq.}
\centering
\begin{tabular}{| l | r | r | r | r |}
\hline
Mesh & ${E}_{tot}$ & $ {N}_{tot}$  & $ N_{tot}^r$  & $w$\\
\hline
\hline
$\mathcal U^{80}_{7}$ &219\,256&254\,851&184\,067&1.02\\
\hline
$\mathcal U^{80}_{8}$&434\,224&505\,641&363\,225&2.05\\
\hline
$\mathcal U^{80}_{9}$&1\,720\,768&1863685&1578389&4.10\\
\hline
$\mathcal U^{80}_{10}$&4\,300\,600&4586823&4014951&8.19\\
\hline
$\mathcal U^{80}_{11}$&12\,900\,400&13473253&12328181&16.38\\
\hline
$\mathcal U^{160}_{6}$&236\,848&272341&202069&1.02\\
\hline
$\mathcal U^{160}_{7}$&451\,816&523131&381227&2.05\\
\hline
$\mathcal U^{160}_{8}$&1\,736\,896&1879729&1594897&4.10\\
\hline
$\mathcal U^{160}_{9}$&431\,6728&4602867&4031459&8.19\\
\hline
$\mathcal U^{160}_{10}$&12\,916\,528&13489297&12344689&16.38\\
\hline
$\mathcal U^{320}_{5}$&310\,360&345631&276431&1.02\\
\hline
$\mathcal U^{320}_{6}$&525\,328&596421&455589&2.05\\
\hline
$\mathcal U^{320}_{7}$&1\,807\,264&1949893&1666085&4.10\\
\hline
$\mathcal U^{320}_{8}$&4\,387\,096&4673031&4102647&8.19\\
\hline
$\mathcal U^{320}_{9}$&12\,986\,896&13559461&12415877&16.38\\
\hline
$\mathcal U^{640}_{4}$&610\,912&645721&578713&1.02\\
\hline
$\mathcal U^{640}_{5}$&825\,880&896511&757871&2.05\\
\hline
$\mathcal U^{640}_{6}$&2101\,312&2243497&1961833&4.10\\
\hline
$\mathcal U^{640}_{7}$&4\,681\,144&4966635&4398395&8.19\\
\hline
$\mathcal U^{640}_{8}$&13\,280\,944&13853065&12711625&16.38\\
\hline
$\mathcal U^{1280}_{3}$&1\,826\,344&1860211&1797635&1.02\\
\hline
$\mathcal U^{1280}_{4}$&2\,041\,312&2111001&1976793&2.05\\
\hline
$\mathcal U^{1280}_{5}$&3\,303\,520&3444781&316750&4.10\\
\hline
$\mathcal U^{1280}_{6}$&5\,883\,352&6167919&5604063&8.19\\
\hline
$\mathcal U^{1280}_{7}$&14\,483\,152&15054349&13917293&16.38\\
\hline
$\mathcal U^{2560}_{2}$&6\,714\,736&6746701&6693037&1.02\\
\hline
$\mathcal U^{2560}_{3}$&\,\,\,6\,929\,704&6997491&6872195&2.05\\
\hline
$\mathcal U^{2560}_{4}$&\,\,\,8\,165\,248&8304625&8036209&4.10\\
\hline
$\mathcal U^{2560}_{5}$&10\,745\,080&11027763&10472771&8.19\\
\hline
$\mathcal U^{2560}_{6}$&19\,344\,880&19914193&18786001&16.38\\
\hline
\end{tabular}
\label{tab:meshcharacteristics1}
\end{table}

\begin{table}[ht]
\caption{Unresolved meshes for the benchmark  \lq\lq single fracture network\rq\rq.}
\centering
\begin{tabular}{| l | r | r | r | r | }
  \hline
Mesh & $n^E$ & $ n^H$  & $ n^r$ & $w$\\
 \hline
   \hline
 $\mathcal{U}^{100}_8$        &254\,740         &295\,757     &214\,173  &2.56 \\
 $\mathcal{U}^{200}_6$        &162\,628        &141\,640      & 121\,903  & 1.28  \\
 $\mathcal{U}^{200}_7$        &283\,720         &324\,775     &  243\,509 & 2.56 \\
 $\mathcal{U}^{200}_8$        &698\,200         &780\,233    & 617\,035    & 5.12   \\
 $\mathcal{U}^{400}_6$        &442\,813         &401\,680     &362\,185 & 2.56  \\
 $\mathcal{U}^{400}_7$        &815\,920         &898\,029      &735\,471& 5.12 \\
 $\mathcal{U}^{400}_8$       &2\,037\,040     &2\,201\,075   &1\,874\,681& 10.24   \\
 $\mathcal{U}^{800}_6$       &1\,291\,360     &1\,373\,625    &1\,212\,347 &5.12    \\
 $\mathcal{U}^{800}_7$       &2\,512\,480     &2\,678\,304    &2\,351\,557 &10.24    \\
  $\mathcal{U}^{800}_8$     &6\,890\,080      &7\,218\,127   & 6\,565\,333  &   20.48\\
  $\mathcal{U}^{1600}_6$   &4\,423\,360     &4\,587\,867   & 4\,265\,313& 10.24 \\
 $\mathcal{U}^{1600}_7$    &8\,800\,960      & 9\,129\,323 &8\,479\,089 & 20.48\\
  \hline
\end{tabular}
\label{tab:mesh_single_unres}
\end{table}

\begin{table}
\centering
\caption{Unresolved meshes for the benchmark   \lq\lq realistic fracture networks\rq\rq.}
\begin{tabular}{| l | r | r | r | r | r |}
\hline
Mesh & $n^E$ & $ n^H$  & $ n^r$ & $w$ \\
\hline
\hline
$\mathcal U^{7}_{7}$&62232&75439&49251 &0.01\\
\hline
$\mathcal U^{7}_{8}$&133128&162364&104142&0.03\\
\hline
$\mathcal U^{7}_{9}$&278223&340375&216354&0.06\\
\hline
$\mathcal U^{7}_{10}$&573276&701910&444958&0.12\\
\hline
$\mathcal U^{14}_{6}$&62232&75439&49251&0.01\\
\hline
$\mathcal U^{14}_{7}$&133128&162364&104142&0.03\\
\hline
$\mathcal U^{14}_{8}$&278223&340375&216354&0.06\\
\hline
$\mathcal U^{14}_{9}$&573276&701910&444958&0.12\\
\hline
$\mathcal U^{28}_{5}$&62286&75492&49321&0.01\\
\hline
$\mathcal U^{28}_{6}$&133182&162417&104212&0.03\\
\hline
$\mathcal U^{28}_{7}$&278277&340428&216424&0.06\\
\hline
$\mathcal U^{28}_{8}$&573330&701963&445028&0.12\\
\hline
\end{tabular}
\label{tab:meshcharacteristicsrealsitc}
\end{table}

\bibliographystyle{plain}
     \bibliography{reference.bib}
\end{document}